\documentclass[leqno,french,italian,german,english,amsthm,oneside]%
{elsart}
\usepackage{babel}
\usepackage[mathscr]{eucal}
\usepackage{amscd,amsfonts,amssymb,amsmath,amsthm,backref,url,
enumerate,mdwlist}

\usepackage{graphicx}

\usepackage{natbib}
\makeatletter
\DeclareRobustCommand\citeauthorpar
   {\begingroup\NAT@swatrue\def\NAT@ctype{1}\NAT@partrue\NAT@citetp}
\makeatother

\usepackage[pagebackref,
            citecolor=cyan,
            linkcolor=cyan,
            pagecolor=cyan]
{hyperref}

\renewcommand*{\backref}[1]{}
\renewcommand*{\backrefalt}[4]{%
\ifnum#1=1
Cited on p. #2.
\else
Cited on pp. #2.
\fi
\par}

\usepackage[all]{hypcap} 

\usepackage[oldenum]{paralist} 

\renewcommand\qedsymbol{\ensuremath\Box}
\newenvironment{unproved}[1]{\def\QQQQ{\string#1}\begin{\QQQQ}}%
                                 {\hfill\qedsymbol\end{\QQQQ}}

\input{fiddles.tex}

\makeatletter
\def\th@plain{%
    \thm@preskip\parskip
    \thm@postskip\z@
    \itshape 
}
\def\th@definition{%
    \thm@preskip\parskip
    \thm@postskip\z@
    \normalfont 
}
\def\th@remark{%
    \thm@headfont{\itshape}%
    \normalfont 
    \thm@preskip\parskip \divide\thm@preskip\tw@
    \thm@postskip\z@
}
\def\swappedhead#1#2#3{%
  \thmnumber{#2}.%
  \thmname{\@ifnotempty{#2}{~}#1}%
  \thmnote{ {\the\thm@notefont(#3)}}}
\makeatother

\swapnumbers
\theoremstyle{plain}
\newtheorem{theorem}{Theorem}

\newtheorem{corollary}[theorem]{Corollary}
\newtheorem{lemma}[theorem]{Lemma}
\newtheorem{proposition}[theorem]{Proposition}
\newtheorem{motto}[theorem]{Motto}

\theoremstyle{definition} 
\newtheorem{definition}[theorem]{Definition}
\newtheorem{definitions}[theorem]{Definitions}

\newtheorem{construction}[theorem]{Construction}

\newtheorem{question}[theorem]{Question}

\newtheorem{questions}[theorem]{Questions}

\newtheorem{example}[theorem]{Example}
\newtheorem{examples}[theorem]{Examples}

\newtheorem*{theorem*}{Theorem}

\newtheorem*{congruence theorem}{Congruence Theorem}

\newtheorem*{corollary to MFW}{Corollary to the %
Morton--Franks-Williams Theorem}

\input{mdefs} 
\begin{document}
\let\NoHyper\relax

\begin{frontmatter}
\title{Knot theory of complex plane curves}
\author{Lee Rudolph\thanksref{support}}
\thanks[support]{During the preparation 
and completion
of this survey,
the author was partially supported
by the Fonds National Suisse 
and by a National Science Foundation
Interdisciplinary Grant in the Mathematical Sciences
(DMS-0308894)
.}
\address{Department of Mathematics \& Computer Science
and Department of Psychology,\\
Clark University, Worcester MA 01610 USA}  
\ead{\href{mailto:lrudolph@black.clarku.edu}{lrudolph@black.clarku.edu}}

\begin{abstract}
%
%
The primary objects of study 
in the ``knot theory of complex plane curves'' are 
\moredef{{\Clink}links}:
links (or knots) cut out of a $3$-sphere in $\C^2$ by complex plane 
curves.  There are two very different classes of {\Clink}links, 
\moredef{transverse} and \moredef{totally tangential}.
Transverse {\Clink}links are naturally oriented.  
There are many natural classes of examples: 
links of singularities; links at infinity;
links of divides, free divides, tree divides, and graph divides;
and---most generally---quasipositive links.
Totally tangential {\Clink}links are unoriented but naturally 
framed; they turn out to be precisely the real-analytic 
Legendrian links, and can profitably be investigated in terms
of certain closely associated transverse {\Clink}links.

The knot theory of complex plane curves
is attractive not only for its own internal
results, but also for its intriguing relationships and
interesting contributions elsewhere in mathematics.  
Within low-dimensional topology, related subjects include 
braids, 
concordance, 
polynomial invariants, 
contact geometry, 
fibered links and open books, and
Lefschetz pencils. 
Within low-dimensional algebraic 
and analytic geometry, related subjects include 
embeddings and injections of the complex line in the complex plane,
line arrangements,
Stein surfaces,
and 
Hilbert's $16$th problem.  
There is even some
experimental evidence that 
nature favors quasipositive knots. %

%
%

\end{abstract}

\begin{keyword}
Arrangement \sep 
braid \sep 
braid monodromy \sep 
braided surface \sep 
Chisini's statement \sep
{\Clink}link \sep 
Hilbert's $16$th problem \sep
Jacobian conjecture \sep 
knot \sep 
labyrinth \sep 
link \sep 
link at infinity \sep 
link of indeterminacy \sep 
Milnor map \sep 
Milnor's question \sep 
monodromy \sep 
plane curve \sep 
quasipositivity \sep
singularity\sep
Thom conjecture\sep
Thurston--Bennequin invariant\sep
unknotting number\sep
Zariski conjecture

\def\MSC{\par\leavevmode\hbox {\it 2000 MSC:\ }}
\MSC Primary 57M25 \sep  
Secondary 14B05 \sep  
20F36 \sep  
32S22 \sep  
32S55 \sep  
51M99
\end{keyword}

\end{frontmatter}

\newpage

%
%
\setcounter{page}{2}\def\thepage{\roman{page}}

\section*{Foreword}
In the past two decades, knot theory in general 
has seen much progress and many changes.
\speak{Classical knot theory}---the study of knots 
as objects in their own right---has taken great strides,
documented throughout this Handbook 
(see the contributions by 
\citeauthor{BirmanBrendle2004},
\citeauthor{Hoste2004},
\citeauthor{Kauffman2004},
\citeauthor{Livingston2004},
and \citeauthor{Scharlemann2004}).
Simultaneously, there have been extraordinarily 
wide and deep developments in what might be called 
\speak{modern knot theory}: the study 
of knots and links in the presence of extra structure,\footnote{%
  Some observers have also detected \speak{postmodern knot theory}:
  the study of extra structure in the absence of knots.} %
for instance, a hyperbolic metric on the knot complement 
(as in the articles by \citeauthor{Adams2004}
and \citeauthor{Weeks2004}) or
a contact structure on the knot's ambient $3$-sphere
(as in the article by \citeauthor{Etnyre2004}). 
In these terms, the knot theory of complex plane curves
is solidly part of modern knot theory---the knots and links
in question are {\Clink}links, and the extra structures 
variously algebraic, analytic, and geometric.

``Some knot theory of complex plane curves'' \citep{Rudolph1983d} 
was a broad view of the state of the art in 1982.  Here I propose 
to look at the subject through a narrower lens, that of 
quasipositivity.
\Secref{intro} is devoted to general notations and definitions.
\Secref{braids and braided surfaces} is a treatment of 
braids and braided surfaces tailored to quasipositivity.\footnote{%
   Consult \citeauthor{BirmanBrendle2004} for a deeper and broader
   account of braid theory.} General transverse {\Clink}links 
are constructed and described in \Secref{transver}, while 
\Secref{alggeom} is a brief look at the special transverse 
{\Clink}links that arise from complex algebraic geometry 
``in the small'' and ``in the large''---to wit, links 
of singularities and links at infinity.\footnote{%
   One consequence of this survey's bias towards quasipositivity 
   is a de-emphasis of other aspects of the knot theory of 
   links of singularities and links at infinity; the reader is 
   referred to \citet{BoileauFourrier1998} (who include sections 
   on both these topics), to the discussions of singularities 
   and their higher-dimensional analogues by \citet[\S2]{Durfee1999} 
   and \citet[\S1]{Neumann2001}, and of course to the extensive
   literature on both subjects---particularly links of 
   singularities---referenced in those articles.}  Totally 
tangential {\Clink}links are constructed and described in \Secref{totaltan}.
The material in \S\Secref{transver}\thru\ref{totaltan} is 
related to other research areas in \Secref{applic}.
In \Secref{future} I give some fairly explicit, somewhat 
programmatic suggestions of directions for future research
in the knot theory of complex plane curves. 

Original texts of some motivating problems in the 
knot theory of complex plane curves are collected in 
an \hyperlink{appendix}{Appendix}.  
This survey concludes
with an \hyperlink{index of definitions}{index of definitions and notations}
and a \hyperlink{references}{bibliography}.
\hyperlink{open questions}{Open questions} are distributed throughout.

%
%

\tableofcontents
\newpage
\setcounter{page}{1}\def\thepage{\arabic{page}}
%
%
\section{Preliminaries}\label{intro}

Terms being defined are set in {\bydeffont{this}} typeface;
mere emphasis is indicated {\shoutfont{thus}}. 
A definition labelled as such is either of greater 
(local) significance, or non-standard to an extent which might
lead to confusion; labelled or not, potentially startling
definitions and notations are flagged with the symbol 
{\raise1ex\hbox{\tinydbend}} in the margins of both
the text and the \hyperlink{index of definitions}{index}.  The end 
or omission 
of a proof is signalled by $\qedsymbol$. 
Both $A\isdefinedas B$ and $B\defines A$ mean ``$A$ is defined as $B$''.
The symbol \danger${\iso}$%
\index{$natural isomorphism$@$\iso$ (natural isomorphism)|indanger} is
reserved for a \moredef{natural isomorphism}%
    \index{natural!isomorphism(s)}%
    \index{isomorphism(s), natural} (in an appropriate category).%
    \index{$non-standard$@ \mbox{{\raise1.25ex\hbox{\dbend}} (non-standard definition)}}%

\subsection{Sets and groups}\label{sets}
The set of real \resp{complex} numbers is $\R$ \resp{$\C$};
write $z\mapsto\conjug{z}$ for \moredef{complex conjugation}%
    \index{conjugate!complex ($\conjug{\phold}$)}%
    \index{$conjugate$@$\conjug{\phantom.\phold\phantom.}$ (complex 
conjugation)} $\C\to\C$, and $\Re$%
    \index{Re@$\Re$ (real part)} \resp{$\Im$}%
    \index{Im@$\Im$ (imaginary part)} for \moredef{real part} 
\resp{\moredef{imaginary part}} $\C\to\R$.
For $x\in\R$, let $\Rge{x}\isdefinedas\{t\in\R\Suchthat t\ge x\}$,
$\Rg{x}\isdefinedas \{t\in\R\Suchthat t>x\}$,
$\Rle{x}\isdefinedas \{t\in\R\Suchthat t\le x\}$,
$\Rl{x}\isdefinedas \{t\in\R\Suchthat t<x\}$.
Let $\EHP\isdefinedas \{w\in\C:\pm\Im{w}\ge 0\}$%
     \index{C lower half plane@$\LHP$ (lower half-plane in $\C$)}%
     \index{C upper half plane@$\UHP$ (upper half-plane in $\C$)}.  Let 
$\N\isdefinedas \Z\cap\Rnonneg $%
     \index{N1@$\N$ (non-negative integers)} and 
$\Nstrict\isdefinedas \N\cap\Rpos$%
     \index{N2@$\Nstrict$ (positive integers)}.  For 
$n\in\N$, let $\intsto{n}\isdefinedas \{k\in\Nstrict\Suchthat k\le n\}$%
     \index{N0@$\intsto{n}$ (integers from $1$ to $n$)}.
Denote projection on the $i$th factor of a cartesian product 
by $\pr_i$.  
The Euclidian norm on $\R^n$ or $\C^n$ is 
$\norm{\phold}$\index{$euclidian$@$\norm{\phold}$ (Euclidian norm)}.
For $\mathbf{u}, \mathbf{v}$ in a real vectorspace,
$\clint{\mathbf{u}}{\mathbf{v}}$ is 
$\{(1-t)\mathbf{u}+t\mathbf{v}\Suchthat 0\le t\le 1\}$.

Let $X$ be a set.  Denote the identity map $X\to X$ by 
$\id{X}$\index{id@$\id{X}$ (identity map)}, and 
the cardinality of $X$ by 
$\card{X}$\index{card()@$\card{\phold}$ (cardinality)}.
A \bydef{characteristic function}%
    \index{function!characteristic} on $X$ 
is an element of $\{0,1\}^X$; for $Y\sub X$, let
$\charfun{Y\into X}\from X\to\{0,1\}$
denote \moredef{the characteristic function of $Y$ in $X$}%
     \index{characteristic function!of $Y$ in $X$}, so
$\charfun{Y\into X}(x)=1$ for $x\in Y$, 
$\charfun{Y\into X}(x)=0$ for $x\in\diff{X}{Y}$.
A \bydef{multicharacteristic function}%
    \index{function!multicharacteristic} on $X$ 
is an element of $\N^X$.  Let $\charfun{}\in\N^X$.
The \bydef{total multiplicity}%
    \index{multicharacteristic function!total multiplicity}%
    \index{multiplicity!total (of a multicharacteristic function)} of 
$\charfun{}$ is $\sum_{x\in X}\charfun{}(x)\in\N\cup\{\infty\}$.
For $n\in\N$, an $n$\nobreakdash-\bydef{multisubset} of $X$ is a pair
$(\supp{\charfun{}},\charfun{}\restr{\supp{\charfun{}}})$ where 
$\charfun{}$ is a multicharacteristic 
function on $X$ of total multiplicity $n$.  
Call the set of all $n$\nobreakdash-multisubsets of $X$ 
the $n$th \moredef{multipower set}%
     \index{multipower!set}%
     \index{MPnX@\protect$\multsets{\phold}n$ ($n$th multipower)!set} of
$X$, denoted $\multsets{X}n$. 
(Note that if $X\ne X'$ then no multicharacteristic 
function on $X$ is a multicharacteristic function on $X'$, whereas 
$\multsets{X}{n}\cap \multsets{X'}{n}\ne\emptyset$ if and only
if $X\cap X'\ne\emptyset$.)  
For any $A$, any $f\from A\to \multsets{X}n$
can be construed as a \moredef{multivalued}%
    \index{function!multivalued}%
     \index{multivalued function} (specifically, 
an \moredef{$n$\nobreakdash-valued})%
     \index{nvalued function@$n$-valued function} \moredef{function
from $A$ to $X$}; typically $\gr{f}%
\isdefinedas \{(a,x)\in A\times X\Suchthat x\in f(a)\}$,
the \moredef{multigraph of $f$}%
    \index{multigraph ($\gr{\phold}$)}, determines 
neither $n$ nor $f$ (unless $n$ is---\shout{and is 
known to be}---equal to $1$), but the notation is still useful.
The \moredef{type}%
    \index{type!of a multicharacteristic function} of $\charfun{}\in\N^X$ 
is $\typeof{\charfun{}}\isdefinedas%
\operatorname{card}\after(\charfun{}^{-1}\restr\Nstrict)\in\N^{\Nstrict}$
(a multicharacteristic function on $\N$).
Identify an $n$\nobreakdash-\moredef{subset}%
     \index{nsubset@$n$-subset} of $X$ (i.e., $Y\sub X$ with $\card{Y}=n$) 
with the $n$\nobreakdash-multisubset $(Y,1)$ of 
type $n\charfun{\{1\}\into\N}$, and the set of \nsubset{n}s of $X$
with the $n$th \moredef{configuration set}%
\index{configuration!set}%
\index{$config$n@\protect$\binom{\phold}{\phold}$ (configuration set or space)}%
\[
\smash{
\subsets{n}{X}\isdefinedas%
\{(\supp{\charfun{}},\charfun{}\restr{\supp{\charfun{}}})\in\multsets{X}n 
\Suchthat \typeof{\charfun{}}=n\charfun{\{1\}\into\N}\}
}
\]
of $X$.
Call $\disc{n}{X}\isdefinedas\diff{\multsets{X}n}{\subsets{n}{X}}$ the 
$n$th \moredef{discriminant set}%
   \index{discriminant!set}%
   \index{DeltanX@\protect$\disc{n}{\phold}$ ($n$th discriminant)!set} of $X$.

Let $G$ be a group.  For $g,h\in G$, let 
$\conjugate{g}{h}$ \resp{$\commute{g}{h}$;$\YBax{g}{h}$} denote
the \moredef{conjugate} %
\resp{\moredef{commutator}; \moredef{yangbaxter}}%
   \index{conjugate!in a group ($\conjugate{g}{h}$)}%
   \index{commutator ($\commute{g}{h}$)}%
    \index{yangbaxter ($\YBax{g}{h}$)} $ghg^{-1}$ %
\resp{$ghg^{-1}h^{-1}=\conjugate{g}{h} h^{-1}$; %
      $ghgh^{-1}g^{-1}h^{-1}=\conjugate{gh}g{h^{-1}}$}.
For $A\sub G$, let $\genby{A}{G}$%
   \index{$group$@$\genby{\phold}{\phold}$ (group generated by)} be 
the subgroup generated by $A$, i.e.,
$\bigcap \{H\Suchthat\text{$A\sub H$ and $H$ is a subgroup of $G$}\}$,
and let $\genbyelt{a}{G}\isdefinedas\genby{\{a\}}{G}$ 
(when $G$ is understood, it may be dropped from these notations).
The \moredef{normal closure}%
    \index{normal!closure} of $A$ in $G$ is 
$\genby{\{\conjugate{g}{a}\Suchthat g\in G, a\in A\}}{G}$.
A \bydef{presentation} of $G$, denoted 
\begin{equation}\label{presentation of a group}
G=\gp{\freegen{i} \,(i\in I)}{\relator{j} \,(j\in J)},
\end{equation}
consists of:  
\begin{inparaenum}
\item
a short exact sequence 
$R\subset F\overset{\pi}{\surjection} G$
in which $F$ is a free group and $R$ is a subgroup of $F$;
\item
a family $\{\freegen{i}\Suchthat i\in I\}\sub F$ of 
free generators of $F$, the 
\moredef{generators of} \dispeqref{presentation of a group}%
  \index{generator(s)!of a presentation}%
  \index{generator(s)!of a group}%
  \index{presentation!generators of}%
  \index{group!generators of}; and
\item
a family $\{\relator{j}\Suchthat j\in J\}\sub R$
with normal closure $R$, the 
\moredef{relators of} \dispeqref{presentation of a group}%
  \index{relator(s)!of a presentation}%
  \index{relator(s)!of a group}%
  \index{presentation!relator(s) of}%
  \index{group!relator(s) of}. 
\end{inparaenum}
(Sometimes the elements $\pi(\freegen{i})$ of $G$ are also, abusively, 
called \moredef{the generators of $G$} with respect to 
\dispeqref{presentation of a group}.) %
The presentation \dispeqref{presentation of a group} is %
\moredef{Wirtinger}%
    \index{Wirtinger presentation}%
    \index{presentation!Wirtinger} in case 
every relator $\relator{j}$ is of the form 
$\conjugate{w(j)}{\freegen{s(j)}}\mskip1mu{\freegen{t(j)}}^{-1}$ 
for some $s, t\from J\to I$ and $w\from J\to F$.
Denote the 
\moredef{free product}%
    \index{free!product ($\freeprod$)}%
    \index{$free$@$\freeprod$ (free product)} of 
groups $G_0$ and $G_1$ by $G_0\freeprod G_1$.

A \bydef{partition} of a set $X$ is the quotient set 
$X/{\equiv}$ of $X$ by an equivalence relation $\equiv$ 
on $X$.  Call $X/{\approx}$ a \bydef{refinement} of
$X/{\equiv}$ in case each $\approx$-class
is a union of $\equiv$-classes.
Given $f\from X\to Y$, write 
$X/f$ for $X/{\equiv}_f$, where $x_0\equiv_f x_1$
iff $f(x_0)=f(x_1)$.  Given a (right) action of a group $G$
on $X$, write $X/G$ for $X/{\equiv_G}$, where $x_0\equiv_G x_1$
iff $x_1=x_0g$ for some $g\in G$; as usual, $xG$ stands for the 
$\equiv_G$-class (i.e., the $G$-orbit) of $x$.
The group $\symgroup{n}$%
    \index{Sn@$\symgroup{n}$ (symmetric group on $n$ letters)} of 
bijections $\intsto{n}\to\intsto{n}$ acts in a standard way on 
$X^n=X^{\intsto{n}}=\{f \Suchthat f\from\intsto{n}\to X\}$.
An \bydef{unordered \protect\ntuple{n}} in $X$ is an element 
$U{\ontuple{x}{n}}\isdefinedas\ontuple{x}{n}\symgroup{n}$ 
of the $n$th \bydef{symmetric power} $\sympower{X}n$ of $X$, 
where $U=\Unorder{X}{n}\from X^n\to\sympower{X}n$ is the 
\moredef{unordering map}%
    \index{unordering map ($U$)}%
    \index{U (unordering map)@$U$ (unordering map)}.
The map
\begin{equation}\label{bijection from sympower to multsets}
\begin{split}
\sympower{X}n & \to\multsets{X}{n}
\suchthat
\smash[t]{%
U(\overbrace{x_1,\dots,x_1}^{n_1},
\overbrace{x_2,\dots,x_2}^{n_2},
\dots,
\overbrace{x_k,\dots,x_k}^{n_k})
}
\\
& \mapsto 
(\{x_1,x_2,\dots,x_k\},(\{x_1,x_2,\dots,x_k\}\to\Nstrict
\suchthat x_i\mapsto n_i))
\end{split}
\end{equation}
(where $x_i\ne x_j$ for $i\ne j$, $n_i>0$, and $n=\sum_{i=1}^k n_i$)
is a bijection. 
\shout{Except} in \dispeqref{bijection from sympower to multsets},
notation is abused in the standard way, so that 
$\untuple{x}{n}$ denotes $U{}{}\ontuple{x}{n}$.  
In case $X$ is ordered by $\preccurlyeq$, 
notations like $\utupsep{x}{1}{k}{\preccurlyeq}{\dotsm}\in\subsets{k}{X}$,
$\{i\preccurlyeq j\}\from C\to\subsets{2}{X}$, and so on, mean 
``$\{x_1,\dots\,x_k\}\in\subsets{k}{X}$ and 
$x_1\preccurlyeq\dotsm\preccurlyeq x_k$'', 
``$\{i(c),j(c)\}\in\subsets{2}{X}$ and $i(c)\preccurlyeq j(c)$ 
for all $c\in C$'', and so on. 
In case $X$ is totally ordered by $\prec$,
call $\smash[t]{\{s\prec t\}, \{s'\prec t'\} \in \subsets{2}{X}}$ 
\bydef{linked} 
\resp{\bydef{unlinked}; \moredef{in} \bydef{touch} \moredef{at $u$}}
iff either $s\prec s'\prec t\prec t'$ or $s'\prec s\prec t'\prec t$
\resp{either $s\prec s'\prec t'\prec t$ or $s'\prec s\prec t\prec t'$;
$\{u\}=\{s,t\}\cap \{s',t'\}$}.

The bijection \dispeqref{bijection from sympower to multsets}
induces a definition of $\typeof{\untuple{x}{n}}$%
    \index{type!of an unordered \protect\ntuple{n}} as 
itself an unordered \ntuple{n}, such that, e.g.,
$\typeof{\{1,1,1,1\}}=\{4\}$, 
$\typeof{\{1,1,2,3\}}=\{1,1,2\}$, etc.

\subsection{Spaces}\label{spaces}
A simplicial complex is not necessarily finite.
A geometric realization of a simplicial complex $\Kscr$ is denoted 
$\geomrel{\Kscr}$\index{$geometric$@$\geomrel{\phold}$ (geometric realization)}. %
A \bydef{triangulation} of a topological space $X$ is a homeomorphism
$\geomrel{\Kscr}\to X$ for some $\Kscr$.
A \bydef{polyhedron} is a topological space $X$ 
which is the target of some triangulation.
Let $X$ be a polyhedron.  The set of components of $X$ 
is denoted $\compsof{X}$\index{pi0@$\compsof{\phold}$ (components)}.  
The Euler characteristic of $X$ is denoted 
$\euler{X}$\index{$euler$@$\chi$ (Euler characteristic)}.
The fundamental group of $X$ with \moredef{base point}%
    \index{base!point ($\basept$)} %
    \index{$base point$@$\basept$ (base-point)} $\basept$ 
is denoted by $\fundgp{X}{\basept}$%
     \index{pi1@$\fundgpnbp{\phold}$ (fundamental group)}%
     \index{pi1@$\fundgpnbp{\phold}$ (fundamental group)!$\fundgp{\phold}{\basept}$ (--- with base point $\basept$)}, or
simply $\fundgpnbp{X}$ in contexts where $\basept$ can be 
safely suppressed.  Call $\fundgpnbp{\diff{X}{K}}$ the 
\moredef{knotgroup}%
    \index{knotgroup (of a subspace)} of $K$ in case $X$ is connected
and $K\sub X$.

Manifolds are smooth ($\mathscr C^\infty$) unless otherwise stated.
A manifold $M$ may have boundary, but corners (possibly cuspidal)
only when so noted.  Denote by 
$\Bd M$ \resp{$\Int{M}$; $\tanbun M$} 
the interior \resp{boundary; tangent bundle}%
     \index{interior!of a manifold ($\Int{}$)}%
     \index{boundary!of a manifold ($\Bd$)}%
     \index{$boundary$@$\Bd$ (boundary of a manifold)}%
     \index{tangent bundle of a manifold ($\tanbun{\phold}$)} of $M$.  
Call $M$ \moredef{closed} \resp{\moredef{open}}%
     \index{closed!manifold}%
    \index{manifold!closed}%
    \index{open!manifold}%
    \index{manifold!open} in case $M$ is compact \resp{$M$ has no
compact component} and $\Bd M=\emptyset$.
Manifolds are (not only orientable, but) 
oriented, unless otherwise stated:
in particular, a complex manifold (e.g., $\C^n$ or $\P_n(\C)$)
has a \moredef{natural orientation}%
    \index{natural!orientation(s)}%
    \index{orientation(s)!natural}, and $\R$, $D^{2n}\isdefinedas %
\{\ontuple{z}{n}\in\C^n \Suchthat \norm{\ontuple{z}{n}}\le 1\}$,
and $S^{2n-1}\isdefinedas \Bd D^{2n}$ have 
\moredef{standard orientations}%
    \index{standard!orientation(s)}%
    \index{orientation(s)!standard}, as do $S^2$ (identified with 
$\Cext\isdefinedas \C\cup\{\infty\}$), 
$\R^3$ (identified with $\diff{S^3}{(0,1)}$) 
and the \bydef{bidisk} $D^2\times D^2$ (with corners 
$S^1\times S^2$).
The tangent space $\tanspace{x}{M}$ to $M$ at $x\in M$ is
an oriented vectorspace.  Let $-M$ \resp{$+M$; $\unorient{M}$}
denote $M$ with orientation reversed \resp{preserved; forgotten};
in case $\mathsf{M}\from M\to M$ is a diffeomorphism reversing 
orientation, $\Mir Q\isdefinedas \mathsf{M}(Q)$ is a 
\moredef{mirror image}%
    \index{mirror image ($\Mir\phold$)} of $Q\sub M$.

Call $Q\sub M$ \moredef{interior}%
     \index{interior!submanifold} \resp{\moredef{boundary}%
     \index{boundary!submanifold}} in case 
$Q\sub \Int{M}$ \resp{$Q\sub \Bd M$}.  In case $Q$ has 
a closed regular neighborhood $\Nb{Q}{M}$%
    \index{Nb@$\Nb{Q}{M}$ (regular neighborhood of $Q$ in $M$)}%
    \index{regular!neighborhood ($\Nb{\phold}{\phold}$)} in 
$(M,\Bd M)$, the \moredef{exterior}%
    \index{exterior ($\Ext{\phold}{\phold}$)} of $Q$ 
in $(M,\Bd M)$ is $\Ext{Q}{M}\isdefinedas\diff{M}{\Int{\Nb{Q}{M}}}$.

A \moredef{connected sum}%
     \index{connected sum ($\connsum$)}%
     \index{sum!connected ($\connsum$)}%
     \index{$connected$@$\connsum$ (connected sum)} \resp{%
\moredef{boundary-connected sum}}%
     \index{boundary!-connected sum ($\bdconnsum$)}%
     \index{sum!boundary-connected ($\bdconnsum$)}%
     \index{$boundary-connected$@$\bdconnsum$ (boundary-connected sum)} of 
\nmanifold{n}s $M_1$ and $M_2$ is denoted 
$M_1\connsum M_2$ \resp{$M_1\bdconnsum M_2$};
notations like $M_1\connsum_{D^n} M_2$, 
$M_1\bdconnsum_{D^{n-1}} M_2$, and so on, can 
be used for greater precision.

\begin{definitions}
A \bydefdanger{stratification} of a topological space $X$ is a  
locally finite partition $X/{\equiv}$ such that
\begin{inparaenum}
\item
each $\equiv$\nobreakdash-equivalence class, 
equipped with the topology induced
from $X$, is a connected, not necessarily oriented, manifold 
(called a $\equiv$\nobreakdash-\bydef{stratum} of $X$, or simply a 
\moredef{stratum} of $X$ when $\equiv$ is understood), and
\item
for every stratum $M$, the closure of $M$ in $X$
is a union of strata.  
\end{inparaenum}
A \bydef{vertex} of $X/{\equiv}$ is 
a point $x$ such that $\{x\}$ is a stratum; 
let $\Verts{X/{\equiv}}$\index{V1@$\Verts{\phold}$ (vertices)} denote 
the set of vertices of $X/{\equiv}$.
An \bydefdanger{edge} of $X/{\equiv}$ is the closure of a 
stratum of dimension $1$; let 
$\Edges{X/{\equiv}}$%
   \index{E1@$\Edges{\phold}$ (edges)} denote
the set of edges of $X/{\equiv}$. 
A \bydef{cellulation} is a stratification such that each 
stratum is diffeomorphic to some $\R^k$.  
The cellulation \moredef{associated to a triangulation}%
    \index{cellulation!associated to a triangulation}
$h\from\geomrel{\Kscr}\to P$ is that with strata 
the $h$\nobreakdash-images of open simplices of $\geomrel{\Kscr}$;
a fixed or understood triangulation $h$ of a polyhedron $P$
determines edges and vertices of $P$ and gives sense to the
notations $\Verts{P}$ and $\Edges{P}$.
\end{definitions}

An \moredef{arc}%
     \index{arc(s)} is a manifold diffeomorphic to $\clint{0}{1}$;
in particular, in a real vectorspace if 
$\clint{\mathbf{u}}{\mathbf{v}}\ne\{\mathbf{u}\}$ then 
$\clint{\mathbf{u}}{\mathbf{v}}$ is an arc oriented 
from $\mathbf{u}$ to $\mathbf{v}$.
An \bydef{edge} is an unoriented arc (this is mildly inconsistent
with the definition above).  
A \moredef{simple \textup(collection of\textup) 
closed curve\textup(s\textup)}%
     \index{simple!closed curve(s)}%
     \index{curve(s)!simple closed}%
     \index{curve(s)!simple closed!collection of}%
     \index{closed!curve(s)!simple}%
     \index{simple!closed curve(s)!collection of} is 
a manifold diffeomorphic to (a disjoint union of copies of) $S^1$.
A \bydef{graph} is a polyhedron $G$ of dimension $\le 1$ 
equipped with a cellulation $G/{\equiv}$ (which need not be 
associated to a triangulation of $G$).  For every
graph $G$, there exists $\valfun{G}\from\Verts{G}\to\N$ such
that, for every triangulation of $G$, 
$\valence{x}{G}=\card{\{\eedge\in\Edges{G}\Suchthat x\in\eedge\}}$;
$\valence{x}{G}$ is the \moredef{valence of $x\in\Verts{G}$ in $G$}%
   \index{valence (of a vertex in a graph, $\valence{\phold}{\phold}$)}.
Call $x\in\Verts{G}$ an \moredef{endpoint}%
   \index{endpoint (of a graph)} \resp{%
\moredef{isolated point}%
   \index{isolated point (of a graph)}; 
\moredef{intrinsic vertex}%
   \index{intrinsic vertex (of a graph)}} of $G$ 
in case $\valence{x}{G}=1$ \resp{$\valence{x}{G}=0$; 
$\valence{x}{G}>2$}.
Call $x\in G$ an \moredef{ordinary point}%
   \index{ordinary point (of a graph)} in case
either $x\in\Verts{G}$ and $\valence{x}{G}=2$ or
$x\notin\Verts{G}$.
A graph embedded in $\C$ is \moredef{planar}%
     \index{planar graph}.  A \bydef{tree} 
is a finite, connected, acyclic graph.
Let $n\in\Nstrict$.  An \nstar{n}%
    \index{star} is a tree with $n+1$ vertices of 
which at least $n$ are endpoints\label{abstract n-star};
an \moredef{\ngon{n}}%
     \index{ngon@$n$\nobreakdash-gon}%
     \index{gon@-gon} is a 
\ndisk{2} $P$ equipped with a cellulation 
having exactly $n$ edges, all in $\Bd P$.

\begin{definitions}\label{definition of surface}
A \bydefdanger{surface} is 
a compact \nmanifold{2} no component of which has empty 
boundary.  
The \moredef{genus}%
   \index{genus!of a surface}%
   \index{surface!genus} of
a connected surface $S$ is denoted $\genus{}{S}$%
     \index{g@$\genus{}{\phold}$ (genus)!of a surface}. A 
surface is \moredef{annular}%
     \index{surface!annular}%
     \index{annular surface} in case
each component is an annulus.
A subset $X$ of a surface $S$ is \bydef{full} provided that no component 
of $\diff{S}{X}$ is contractible.
\end{definitions}

\begin{figure}
\centering
\includegraphics[width=.3\textwidth]{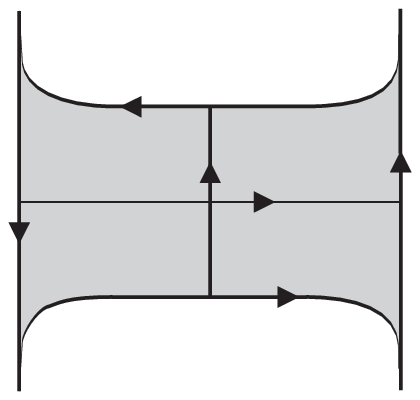}
\caption{The standard \protect\nhandle{1} and associated arcs.
\label{1-handle figure}}
\end{figure}

The \moredef{standard} (\ndimensional{2}) \moredef{\nhandle{0}}%
     \index{standard!\nhandle{0} ($\h0{}$)}%
     \index{h(0)@$\h{0}{}$ (\nhandle{0})!standard} is
$\h0{}\isdefinedas D^2$.  
Fix some continuous function $H\from\clint{-2}{2}\to\clint{1}{2}$ 
such that:
\begin{inparaenum}
\item
$H$ is even;
\item
$H(x)=1$ for $|x|\le 1$, $H(2)=2$, and $H\restr{\clint{1}{2}}$ is strictly
increasing;
\item
$H\restr{\opint{-2}{2}}$ and 
$(H\restr{\opclint{1}{2}})^{-1}$ are smooth;
\item
for $n\in\Nstrict$, $D^n((H\restr{\opclint{1}{2}})^{-1})(2)=0$. 
\end{inparaenum}
The \moredef{standard} (\ndimensional{2}) \moredef{\nhandle{1}}%
     \index{standard!\nhandle{1} ($\h1{}$)}%
     \index{h(1)@$\h{1}{}$ (\nhandle{1})!standard} is
$\h1{}\isdefinedas \{z\in\C\Suchthat 
|\Im{z}|\le H(\Re{z}),\medspace |\Re{z}|\le 2\}$,
a \nmanifold{2} with cuspidal corners.
The \moredef{attaching}%
     \index{attaching!arc(s)!of $\h1{}$}%
     \index{arc(s)!attaching!of $\h1{}$} \resp{\moredef{free}%
     \index{free!arc(s), of $\h1{}$}%
     \index{arc(s)!free, of $\h1{}$}} \moredef{arcs} of $\h1{}$ 
are the intervals $\clint{\pm(2-2\imunit)}{\pm(2+2\imunit)}\sub\C$
\resp{the arcs $\mp\{z\in\C\Suchthat \Im{z}=\pm H(\Re{z}),
|\Re{z}|\linebreak[1]\le 2\}$}.
The union of the attaching arcs of $\h1{}$ is
the \moredef{attaching region} $\attach{\h1{}}$%
     \index{attaching!region ($\attach{\phold}$)!of $\h1{}$}%
     \index{A()@$\attach{\phold}$ (attaching region)!of $\h1{}$} of $\h1{}$.  
The \moredef{standard core} \resp{\moredef{transverse}} \moredef{arc}%
     \index{standard!core arc!of $\h1{}$}%
     \index{core arc(s)!of $\h1{}$!standard}%
     \index{standard!transverse arc of $\h1{}$}%
     \index{transverse!arc(s) of $\h1{}$!standard!} of $\h1{}$ 
is $\core{\h1{}}\isdefinedas\clint{-2}{2}$ %
\resp{$\transv{\h1{}}\isdefinedas\clint{-\imunit}{\imunit}$}; 
a \moredef{core arc}%
    \index{core arc(s)!of $\h1{}$}%
    \index{arc(s)!core!of $\h1{}$} \resp{\moredef{transverse arc}}%
   \index{arc(s)!transverse!of $\h1{}$}%
   \index{transverse!arc(s) of $\h1{}$} of $\h1{}$ 
is any arc isotopic in $\h1{}$ to $\core{\h1{}}$ \resp{$\transv{\h1{}}$}.
(See \Figref{1-handle figure}.)

\begin{definitions}\label{standard bowtie}
Fix some smooth function $g\from\R\to\clint{0}{\pi}$ 
such that:
\begin{inparaenum}
\item
$g$ is odd, and periodic with period $8$;
\item
$g'(x)>0$ for $x\in\opint01$, $g(x)=\pi$ for $x\in\clint{1}{3}$,
$g'(x)<0$ for $x\in\opint34$, and $g(4)=0$.
\end{inparaenum}
The \moredef{standard} \bydefdanger{bowtie}%
     \index{standard!bowtie ($\bow{}$)}%
     \index{$standard bowtie$@ $\bow{}$ (standard bowtie)}%
     \index{bowtie!standard ($\bow{}$)} is 
$\bow{}\isdefinedas\Theta(\h1{})$, where
$\Theta\from\C\to\C\suchthat z\mapsto 
\Re{z}+\imunit\cos(\pi+g(\Re{z}))\Im{x}$.
Give each of
the two halves $\bow{}\cap\imunit\EHP$ of $\bow{}$ the orientation
induced on it from $\h1{0}$ by $\Theta$.  
The \moredef{attaching region}%
     \index{attaching!region ($\attach{\phold}$!of $\bow{}$}%
     \index{A()@$\attach{\phold}$ (attaching region)!of $\bow{}$}%
     \index{arc(s)!attaching!of $\bow{}$} of $\bow{}$ is
$\attach{\bow{}}\isdefinedas \Theta(\attach{\h1{}})$;
the \moredef{attaching arcs}%
     \index{attaching!arc(s)!of $\bow{}$} of 
$\bow{}$ are the components of $\attach{\bow{}}$.
The \moredef{standard core arc}%
     \index{standard!core arc!of $\bow{}$}%
     \index{core arc(s)!of $\bow{}$ (standard)}
     \index{arc(s)!core!of $\bow{}$} of $\bow{}$ is 
$\Theta(\core{\h1{}})$.
The \moredef{crossed arcs}%
     \index{crossed arc(s)!of $\bow{}$}%
     \index{arc(s)!crossed!of $\bow{}$} of $\bow{}$ are 
the $\Theta$-images of the free arcs of $\h1{0}$,
and the \moredef{crossing}%
    \index{crossing!of a bowtie} of $\bow{}$
is their point of intersection.
(See \Figref{bowtie figure}.)

\begin{figure}
\centering
\includegraphics[width=.3\textwidth]{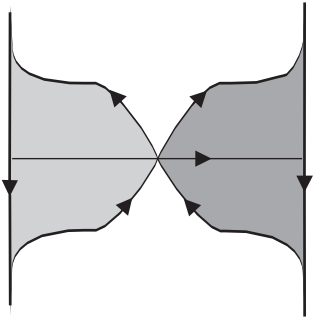}
\caption{The standard bowtie and associated arcs.
\label{bowtie figure}}
\end{figure}

\end{definitions}

Let $S$ be a \nmanifold{2}.  Call 
\begin{equation}\label{handle decomposition of a surface}
\smash[t]{S=\bigcup_{x\in X}\h{0}{x} \cup \bigcup_{t\in T}\h{1}{t}}
\end{equation}
a \moredef{\protect\nhandle{(0,1)} decomposition}%
     \index{handle decomposition} of $S$ provided that:
\begin{inparaenum}[(1)]
\item\label{hd:finite}
$X$ and $T$ are finite;
\item\label{hd:0-handles}
each 
\moredef{\nhandle{0}}%
     \index{handle}%
     \index{h(0)@$\h{0}{}$ (\nhandle{0})} %
$\h0x$ is diffeomorphic to the standard \nhandle{0} $D^2$;
\item\label{hd:1-handles}
each \moredef{\nhandle{1}}%
     \index{handle}%
     \index{h(1)@$\h{1}{}$ (\nhandle{1})} %
$\h1t$ is diffeomorphic, as a manifold with corners, 
to the standard \nhandle{1} $\h1{}$,
and thereby equipped with an attaching region
$\attach{\h1t}$, attaching arcs, free arcs,
a standard core arc $\core{\h1t}$ and other core arcs, 
a standard transverse arc $\transv{\h1t}$ and other transverse arcs;
\item\label{hd:disjoint}
the \nhandle{0}s are pairwise disjoint, as are the \nhandle{1}s;
\item
$\h{1}{t} \cap \bigcup_{x\in X}\h{0}{x} =
\h{1}{t} \cap \Bd \bigcup_{x\in X}\h{0}{x}=\attach{\h1t}$
for each $t\in T$; and
\item
the orientations of $S$ and all the $0$\nhandledash~and \nhandle{1}s are
compatible.
\end{inparaenum}
A \nmanifold{2} $S$ has a \nhandle{(0,1)} decomposition if and only if 
$S$ is a surface.

\subsection{Smooth maps}\label{smooth maps}
Maps between manifolds are smooth%
    \index{map!smooth} ($\mathscr C^\infty$), 
and isotopies are ambient, unless otherwise stated.
Given manifolds $M$ and $N$, let $\dfuns{M}{N}$%
    \index{diff@$\dfuns{M}{N}$ (smooth maps $M\to N$)} denote
the set of maps $M\to N$, with a suitable topology.
Let $f\in\dfuns{M}{N}$.  
The \moredef{$f$\nobreakdash-multiplicity}%
     \index{multiplicity!of a point with respect to a map}%
     \index{f-multiplicity@$f$-multiplicity} %
$x\in M$ is $\card{f^{-1}(x)}$.  A point of
$f$\nobreakdash-multiplicity $1$ \resp{$2$; at least $2$} is a 
\moredef{simple point}%
     \index{simple!point(s)}
\resp{\moredef{double point}%
     \index{double!point(s)}; 
\moredef{multiple point}%
     \index{multiple!point(s)}} of 
$f$; the image by $f$ of a simple \resp{double; multiple} point
of $f$ is a \moredef{simple} \resp{\moredef{double}; \moredef{multiple}}
\moredef{value}%
     \index{simple!value(s)}%
     \index{double!value(s)}%
     \index{multiple!value(s)}%
     \index{value(s)!simple}%
     \index{value(s)!double}%
     \index{value(s)!multiple} of $f$.  Let 
$\simplepts{f}\isdefinedas
\{x\in M\Suchthat x \text{ is a simple point of }f\}$,
$\doublepts{f}\isdefinedas 
\{x\in M\Suchthat x \text{ is a double point of }f\}$,
$\multiplepts{f}\isdefinedas 
\{x\in M\Suchthat x \text{ is a multiple point of }f\}$.
Denote by $\deriv{f}{x}\from\tanspace{x}{M}\to\tanspace{f(x)}{N}$
the derivative of $f$ at $x\in M$, and 
by $\derbun{f}\from\tanbun{M}\to\tanbun{N}$ 
the map induced by $f$ on tangent bundles.
A \moredef{critical point}%
    \index{critical!point(s) ($\critpts{\phold}$)} of $f$ is any $x\in M$ 
with $\rank(\deriv{f}{x})<\dim\tanspace{f(x)}{N}$.
The set of critical points of $f$ is denoted 
$\critpts{f}$, so $\critvals{f}$\index{critical!value(s)}%
     \index{value(s)!critical} is
the set of \moredef{critical values} of $f$.
As usual, $y\in \diff{N}{\critvals{f}}$ is called 
a \moredef{regular value}%
    \index{value(s)!regular}%
    \index{regular!value(s)} of $f$ (even if $y\notin f(M)$).
For $\dim(N)=1$, the index of $f$ at $x\in\critpts{f}$ is
denoted $\ind{f}{x}$\index{ind()@$\ind{f}{x}$ (index of $f$ at $x$)}.
Call $f$ a \moredef{Morse function}%
     \index{Morse!function} \resp{\moredef{Morse map}}%
     \index{Morse!map}%
     \index{map!Morse} in case 
$N=\R$ \resp{$N=S^1$}, 
$f$ is constant on $\Bd M$,
and every $x\in\critpts{f}$ is non-degenerate; in case 
also $\ind{f}{x}<\dim M$ for all $x\in\critpts{f}$, 
call $f$ \moredef{topless}%
     \index{Morse!function!topless}%
     \index{Morse!map!topless}%
     \index{topless!Morse function}%
     \index{topless!Morse map}.
Call $f\in\dfuns{M}{N}$ an \bydef{immersion} in case 
$\deriv{f}{x}\from\tanspace{x}{M}\to\tanspace{f(x)}{N}$
is injective for every $x\in M$.  
An \bydef{embedding} is an immersion that is a homeomorphism
onto its image.  Write $f\from M\imto N$ \resp{$f\from M\embto N$}%
    \index{$immersion$@$\imto$ (immersion)}%
    \index{$embedding$@$\embto$ (embedding)} to indicate 
that $f$ is an immersion \resp{embedding}.  The 
\moredef{normal bundle}%
     \index{normal!bundle!of an immersion} of $f\from M\imto\N$ is 
denoted $\normbun{f}$; given a submanifold $M\sub N$, the
\moredef{normal bundle of $M$ in $N$}%
     \index{normal!bundle!of a submanifold} is
$\normbun{\inclusion{M}{N}}$, where $\inclusion{M}{N}$ denotes inclusion.%
    \index{i@$\inclusion{Y}{X}$ (inclusion map)} %

Here are various constructions with normal bundles, 
in the course of which assorted notations 
and definitions are established.

\begin{definitions}\label{normal bundle stuff}
Let $M$ be a manifold of dimension $m$.
\begin{asparaenum}[\bf \thedefinitions.1.]
\item\label{collaring}
Let $Q\sub M$ be a submanifold of dimension $m-1$ with trivial 
normal bundle.  
A \bydef{collaring} of $Q$ in $M$ is an orientation-preserving
embedding $\col{Q}{M}\from Q\times \clint{0}{1}\embto \Nb{Q}{M}$ 
with $\col{Q}{M}(q,0)=q$ for all $q$;
a \bydef{collar} of $Q$ in $M$ is the image 
$\collar{Q}{M}\isdefinedas \col{Q}{M}(Q\times \clint{0}{1})$
of a collaring.  Let 
$\cut{M}{Y}\isdefinedas\diff{M}{\col{Q}{M}(Q\times \opint{0}{1})}$ 
be called $M$ \moredef{cut}%
    \index{cut ($\cut{}{}$)}%
    \index{$cut$@$\cut{}{}$ (cut)} along $Y$.
The \moredef{push-off map}%
    \index{push-off!map} of $Q$ 
is $Q\embto \diff{M}{Q}\from q\mapsto \col{Q}{M}(q,1)$.
Call the image $\pushoff{Y}$ of $Y\sub Q$ by the 
push-off map the \bydef{push-off} of $Y$.  (Note
that $\pushoff{Q}\sub \Bd\collar{Q}{M}$ has the ``outward 
normal'' orientation, whereas $Q\into \Bd\collar{Q}{M}$ 
reverses orientation.) 
It is convenient to define various 
\moredef{standard}%
     \index{standard!collaring}%
     \index{collaring!standard} collarings, thus.
\begin{inparaenum}
\item 
Given $Q\sub S^n=\Bd D^{n+1}$, let
$\col{Q}{D^{n+1}}\suchthat (x,t)\mapsto \e(1-t)x$
for a suitable $\e\in\opint01$.
\item
Given a manifold $W$, $u\in\Rpos\cup\{\infty\}$, and
$Q\sub \Int{W}\times\{0\}\sub\Bd W\times\clopint{0}{\pm u}$, let
$\col{Q}{\Bd W\times{[}0,\pm u{[}}\suchthat (x,t)\mapsto (x,\pm\e t)$
for a suitable $\e\in\Rpos$.
\item
For a suitable $\e\from\opint{-2}{2}\to\opint{0}{1}$,
$(\pm2+\imunit y,t)\mapsto \pm2+\imunit\e(t)y$ is
a collaring $\col{\Int{\attach{\h1{}}}}{\h1{}}$.
(The cusps of $\h1{}$ prevent the existence of a 
collaring of $\attach{\h1{}}$ in $\h1{}$.)
\end{inparaenum}%

\item\label{meridional disks}
Let $Q$ be a manifold of dimension $q<m$,
$j\from Q\imto M$ an immersion, $B\isdefinedas j(Q)$.  
For $x\in\diff{B}{j\of{\Bd Q\cup\multiplepts{j}}}$,
a \moredef{meridional \ndisk{(m-q)}}%
    \index{meridional disk ($\merid{\phold}{\phold}$)}%
    \index{disk!meridional ($\merid{\phold}{\phold}$)}%
    \index{$meridional$@$\merid{\phold}{\phold}$ (meridional disk)} %
\moredef{of $B$ at $x$} is
the image $\merid{B}{x}$ of an embedding $f\from D^{m-q}\embto M$
such that $f^{-1}(B)=f^{-1}(x)=0$, $f$ is transverse to $j$, 
and $\merid{B}{x}$ intersects $B$ positively (with respect to 
given orientations of $M$ and $Q$).  In case $M$ is connected
and $q=m-2$, any element of $\fundgpnbp{\diff{M}{B}}$
represented by a loop freely homotopic to $\Bd\merid{B}{x}$
in $\diff{M}{B}$ is called a \moredef{meridian}%
    \index{meridian (in a knotgroup)} in that knotgroup.
\item\label{adapted maps to the circle}
Let $B\sub M$ be a submanifold of dimension $m-2$ with
trivial normal bundle.  Let $n\from B\times\C\to\normbun{B}$
be a fixed trivialization 
of $\normbun{B}$.  In the standard way, 
using (inexplicit) metrics, etc., identify 
$\Nb{B}{M}$ to $n(B\times D^2)$
in such a way that if $x\in B$, then the image of 
$D^2\to M\suchthat z\mapsto n(x,z)$ is a meridional 
\ndisk{2} $\merid{B}{x}$.
Call $f\from\diff{M}{B}\to S^1$ \bydef{weakly adapted}%
    \index{adapted!weakly} to $n$ in case 
for every $Q\in\compsof{B}$ there is an integer
$d(Q)$ such that,
if $d(Q)\ne 0$, then $f(n(x,z))=(\argument{z})^{d(Q)}$ 
for all $x\in Q$, $z\in \diff{D^2}{\{0\}}$.
Call $f\from\diff{M}{B}\to S^1$ \bydef{adapted} in 
case $f$ is weakly adapted and, in addition,
if $d(Q)=0$, then $f$ extends to 
$f_Q\in\dfuns{\diff{M}{\diff{B}{Q}}}{S^1}$,
$f_Q\restr Q\from Q\to S^1$ is an immersion,
and $f_Q\restr\merid{B}{x}$ is constant for each $x\in Q$.
\end{asparaenum}
\end{definitions}

\begin{definitions}\label{definition of proper}
Let $f\from M\to N$ be smooth; let $Q$ be a \ncodimension0 
submanifold of $\Bd M$.  
Call $f$ \moredef{proper along $Q$}%
     \index{proper!along a boundary submanifold}%
     \index{immersion!proper}%
     \index{immersion!proper!along boundary submanifold}, %
\moredef{relative to} $\col{\Bd N}{N}$ with 
$\collar{\Bd N}{N}\sub \Nb{\Bd N}{N}$ 
and $\col{Q}{M}$, provided that:
\begin{inparaenum}
\item
$f(Q)\sub\Bd N$ and $f(\diff{M}{Q})\sub\Int{N}$;
\item
$f^{-1}(\Nb{\Bd N}{N})\sub M$ is a submanifold, and  
$f\restr{f^{-1}(\Nb{\Bd N}{N})}$ is an embedding;
\item
$f\after\col{Q}{M}\from\collar{Q}{M}\embto N$ 
is an embedding into $\Nb{\Bd N}{N}$ and
$\pr_2\after(\col{\Bd N}{N})^{-1}\after f\after\col{Q}{M}=
\pr_2\from\collar{Q}{M}\to\clint{0}{1}$.
\end{inparaenum}
Call $f$ \moredef{proper along $Q$} in case there exist
collarings $\col{Q}{M}$ and $\col{\Bd N}{N}$ such that
$f$ is proper along $Q$ relative to $\col{Q}{M}$ and $\col{\Bd N}{N}$
(``along $Q$'' is dropped when $Q=\Bd M$).
A properly embedded submanifold is simply \moredef{proper}%
     \index{submanifold!proper}%
     \index{proper!submanifold}%
     \index{embedding!proper}.
\end{definitions}

Some special cases of low-dimensional immersions and embeddings 
of particular interest, and associated ancillary constructions, 
need extra terminology.

\begin{definitions}\label{low-dim immersions}
Let $M$ be a compact \nmanifold{m},
$N$ an \nmanifold{n}.  Let $f\from M\imto N$ 
be an immersion with $\multiplepts{f}=\doublepts{f}$.
\begin{asparaenum}[\bf \thedefinitions.1.]
\item\label{half-proper}
Let $m=1$, and suppose $M$ is an arc or an edge.
Call $f$ \moredef{half-proper}%
     \index{immersion!half-proper}%
     \index{half-proper!immersion (of an arc or edge)}%
     \index{half-proper!arc or edge} provided that it is
proper along a single component of $\Bd M$.
A \moredef{half-proper} arc or edge is the image of
a half-proper embedding.  An \nstar{n} 
embedded in $N$ is \moredef{proper}%
    \index{proper!star}%
    \index{star!proper} 
provided each of its edges is half-proper.
\item\label{normal immersions of curves}%
Let $m=1$, $n=2$.  Call $f$ \moredef{normal}%
     \index{normal!immersion!of $M^1$ in $N^2$} %
provided that $f$ is proper
(whence $\doublepts{f}\sub\Int{M}$),
$\doublepts{f}$ is finite, and 
if $f(x_1)=f(x_2)$ with $x_1\ne x_2$, then 
the lines $\deriv{f}{x_1}(\tanspace{x_1}{M})$ and 
$\deriv{f}{x_2}(\tanspace{x_2}{M})$ in $\tanspace{f(x_1)}{N}$
are transverse.
A \moredef{crossing}%
     \index{crossing!of a normal immersion $M^1\imto N^2$} of
a normal immersion $f\from M\imto N$ is a double value of $f$; 
the \moredef{crossing number}%
     \index{crossing!number ($\crossings{\phold}$)!of a normal immersion}%
     \index{normal!immersion!crossing number}%
     \index{immersion!normal (of $M^1$ in $N^2$)!crossing}%
     \index{immersion!normal (of $M^1$ in $N^2$)!crossing number} of $f$
is $\crossings{f}\isdefinedas\card{f(\doublepts{f})}$.
A \moredef{branch}%
     \index{branch!at a crossing of $M^1\imto N^2$} of $f$ 
at a crossing $y$ is the 
(germ at $y$ of the) $f$-image of either component of 
$\Nb{f^{-1}(y)}{\diff{M}{(\diff{\doublepts{f}}{f^{-1}(y)})}}$.  
The image of a normal immersion of $S^1$ \resp{%
an arc; 
a finite disjoint union of copies of $S^1$; 
a finite disjoint union of arcs} is a \moredef{normal closed curve}%
     \index{normal!closed curve(s)}%
     \index{curve(s)!normal!closed}%
     \index{closed!curve(s)!normal}
\resp{a \moredef{normal arc}; 
a \moredef{normal collection of closed curves};
a \moredef{normal collection of arcs}}%
     \index{normal!closed curve(s)!collection of}%
     \index{normal!arc(s)}%
     \index{normal!arc(s)!collection of}%
     \index{arc(s)!normal}%
     \index{arc(s)!normal!collection of}; the 
\moredef{crossing number}%
     \index{crossing!number ($\crossings{\phold}$)!of normal closed curve(s)}%
     \index{crossing!number ($\crossings{\phold}$)!of normal arc(s)} 
of the normal collection of closed curves or arcs $f(M)$ is 
$\crossings{f(M)}\isdefinedas\crossings{f}$. 
A \bydef{bowtie} in $N$ is the image of $\bow{}\sub\h1{}$ 
by an \moredef{embedding}%
     \index{embedding!of the standard bowtie $\bow{}$} of 
$\bow{}$ in $N$, that is, a map $\delta\from\bow{}\to N$ 
that extends to an embedding $\h1{}\embto N$.  
Let $C\isdefinedas f(M)$ be a normal collection of closed curves.
For each $y\in f(\doublepts{f})$,
let $\bow{y}=\delta_y(\bow{})$ be a bowtie with 
$\bow{y}\sub\Nb{\{y\}}{\diff{N}{\diff{f(\doublepts{f})}{\{y\}}}}$,
such that the \moredef{crossed arcs}%
     \index{crossed arc(s)!of a bowtie}%
     \index{arc(s)!crossed!of a bowtie} of $\bow{y}$ (i.e., 
the $\delta_y$-images of the crossed arcs of $\bow{}$)
are the branches of $f$ at $y$, correctly oriented. %
\begin{figure}
\centering
\includegraphics[width=.8\textwidth]{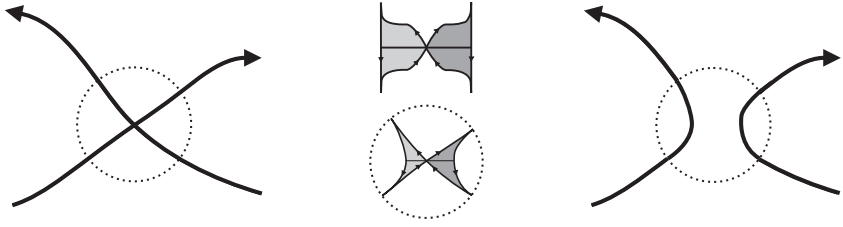}
\caption{Local smoothing, using a bowtie.
\label{smooloc figure}}
\end{figure}%
The \moredef{local smoothing}%
     \index{smoothing!of normal closed curve(s)!local}%
     \index{normal!closed curve(s)!local smoothing of} of $C$ 
at $y$ 
(see \Figref{smooloc figure})
is the normal collection of closed curves 
$\smoothloc{C}{y}\defines \smoothloc{f}{y}(\smoothloc{M}{y})$ 
created by replacing the crossed arcs of $\bow{y}$ 
with $-\delta_y(\attach{\bow{}})$. 
Here, $\smoothloc{M}{y}$ is unique up to
diffeomorphism, and $\smoothloc{C}{y}$ up to 
(arbitrarily small) isotopy;
and $\crossings{\smoothloc{C}{y}}=\crossings{C}-1$.  
The \moredef{smoothing}%
     \index{smoothing!of normal closed curve(s)}%
     \index{normal!closed curve(s)!smoothing of} of $C$
(see \Figref{smooth figure}) is 
the simple collection of closed curves $\smooth{C}\isdefinedas %
\smoothloc{\dotsm\smoothloc{\smoothloc{C}{y_1}}{y_2}\dotsm}{y_{\crossings{C}}}
\sub N$, independent of the enumeration $\{y_1,\dots,y_{\crossings{C}}\}$
of $f(\doublepts{f})$. 
\begin{figure}
\centering
\includegraphics[width=.6\textwidth]{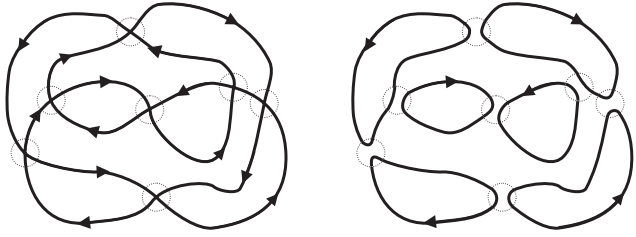}
\caption{A normal closed curve in $\R^2$ and its smoothing.
\label{smooth figure}}
\end{figure}%
\item\label{clasp immersions}
Let $m=2$, $n=3$.
Call $f$ \moredef{clasp}%
     \index{immersion!clasp (of $M^2$ in $N^3$)}%
     \index{clasp!immersion} provided that 
$\doublepts{f}$ is the union of finitely many
pairwise disjoint edges $A'_1,A''_1,\dots,A'_s,A''_s$ with
$f(A'_i)=f(A''_i)\sub\Int{N}$, both $A'_i$ and $A''_i$ 
half-proper ($i=1,\dots,s$); call $s$ the 
\moredef{clasp number}%
     \index{clasp!number!of a clasp immersion}%
     \index{clasp!immersion!clasp number}%
     \index{immersion!clasp (of $M^2$ in $N^3$)!clasp number} of the 
clasp immersion $f$, and denote it by $\clasp{f}$.  
The image $f(S)$ is a \moredef{clasp surface}%
     \index{clasp!surface}%
     \index{surface!clasp} (see \Figref{clasp figure});
the \moredef{clasp number}%
     \index{clasp!number!of a clasp surface} of 
$f(S)$ is $\clasp{f(S)}\isdefinedas\clasp{f}$.
\begin{figure}
\centering
\includegraphics[width=.6\textwidth]{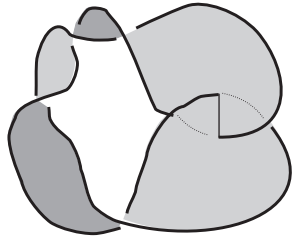}
\caption{A clasp surface in $\R^3$.
\label{clasp figure}}
\end{figure}
\item\label{ribbon immersions}
Let $m=2$, $n=3$.
Call $f$ \moredef{ribbon}    \index{immersion!ribbon (of $M^2$ in $N^3$)}%
    \index{ribbon!immersion} provided that 
$\doublepts{f}$ is the union of finitely many
pairwise disjoint edges $A'_1,A''_1,\dots,A'_s,A''_s$ with
$A'_i$ proper, $A''_i$ interior, and $f(A'_i)=f(A''_i)\sub \Int{N}$. The 
image of a ribbon immersion of a surface is a
\moredef{ribbon surface}\index{ribbon!surface!in $N^3$}%
 \index{surface!ribbon!in $N^3$} (see \Figref{ribboth figure}).%
\begin{figure}
\centering
\includegraphics[width=.6\textwidth]{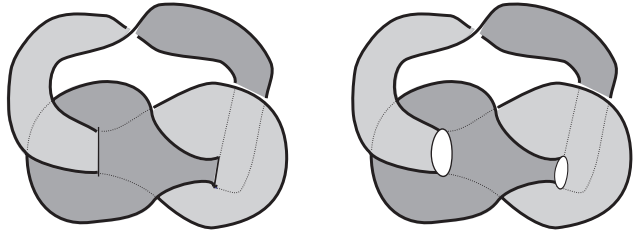}
\caption{A ribbon surface in $\R^3$ and its smoothing.
\label{ribboth figure}}
\end{figure}%
For connected $S$, the \moredef{genus}%
     \index{genus!of a ribbon surface} 
of the ribbon surface $f(S)$ is 
$\genus{}{f(S)}\isdefinedas\genus{}{S}$. The \moredef{smoothing}%
     \index{smoothing!of a ribbon surface}%
     \index{ribbon!surface!smoothing of} of a 
ribbon surface $R=f(S)$ is an embedded surface
$\smooth{R}\sub N$, unique up to isotopy, 
constructed as follows.  Let $S_0\isdefinedas %
\cut{S}{(\bigcup_{i=1}^s A'_i \cup \bigcup_{i=1}^s \Int{A''_i})}$.
Let ${\equiv}$ be the equivalence relation on $S_0$ having as its
non-trivial equivalence classes $\{x,y\}$ and 
$\{\pushoff{x},\pushoff{y}\}$ ($x\in \Int{A'_i}$, $f(x)=f(y)$) 
and $\{x,\pushoff{x},y\}$ ($x\in \Bd A'_i$, $f(x)=f(y)$), $i=1,\dots,s$. %
There is a natural way to smooth the quotient manifold  
$S_1\isdefinedas S_0/{\equiv}$, and a natural map 
$f_1\from S_1\to N$ with $f_1(S_1)=f(S)$, which after
an arbitrarily small perturbation yields an embedding 
$S_1\embto N$ with image $\smooth{f(S)}$.\footnote{%
     \citet{Fox1962} introduced the word ``ribbon'' into knot theory
     (specifically, in the context of ribbon-immersed \ndisk{2}s).  
     His usage was soon generalized \citep{Tristram1969} and widely 
     adopted.  In a distinct chain of development, the biologist 
     \citet{Crick1976}, followed by physicists 
     \citep*{Grundberg1989} and other scientists applying 
     mathematics, gave the word quite a different meaning
     (perhaps closer to its everyday use), essentially 
     to refer to twisted \ndimensional{2} bands.  
     More recently this conflicting usage has been adopted by some 
     knot-theorists \citep[see, e.g.,][]{ReshetikhinTuraev1990},
     particularly of a categorical bent.  \foreign{Caveat lector.}%
}%
\item\label{nodal immersions of surfaces}
Let $m=2$, $n=4$.  
Call $f$ \moredef{nodal}%
     \index{nodal!immersion!of $M^2$ in $N^4$} %
provided that $f$ is proper, 
$\doublepts{f}$ is finite,
and if $f(x_1)=f(x_2)$ and $x_1\ne x_2$, 
then the planes 
$\deriv{f}{x_1}(\tanspace{x_1}{M})$ and 
$\deriv{f}{x_2}(\tanspace{x_2}{M})$ in $\tanspace{f(x_1)}{N}$
are transverse.
A \bydef{node} of a nodal immersion $f$
is a double value of $f$; the \moredef{node number}%
    \index{number ($\node{\phold}$)!of a nodal immersion $M^2\imto N^4$}%
     \index{nodal!immersion!node number}%
     \index{immersion!nodal (of $M^2$ in $N^4$)!node number} of $f$
is $\node{f}\isdefinedas\card{f(\doublepts{f})}$. A \moredef{branch}%
     \index{branch!at a node of $M^2\imto N^4$} %
$f$ at a node $y$ is the (germ at $y$ of the) 
$f$-image of either component of 
$\Nb
    {
     f^{-1}(y)
    }
    {
     \diff{M}
          {
           \diff{\doublepts{f}
                }
                {f^{-1}(y)
                }
          }
    }$.  
The \moredef{sign}%
     \index{sign ($\sign{\phold}$)!of a node}%
     \index{node!sign of ($\sign{\phold}$)}%
     \index{positive!node}%
     \index{node!positive}%
     \index{negative!node}%
     \index{node!negative} %
$\sign{y}$ of the node $y$ is the sign ($+$ or $-$) 
of the given orientation of $\tanspace{y}{N}$ with respect to 
its orientation as the direct sum of the oriented $2$--planes 
$\deriv{f}{x_1}(\tanspace{x_1}{M})$ and $\deriv{f}{x_2}(\tanspace{x_2}{M})$ 
(in either order), where $f(x_1)=f(x_2)=y$, $x_1\ne x_2$.
The image of a nodal immersion of a surface is a
\moredef{nodal surface}%
    \index{nodal!surface}%
     \index{surface!nodal}.  
The \moredef{node number}%
    \index{node!number ($\node{\phold}$)!of a nodal surface} of 
a nodal surface $f(M)$ is $\node{f(M)}\isdefinedas\node{f}$,
and its \moredef{smoothing}%
     \index{smoothing!of a nodal surface}%
     \index{nodal!surface!smoothing of} is an embedded surface
$\smooth{f(M)}\sub N$, unique up to isotopy, constructed
by replacing $\Nb{f(\doublepts{f})}{f(M)}$ with 
$\card{f(\doublepts{f})}$ annuli embedded in 
$\Nb{f(\doublepts{f})}{N}$ in a standard way.
(In appropriate local coordinates on $N$, the neighborhood on $f(M)$
of a positive node is $\{(z,w)\in D^4\sub\C^2\Suchthat zw=0\}$, 
and is replaced by 
$\{(z,w)\in D^4\Suchthat zw=\e(1-\norm{(z,w)}^4)\}$, 
where $\e\from\clint01\to\Rnonneg$ is smooth, 
$0$ near $0$, positive near $1$, and sufficiently small.)%
\item\label{slice embeddings}
Let $m=2$, $n=4$.  Call $f$ %
\moredef{slice}%
     \index{immersion!slice (of $M^2$ in $N^4$)}%
     \index{slice!immersion!of $M^2$ in $N^4$} provided that $f$ 
is proper.  A \moredef{slice surface}%
     \index{slice!surface}%
     \index{surface!slice} is 
the image of a slice embedding of a surface (i.e.,
a proper surface).\footnote{%
   The redundant term ``slice surface'' has been 
   retained for the sake of tradition. 
   See \citet[\S1]{Rudolph1993} for a history of the use 
   of the word ``slice'' in knot theory.}%
\item\label{ribbon embeddings}
Let $m=2$, $n=4$, and suppose that $N$ is compact and 
$\rho\from N\to\R$ is a topless Morse function.  
Call $f$ \moredef{$\rho$-ribbon}%
     \index{immersion!rho-ribbon@$\rho$-ribbon (of $M^2$ in $N^4$)}%
     \index{rho-ribbon@$\rho$-ribbon immersion (of $M^2$ in $N^4$)} %
provided that $f$ is slice and 
$\rho\after f$ is a topless Morse function on $M$.
In case $N=D^4$ and $\card{\critpts{\rho}}=1$ (so that,
up to diffeomorphisms of $N$ and $\R$, $\rho=\norm{\phold}^2$),
a $\rho$-ribbon embedding is called simply a 
\moredef{ribbon embedding}%
     \index{embedding!ribbon}%
     \index{ribbon!embedding}, and 
its image is called a \moredef{ribbon surface}%
     \index{ribbon!surface!in $D^4$}%
     \index{surface!ribbon!in $D^4$} in $D^4$.
\end{asparaenum}
\end{definitions}

There is a close relation between ribbon surfaces in $D^4$ and in $S^3$.

\begin{unproved}{proposition}\label{ribbon=ribbon}
Let $M$ be a surface.  
If $f\from M\imto S^3=\Bd D^4$ is a ribbon immersion, then there 
is a non-ambient isotopy $\{f_t\from M\imto D^4\}_{t\in[0,1]}$
such that $f_0=f$, 
$f_t\restr{\Bd M}=f_0\restr{\Bd M}$ for $t\in\clint{0}{1}$,
and $f_t\from M\embto D^4$ is a ribbon embedding 
for $t\in\opclint{0}{1}$;
conversely, if $g\from M\embto \Bd D^4$ is ribbon,
then $g=f_1$ for some such non-ambient isotopy 
$\{f_t\from M\imto D^4\}_{t\in[0,1]}$ 
with $f_0\from M\imto S^3$ ribbon. 
\textup(Although the first of these non-ambient isotopies
is unique up to ambient isotopy, the second enjoys
no such uniqueness.\textup)
\end{unproved}

See \citet{Tristram1969}, \citet{Hass1983}, or \citet{Rudolph1985} 
for more detailed statements and proofs.

A (smooth) \moredef{covering map}%
    \index{covering!map}%
    \index{map!covering} is an orientation-preserving 
immersion $f\from M\imto N$ that is also a topological 
covering map (it is not required that the domain of a topological 
covering map 
be connected).
The usual theory of covering maps is assumed---in particular, 
the well-behaved, albeit many--many, correspondence for 
connected $N$ between permutation 
representations $\rho{}\from\fundgp{N}{\basept}\to\symgroup{s}$
and covering maps $f$ with target $N$ and \moredef{base fiber}%
    \index{base!fiber} $f^{-1}(\basept)=\intsto{s}$. %
Given $s\in\Nstrict$ and $\perm{g}\in\symgroup{s}$, 
let $\repn{\perm{g}}\from\fundgp{S^1}{1}\to\symgroup{s}$ be
the permutation representation $\circmap{k}\mapsto\perm{g}^k$,
where $\circmap{}$ denotes the \moredef{positive generator}%
    \index{positive!generator of $\fundgp{S^1}{1}$ ($\circmap{}$)}%
\index{$positive$@$\circmap{}$ (positive generator of $\fundgp{S^1}{1}$)} %
of $\fundgp{S^1}{1}\iso\fundgp{\diff{D^2}{\{0\}}}{1}$, 
that is, the homotopy class of $\id{S^1}$.
Call $\cover{\repn{\perm{g}}}\isdefinedas(\diff{D^2}{\{0\}})\times 
\intsto{s}/\genby{\perm{g}}{}$ the \moredef{standard covering 
space of $D^2$ of type $\perm{g}$}%
\index{standard!covering space of $D^2$ ($\cover{\repn{\perm{\phold}}}$)} and
$\covmap{\repn{\perm{g}}}\from\cover{\repn{\perm{g}}}\to D^2%
\suchthat (z,t\genby{\perm{g}}{})\mapsto 
z^{\card{t\genby{\perm{g}}{}}}$ the \moredef{standard covering 
map of $D^2$ of type $\perm{g}$}%
\index{standard!covering space of $D^2$ ($\cover{\repn{\perm{\phold}}}$)}%
\index{type!of a covering (space or map)}.

The theory of \moredef{branched covering maps}%
    \index{branched covering!map}%
    \index{covering!branched} originated in complex analysis and
algebraic geometry.  Following earlier work of 
\citet{Heegard1898} and \citet{Tietze1908}, the theory
was adapted to combinatorial manifolds by 
\citet{Alexander1920} and \citet{Reidemeister1926}, 
then to more general spaces by \citet{Fox1957}.
\citet{DurfeeKauffman1975} made the construction more 
precise in the smooth category.  The \foreign{ad hoc} 
approach of \eqref{branched covers} and 
\eqref{branched cover of 2-manifold}, below, 
suffices to handle the cases that are most 
important for the knot theory of complex plane curves.

\begin{definitions}\label{branched covers}
Let $s$, $\perm{g}$, etc., be as above.  The 
\moredef{standard branched covering space} \resp{\moredef{map}}
\moredef{of $D^2$ of type $\perm{g}$}%
\index{standard!branched covering!space of $D^2$ ($\brcover{\repn{\perm{\phold}}}$)}%
\index{standard!branched covering!map of $D^2$}%
\index{branched covering!space!standard, of $D^2$ ($\brcover{\repn{\perm{\phold}}}$)}%
\index{branched covering!map!standard, of $D^2$ ($\brcover{\repn{\perm{\phold}}}$)} %
is $\brcover{\repn{\perm{g}}}\isdefinedas%
{D^2}\times\intsto{s}/\genby{\perm{g}}{}$ \resp{%
$\brcovmap{\repn{\perm{g}}}\from%
\brcover{\repn{\perm{g}}}\to{D^2}\suchthat
(z,t\genby{\perm{g}}{})\mapsto 
z^{\card{t\genby{\perm{g}}{}}}$}.  
Let $N$ be a connected \nmanifold{n} equipped with
a stratification $N/{\equiv}$ such that:
\begin{inparaenum}%
\item
every stratum of $N/{\equiv}$ is smoothly immersed in
$N$;
\item
no stratum of $N/{\equiv}$ has dimension $n-1$ (whence
there is a unique stratum $N_0$ of dimension $n$) or $n-3$; and
\item
$\diff{N}{N_0}$ is the closure in $N$ of the union $B$ of all 
\ndimensional{(n-2)} strata.
A smooth map $f\from M\to N$ is a \moredef{branched 
covering map of $N$}%
     \index{branched covering!map} \moredef{branched over $B$}, and 
$M$ is a \moredef{branched covering space of $N$}%
     \index{branched covering!space} \moredef{branched over $B$}, 
provided that:
\item
$\critvals{f}\sub B$;
\item
$f\restr(f^{-1}(\diff{N}{B})$ is a covering map of degree $s$;
and
\item
for every $x\in B$, there exist $\perm{g}(x)\in\symgroup{s}$
(the \moredef{type of $f$ at $x$}%
    \index{type!of a branched covering (space or map)}), an 
embedding $\phi\from D^2\embto N$ onto a meridional 
disk $\merid{B}{x}$, and an embedding 
$\Phi\from\brcover{\repn{\perm{g}(x)}}\embto M$
with $(f\restr(f^{-1}(\merid{B}{x})))\after\Phi=
\phi\after\brcovmap{\repn{\perm{g}(x)}}$.  
\end{inparaenum}
The conjugacy class of $\perm{g}(x)$ in $\symgroup{s}$ is 
constant on each component of $B$, and 
$B=\{x\in\critvals{f}\Suchthat\perm{g}(x)\ne\id{\intsto{s}}\}$.
The branched covering $f$ is called \moredef{simple}%
    \index{branched covering!simple} in case $\perm{g}(x)$ is
a transposition for each $x\in B$.
\end{definitions}

\begin{construction}\label{branched cover of 2-manifold}
Let $N$ be a connected \nmanifold{2}, $B\sub\Int{N}$ a non-empty
finite subset.  Fix an enumeration $B\defines\{x_1,\dots,x_q\}$.  
Let $\basept_0\in\Ext{B}{\Int{N}}$.
The regular neighborhood $\Nb{B}{N}$ is a union of pairwise
disjoint meridional \ndisk{2}s $D^2_i\isdefinedas\merid{B}{x_i}$, 
$i\in\intsto{q}$.  Let $\basept_i\in\Bd D_i$.
Fix a proper \nstar{q} $\psi\sub\Ext{B}{N}$ with 
$\Verts{\psi}=\{\basept_0,\basept_1,\dots,\basept_q\}$
and $\Edges{\psi}\defines\{\eedge_1,\dots,\eedge_q\}$, such that
$\Bd{\eedge_k}=\{\basept_0,\basept_k\}$ and the cyclic order 
of $\Edges{\psi}$ at $\basept_0$ (with respect to the orientation
of $N$) is the cyclic order of their indices $1,\dots,q$.
Let $\Psi\isdefinedas\psi\cup\bigcup_{i=1}^q(\diff{D_i}{x_i})$.
Let $\freegen{i}\in \fundgp{\Psi}{\basept_0}$
be the element represented by a loop 
that traverses the edge $\eedge_{i}$ of $\psi$ 
from $\basept_0$ to $\basept_i$,
represents $\circmap{}$ in 
$\fundgp{\Bd D^2_i}{\basept_i}\iso\fundgp{S^1}{1}$,
and returns to $\basept_0$ along $\eedge_{i}$.  
Evidently $\fundgp{\Psi}{\basept_0}$ is the free group
$\gp{\freegen{i}, i\in\intsto{q}}{\emptyset}\iso
\fundgpnbp{\Bd D^2_1}\freeprod\dotsm\freeprod\fundgpnbp{\Bd D^2_q}$.
Call $\monogen{}\defines(\monogen{1},\dots,\monogen{q})\in\symgroup{s}^q$ 
\bydef{compatible} \moredef{with $\psi$} (or $\Psi$) 
in case there exists a permutation representation 
$\repn{\monogen{}}\from\fundgp{\Ext{B}{N}}{\basept_0}\to\symgroup{s}$ 
such that 
$\repn{\monogen{}}(\freegen{i}')=\monogen{i}$ ($i\in\intsto{q}$),
where $\freegen{}\mapsto\freegen{}'$ is the inclusion-induced
homomorphism $\fundgp{\Psi}{\basept_0}\to%
\fundgp{\diff{N}{B}}{\basept_0}\iso\fundgp{\Ext{B}{N}}{\basept_0}$.
(If $N$ is closed then compatibility is a genuine restriction;
if $N$ is not closed, then every $\monogen{}$ is compatible with 
$\psi$, but $\repn{\monogen{}}$ may not be unique if $N$ is not 
contractible.)  Finally, given $\Psi$, 
a compatible \ntuple{q} $\monogen{}$, and 
$\repn{\monogen{}}$, construct 
$\brcover{\repn{\monogen{}}}$ as the identification space
(with an appropriate, easily defined smooth structure)
of the disjoint union of copies of $\brcover{\repn{\monogen{i}}}$
($i\in\intsto{q}$) and 
$\cover{\repn{\monogen{}}}$ along their tautologously diffeomorphic
boundaries $\brcovmap{\repn{\monogen{i}}}^{-1}(\Bd D^2_i)=%
\Bd\brcover{\repn{\monogen{1}}}$ ($i\in\intsto{q}$) and 
$\covmap{\repn{\monogen{}}}^{-1}(\Bd D^2_i)\sub%
\Bd\cover{\repn{\monogen{}}}$; let 
$\brcovmap{\repn{\monogen{}}}\restr\brcover{\repn{\monogen{i}}}%
=\brcovmap{\repn{\monogen{i}}}$ and 
$\brcovmap{\repn{\monogen{}}}\restr\cover{\repn{\monogen{}}}%
=\covmap{\repn{\monogen{}}}$.
\end{construction}

\subsection{Knots, links, and Seifert surfaces}
\label{knots, links, and Seifert surfaces}

A \moredef{link}%
    \index{link(s)} is a simple collection of closed curves embedded
in $S^3$; except where otherwise stated, isotopic links are 
treated as identical.  A \ncomponent{1} link is a 
\moredef{knot}%
    \index{knot(s)}.  For $n\in\N$, $O^{(n)}$
    \index{unlink}
     \index{O(n)@$O^{(n)}$ (\protect\ncomponent{1} unlink)} denotes 
the \moredef{\ncomponent{n} unlink}, that is, the boundary of 
$n$ pairwise disjoint copies of $D^2$ embedded in $S^3$; the
\bydef{unknot} is $O\isdefinedas O^{(1)}$\index{O (unknot)@$O$ (unknot)}. 

\begin{definitions}\label{diagram definition}
A \bydefdanger{link diagram}%
    \index{D(L)@$\diagram{\phold}$ (link diagram)}%
    \index{diagram!link}%
    \index{diagram!link!information of ($\info{\phold}$)}%
    \index{diagram!link!picture of ($\pict{\phold}$)} is 
a pair $\diagram{L}\defines(\pict{L},\info{L})$%
    \index{I@information!about a link ($\info{\phold}$)}%
    \index{P@picture!of a link ($\pict{\phold}$)}, where
\begin{inparaenum}
\item
$L\sub\R^3\sub S^3$ is a link (called \moredef{the link of}%
     \index{link(s)!of a link diagram}%
     \index{link diagram!link of} the diagram);
\item\label{diagram: image}
$\pict{L}$, the $\diagram{L}$-\moredef{picture} of $L$,
is the image of $L$ by an affine projection 
$\proj{L}\from\R^3\to\C$; and 
\item\label{diagram: information}
$\info{L}$, the $\diagram{L}$-\moredef{information} 
about $L$, is information sufficient to reconstruct 
from $\pict{L}$ the embedding of $L$ into $\R^3$, 
up to isotopy respecting the fibers of $\proj{L}$.  
\end{inparaenum}
The \moredef{mirror image}%
     \index{mirror image ($\Mir\phold$)!of a link diagram} of
$\diagram{L}$ is the link diagram $\diagram{\Mir L}$,
where $\pict{\Mir L}\isdefinedas\conjug{\pict{L}}$, and
$\proj{L}$ and $\info{L}$ are modified accordingly to produce
$\proj{\Mir L}$ and $\info{\Mir L}$; of course the link of 
$\diagram{\Mir L}$ is the mirror image $\Mir L$ of $L$ as
already defined.
\end{definitions}

The nature of the information $\info{L}$ can be of various sorts, 
depending on context.  For instance, $\diagram{L}$ is 
a \moredef{standard link diagram}%
    \index{standard!link diagram}
    \index{link diagram!standard}%
    \index{diagram!link!standard} for $L$ provided that 
\begin{inparaenum}
\item
$\pict{L}$ is a normal collection of closed curves, and
\item
$\info{L}$ consists of \begin{inparaenum}
\item\label{standard diagram: global information}
the global information that $\proj{L}\restr{L}\from L\to\C$ 
is a normal immersion with image $\pict{L}$,
supplemented by 
\item\label{standard diagram: local information}
local information at each crossing specifying
which branch is ``under'' and which is ``over''---equivalently, 
which of the two points of $L\cap\proj{L}^{-1}(z)$
is the \moredef{undercrossing}%
    \index{undercrossing ($\underx{\phold}$)}%
    \index{$under$@$\underx{\phold}$ (undercrossing)} $\underx{z}$ 
(initial endpoint) and which is the \moredef{overcrossing}%
    \index{overcrossing ($\overx{\phold}$)}%
    \index{$over$@$\overx{\phold}$ (overcrossing)} $\overx{z}$ 
(terminal endpoint) of the interval between them on $\proj{L}^{-1}(z)$,
when the standard orientations of $\R^3$ and $\C$ are used
to orient the fibers of $\proj{L}$.  
(It is usual to depict crossings in the style of the left half of 
\Figref{positive crossing}.)
\begin{figure}
\centering
\includegraphics[width=.6\textwidth]{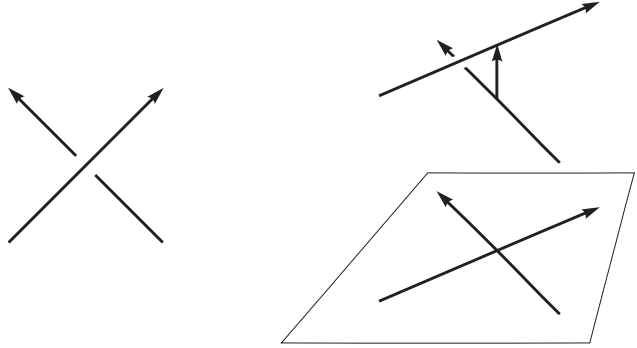}
\caption{A positive crossing in a standard link diagram.
\label{positive crossing}}
\end{figure}
\end{inparaenum} 
\end{inparaenum} 
Every link has many different standard link diagrams 
(some satisfying further conditions), 
as well as non-standard link diagrams of various types, some 
of which will be introduced later as needed.

\begin{definitions}\label{apparatus of standard diagram}
Let $\diagram{L}$ be a standard link diagram.
\begin{asparaenum}[\bf \thedefinitions.1.]
\item\label{Seifert cycle of standard diagram}
A \moredefdanger{Seifert cycle}\footnote{%
  The more commonly used term ``Seifert circle'' seems to have
  been popularized, if not coined, by Fox (\citeyear{Fox1962}; 
  see also a \citeyear{Fox1963} review in which Fox glosses 
  \citeauthor{Murasugi1960}'s ``standard loops'' as 
  ``Seifert circles'').	Certainly \citeauthor{Seifert1934}'s 
  term ``Kreis'' does mean ``circle'', but it can also be   
  translated as ``cycle'', and in the expositon of Seifert's 
  construction the latter translation has two apparent advantages over the 
  former: it does not connote geometric rigidity, and does connote 
  intrinsic orientation.}%
     \index{Seifert!cycle (of a standard link diagram)}%
     \index{link diagram!standard!Seifert cycle of ($\Seifo$)}%
     \index{link diagram!standard!set of Seifert cycles of ($\SeifO{}$)}%
     \index{O@$\SeifO{}$ (set of Seifert cycles)}%
     \index{o@$\Seifo$ (Seifert cycle)}%
     \index{cycle, Seifert|indanger}%
     \index{Seifert!circle|see {Seifert cycle}}%
     \index{circle, Seifert|see {Seifert cycle}} of $\diagram{L}$ is any 
$\Seifo\in\SeifO{\diagram{L}}\isdefinedas\compsof{\smooth{\pict{L}}}$.  The
\moredef{inside}%
     \index{inside, of a Seifert cycle ($\inside{\phold}$)}%
     \index{Seifert!cycle (of a standard link diagram)!inside of ($\inside{\phold}$)}%
     \index{cycle, Seifert!inside of ($\inside{\phold}$)} of 
a Seifert cycle $\Seifo$ is the \ndisk{2} $\inside{\Seifo}\sub\C$
oriented so that $\Bd{\inside{\Seifo}}=\Seifo$.  
The \moredef{sign}%
     \index{sign ($\sign{\phold}$)!of a Seifert cycle}%
     \index{Seifert!cycle (of a standard link diagram)!sign of}%
     \index{cycle, Seifert!sign of}%
     \index{positive!Seifert cycle}%
     \index{Seifert!cycle (of a standard link diagram)!positive}%
     \index{negative!Seifert cycle}%
     \index{Seifert!cycle (of a standard link diagram)!negative} $\sign{\Seifo}$ of 
$\Seifo$ is the sign, positive ($+$) or 
negative ($-$), of the orientation of $\inside{\Seifo}$
with respect to the standard orientation of $\C$.
Let $\SeifO{\diagram{L}}^{+}$ \resp{$\SeifO{\diagram{L}}^{-}$}
be the set of positive \resp{negative} Seifert cycles of  
$\diagram{L}$.%
\index{O@$\SeifO{}$ (set of Seifert cycles)!$\SeifO{}^{+}$ (--- positive ---)}%
\index{O@$\SeifO{}$ (set of Seifert cycles)!$\SeifO{}^{-}$ (--- negative ---)} %
Given $\Seifo,\Seifo'\in\SeifO{\diagram{L}}$, say
that $\Seifo$ \moredef{encloses}%
    \index{encloses ($\encloses$)}%
    \index{$encloses$@$\encloses$ (encloses)} $\Seifo'$,
and write $\Seifo\encloses\Seifo'$, 
in case $\Int{\inside{\Seifo}}\supset\Seifo'$.
Call $\diagram{L}$ \moredef{nested}\label{define:nested}%
    \index{nested standard link diagram}%
    \index{link diagram!standard!nested} \resp{\moredef{scattered}}%
    \index{scattered standard link diagram}%
    \index{link diagram!standard!scattered} in case
$\SeifO{\diagram{L}}$ is an 
${\encloses}$\nobreakdash-chain %
\resp{${\encloses}$\nobreakdash-antichain}. 
\item\label{crossings of standard diagram}
A \moredef{crossing} of $\diagram{L}$ is any 
$\Seifx\in\SeifX{\diagram{L}}\isdefinedas\proj{L}(\doublepts{\proj{L}})$%
     \index{crossing!of a standard link diagram}%
     \index{link diagram!standard!crossing of ($\Seifx$)}%
     \index{X@$\SeifX{}$ (set of crossings)}%
     \index{x@$\Seifx$ (crossing)}%
     \index{link diagram!standard!set of crossings of ($\SeifX{}$)}. 
Let $\Seifx\in\SeifX{\diagram{L}}$.
An \moredef{over-arc of $\diagram{L}$ on $L$ at $\Seifx$}%
\index{over-arc of a standard link diagram ($\overarc{\Seifx}{\phold}$)} is
any arc $\overarc{\Seifx}{L}\sub L$ such that 
\begin{inparaenum}
\item\label{over-arc at x is over x}
$\overx{\Seifx}\in\overarc{\Seifx}{L}$,
\item\label{over-arc at x isn't under anything}
$\overarc{\Seifx}{L}\cap \{\underx{\Seify}\Suchthat 
\Seify\in\SeifX{\diagram{L}}\}=\emptyset$, and 
\item
among arcs $\a\sub L$ satisfying (\ref{over-arc at x is over x}) 
and (\ref{over-arc at x isn't under anything}), 
$\overarc{\Seifx}{L}$ maximizes 
$\card{\a\cap\{\overx{\Seify}\Suchthat 
\Seify\in\SeifX{\diagram{L}}\}}$.
\end{inparaenum}
An \moredef{over-arc of $\diagram{L}$ in $\pict{L}$ at $\Seifx$} 
is any arc $\proj{L}(\overarc{\Seifx}{L})\sub\pict{L}$.
Let $\overx{\mathbf{u}}$ \resp{$\underx{\mathbf{u}}$} be 
a positively oriented basis vector for $\tanspace{\overx{x}}{L}$
\resp{$\tanspace{\underx{x}}{L}$}.
The \moredef{sign}%
     \index{sign ($\sign{\phold}$)!of a crossing}%
     \index{crossing!of a standard link diagram!sign of}%
     \index{positive!crossing (of a link diagram)}%
     \index{crossing!of a standard link diagram!positive}%
     \index{negative!crossing (of a link diagram)}%
     \index{crossing!of a standard link diagram!negative} $\sign{\Seifx}$ of 
$\Seifx$ is the sign, positive ($+$) or negative ($-$),
of the frame 
$(\underx{\mathbf{u}},\overx{\Seifx}-\underx{\Seifx},\overx{\mathbf{u}})$
with respect to the standard orientation of $\R^3$.
(The crossing in \Figref{positive crossing} is positive.)
Let $\SeifX{\diagram{L}}^{+}$ \resp{$\SeifX{\diagram{L}}^{-}$}
be the set of positive \resp{negative} crossings of 
$\diagram{L}$.%
\index{X@$\SeifX{}$ (set of crossings)!$\SeifX{}^{+}$ (--- positive ---)}%
\index{X@$\SeifX{}$ (set of crossings)!$\SeifX{}^{-}$ (--- negative ---)} 
\item\label{adjacency of crossings and cycles}
Call $\Seifo\in\SeifO{\diagram{L}}$ \bydef{adjacent} to 
$C\sub\SeifX{\diagram{L}}$ in case, for some $\Seifx\in C$,
$\Seifo$ contains an attaching arc of the bowtie 
$\bow{\Seifx}$ used to construct $\smooth{\pict{L}}$.
Let $\SeifO{\diagram{L}}^{\scriptscriptstyle\ge}$ 
\resp{$\SeifO{\diagram{L}}^{\scriptscriptstyle<}$;
$\SeifO{\diagram{L}}^{\scriptstyle\emptyset}$} denote the
set of Seifert cycles adjacent to $\SeifX{\diagram{L}}^{+}$
\resp{not adjacent to $\SeifX{\diagram{L}}^{+}$; 
not adjacent $\SeifX{\diagram{L}}$}.%
\index{O@$\SeifO{}$ (set of Seifert cycles)!$\SeifO{}^{\scriptscriptstyle\ge}$ (--- adjacent to $\SeifX{}^{+}$)}%
\index{O@$\SeifO{}$ (set of Seifert cycles)!$\SeifO{}^{\scriptscriptstyle<}$ (--- not adjacent to $\SeifX{}^{+}$)}%
\index{O@$\SeifO{}$ (set of Seifert cycles)!$\SeifO{}^{\scriptstyle\emptyset}$ (--- not adjacent to $\SeifX{}$)}%
\end{asparaenum}
\end{definitions}

\begin{figure}
\centering
\includegraphics[width=.4\textwidth]{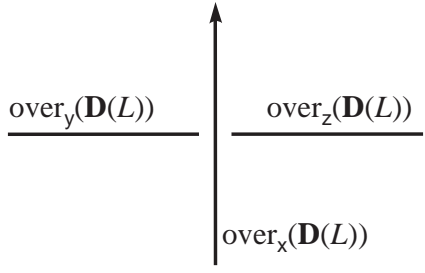}
\caption{Labelled over-arcs (the unlabelled crossing 
need not be $\Seifx$).
\label{Wirtinger relation figure}}
\end{figure}

\begin{theorem}\label{knotgroups}
Let $\diagram{L}$ be a standard link diagram.
If $C\sub\SeifX{L}$ is any set of minimal cardinality
such that $\SeifX{L}\sub\cup_{\Seifx\in C}\proj{L}(\overarc{\Seifx}{L})$,
then the knotgroup of $L$ has a Wirtinger presentation 
\begin{equation}\label{Wirtinger presentation of knotgroup}
\smash{%
\quad
\fundgpnbp{\diff{S^3}{L}}=
   \gp{
      \freegen{\Seifx} \, (\Seifx\in C),
      \freegen{\Seifo} \, (\Seifo\in 
            \SeifO{\diagram{L}}^{\scriptstyle\emptyset})
   }
   {%
      \relator{\Seifx} \,(\Seifx\in\SeifX{L})
   }
}
\end{equation}
where $\relator{\Seifx}\isdefinedas 
\conjugate{\freegen{\Seifx}}{\freegen{\Seify}}%
\freegen{\Seifz}^{-1}$ in case $\pict{L}$ looks locally
like \Figref{Wirtinger relation figure} near $\Seifx$.
\end{theorem}

As \citet{Epple1995} points out, \citet{Wirtinger1905} discovered
\dispeqref{Wirtinger presentation of knotgroup} in the course of
a study of the topology of holomorphic curves.  A proof of 
\eqref{knotgroups} is given by \citet{CrowellFox1977}.

A \moredef{Seifert surface}%
     \index{Seifert!surface}%
     \index{surface!Seifert} is a surface $S\sub S^3$.  
The boundary of a Seifert surface $S$ is a link $L$, 
and $S$ is called a Seifert surface \moredef{for} 
$L$\index{surface!Seifert!for a link}.  Similarly, 
a ribbon \resp{clasp} surface $S$ with $L=\Bd S$ 
is called a \moredef{ribbon} \resp{\moredef{clasp}} \moredef{surface for}%
     \index{surface!ribbon!for a link}%
     \index{surface!clasp!for a link}%
     \index{ribbon!surface!for a link}%
     \index{clasp!surface!for a link} $L$.  
It is a well known fact \citep[apparently first stated by][]%
{FranklPontrjagin1930} that, if $L$ is a link, then there exist 
Seifert surfaces (and, \foreign{a fortiori}, ribbon surfaces 
and clasp surfaces) for $L$.  For the purposes of this survey
the construction sketched by \citet{Seifert1934,Seifert1935} is 
especially convenient.  

\begin{construction}\label{Seifert's construction}
Let $\diagram{L}$ be a standard link diagram for a link 
$L\sub\R^3\sub S^3$, equipped with a fixed smoothing 
$\smooth{\pict{L}}$ of $\pict{L}$ given by a fixed 
family $\{\delta_\Seifx\from\bow{}\to\C %
\Suchthat \Seifx \in \SeifX{\diagram{L}}\}$
of embeddings of $\bow{}$.  To implement Seifert's construction,
choose embeddings
$\hparam{\Seifo}\from\inside{\Seifo}\embto \C\times\Rnonpos$
for $\Seifo\in\SeifO{\diagram{L}}$, and 
$\hparam{\Seifx}\from\h1{}\embto\C\times\Rnonneg$
for $\Seifx\in\SeifX{\diagram{L}}$, subject to the following
conditions.
\begin{figure}
\centering
\includegraphics[width=.3825\textwidth]{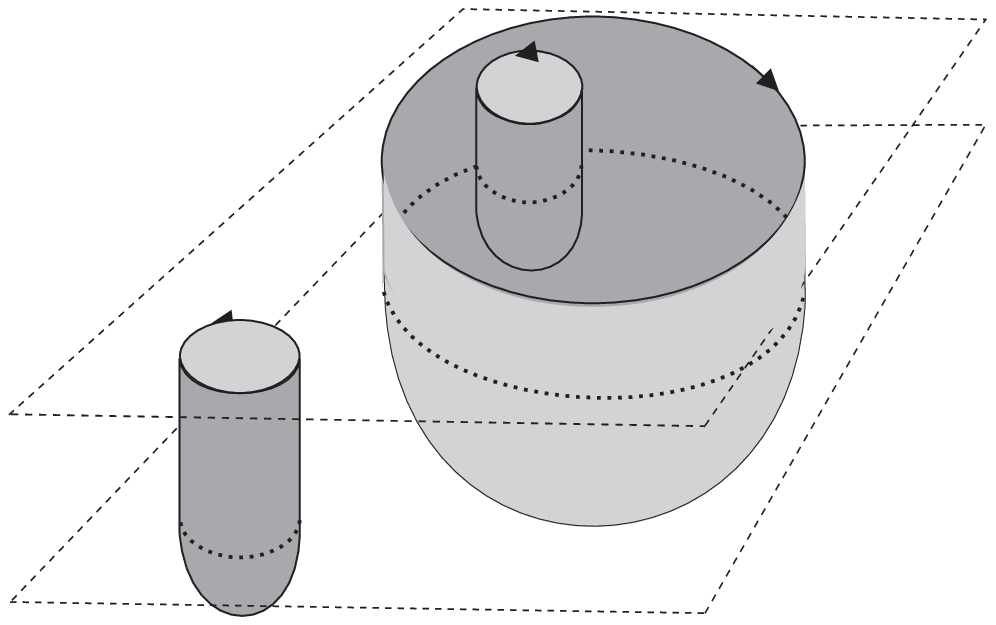}
\caption{Seifert cycles $\Seifo_i$ and 
proper $2$-disks $\hparam{\Seifo_i}(\inside\Seifo_i)$, with
$\Seifo_1\encloses\Seifo_2$, $\Seifo_1\not\encloses\Seifo_3$.
\label{Seifert h0}}
\end{figure}
{\parskip6pt
As suggested in \Figref{Seifert h0},
for each Seifert cycle $\Seifo$,
\begin{enumerate*}
     \item\label{Seifert: 0-handles proper} 
     $\hparam{\Seifo}$ is proper relative to the standard collarings 
     of $\Seifo$ in $\inside{\Seifo}$ and 
     $\C\times\{0\}$ in $\C\times\Rnonpos$,
     \item
     $\pr_1(\hparam{\Seifo}(\diff{\inside{\Seifo}}%
          {\Int{\collar{\Seifo}{\inside{\Seifo}}}})=\inside{\Seifo}
     \sub \C\times\{0\} \sub \C\times\Rnonpos$, and
     \item\label{Seifert: 0-handles disjoint} 
     if $\Seifo'\ne\Seifo$, then $\hparam{\Seifo}(\inside{\Seifo})$ 
     and $\hparam{\Seifo'}(\inside{\Seifo'})$ are disjoint.  
\suspend{enumerate*}
As suggested in \Figref{Seifert h1},
for each crossing $\Seifx$,%
\begin{figure}
\centering
\includegraphics[width=.8\textwidth]{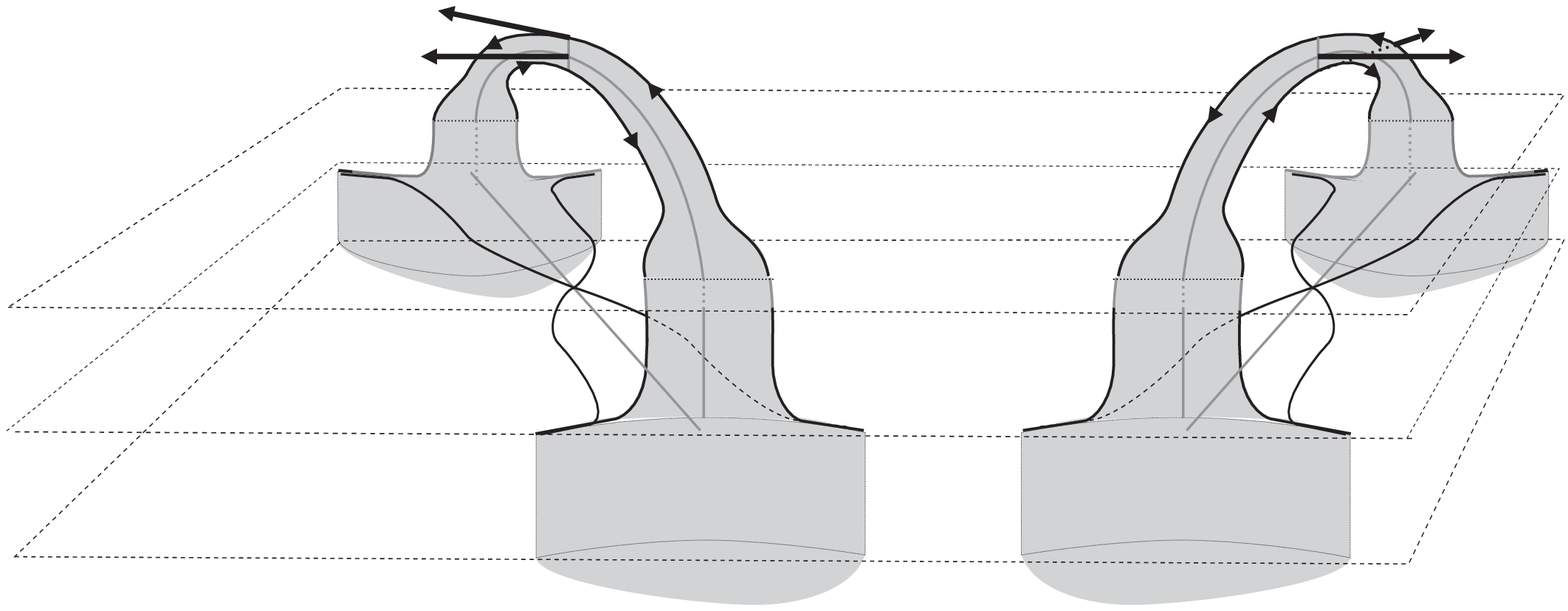}
\caption{Attaching $\hparam{\Seifx_i}(\h1{})$ to 
$\bigcup_{\Seifo\in\SeifO{\diagram{L}}}\Seifo$
(left, positive; right, negative).
\label{Seifert h1}}
\end{figure}%
\resume{enumerate*}
     \item\label{Seifert: 1-handles proper} 
     $\hparam{\Seifx}\from\h1{}\embto\C\times\Rnonneg$
     is proper relative to the standard 
     collarings of $\Int{\attach{\h1{}}}$ in $\h1{}$
     and $\C\times\{0\}$ in $\C\times\Rnonneg$,
     \item\label{Seifert: 1-handles map to bowties} 
     $\pr_1(\hparam{\Seifx}(\diff{\h1{}}%
                            {\Int{ 
                             \collar{\attach{
                                          \h1{}
                                         }
                                 }
                                 {\h1{}
                                 }}
                            })
     =\delta_{\Seifx}(\bow{})$, 
     \item\label{Seifert: 1-handle core arcs}
     $\hparam{\Seifx}(\h1{})\cap\proj{L}^{-1}(\Seifx)=\core{\h1{}}$, 
     \item\label{Seifert: 1-handles oriented} 
     $\hparam{\Seifx}\restr{\attach{\h1{}}}\from
              \hparam{\Seifx}(\attach{\h1{}})\embto
              \bigcup_{\Seifo\in\SeifO{\diagram{L}}}\Seifo$ is orientation
     reversing, and
     \item\label{Seifert: sign condition}
     the sign of the crossing of 
     $\hparam{\Seifx}(\diff{\Bd\h1{}}{\attach{\h1{}}})$ 
     is equal to $\sign{\Seifx}$.
\end{enumerate*}
}
It follows that 
\begin{equation}\label{equation:Seifert's construction}
\smash{%
\Sigma\isdefinedas
\bigcup_{\Seifo\in\SeifO{\diagram{L}}}\hparam{\Seifo}(\inside{\Seifo})
\cup 
\bigcup_{\Seifx\in\SeifX{\diagram{L}}}\hparam{\Seifx}(\h1{})
\sub \C\times\R
}
\end{equation}
is a \nhandle{(0,1)} decomposition 
\dispeqref{handle decomposition of a surface} of a surface.
\end{construction}

\begin{unproved}{proposition}\label{Seifert's construction works}
\begin{inparaenum}
\item\label{works:existence}
There exists a diffeomorphism $\delta\from \C\times\R\to\R^3$ 
such that $\proj{L}\after\delta=\pr_1\from\C\times\R\to\C$
and $L=\delta(\Bd\Sigma)$.
\item\label{works:uniqueness}
The Seifert surface 
$\delta(\Sigma)\sub \R^3\sub S^3$ for $L$ is independent of $\delta$,
up to isotopy fixing $L$ pointwise.
\end{inparaenum}
\end{unproved}

\begin{figure}
\centering
\includegraphics[width=.8\textwidth]{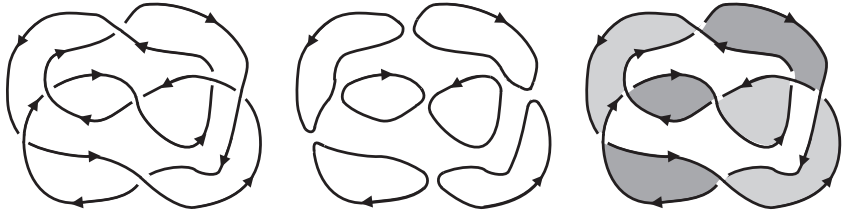}
\caption{Seifert's construction applied to a scattered diagram.
\label{scattered}}
\end{figure}

\begin{figure}
\centering
\includegraphics[width=.8\textwidth]{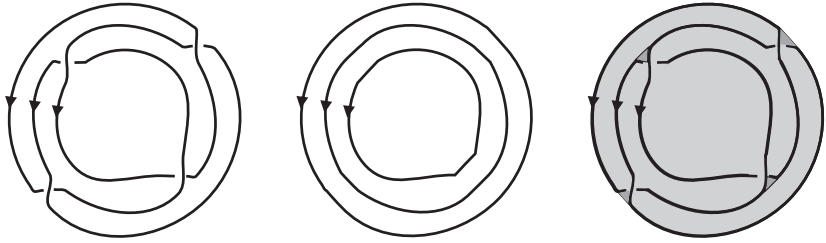}
\caption{Seifert's construction applied to a nested diagram.
\label{nested}}
\end{figure}

Any Seifert surface for $L$ of the form $\delta(\Sigma)$ in
\eqref{Seifert's construction works}(\ref{works:uniqueness}) is 
denoted $\Seifsurf{L}$, and called a
\bydef{diagrammatic Seifert surface}%
    \index{Seifert!surface!diagrammatic}%
    \index{surface!Seifert!diagrammatic}%
    \index{S(D())@$\Seifsurf{\phold}$ (diagrammatic Seifert surface)} for 
$L$.  (A Seifert surface need not be isotopic to any diagrammatic Seifert 
surface.\label{not every Seifert surface is diagrammatic}) \Figref{scattered} 
and \Figref{nested} depict two diagrammatic
Seifert surfaces.

Many operations on links are (most conveniently, and sometimes
necessarily) defined using Seifert surfaces.  Here are two
examples: 
a \moredef{connected sum}%
     \index{connected sum ($\connsum$)!of links}%
     \index{sum!connected ($\connsum$)!of links}%
     \index{$connected$@$\connsum$ (connected sum)} 
of links $L_1, L_2$ bounding Seifert surfaces $S_1, S_2$
is $L_1\connsum L_2\isdefinedas 
\Bd (S_1\bdconnsum S_2)\sub S^3\connsum S^3\iso S^3$;
the \bydef{split sum}%
     \index{sum!split ($\splitsum$), of links}%
     \index{$split$@$\splitsum$ (split sum of links)} of 
$L_1, L_2$ is $L_1\splitsum L_2\isdefinedas \Bd(S_1\connsum S_2)$.  
It is well-known that if $L_1=K_1$ and $L_2=K_2$ are knots, 
then $K_1\connsum K_2$ is well-defined up to isotopy, and 
independent of $S_1$ and $S_2$; in any case, $L_1\splitsum L_2$ 
is well-defined.  In particular, for $n\in\N$ the 
\ncomponent{n} unlink is the split sum of $n$ unknots.
For any link $L$, let $L^{(n)}$%
   \index{Ln@$L^{(n)}$ (connected sum $L\connsum O^{(n)}$)} denote 
the (well-defined) link $L\connsum O^{(n)}$.  

Similarly, many knot and link invariants are defined 
using Seifert surfaces.  

\begin{unproved}{proposition}
If $L$ and $L'$ are disjoint links,
then the algebraic number of intersections of $L'$ with a Seifert
surface $S$ for $L$ is independent of $S$ \textup(provided only
that $L'$ intersects $S$ transversely\textup).
\end{unproved}

This integer invariant of the pair $(L,L')$, denoted 
$\link{L}{L'}$\index{link@$\link{\phold}{\phold}$ (linking number)} and
called the \bydef{linking number} of $L$ and $L'$, satisfies 
$\link{L}{L'}=\link{L'}{L}=-\link{-L}{L'}\linebreak[3]
=-\link{\Mir L}{\Mir L'}$.

\begin{definitions}\label{genus and other invariants}
Let $K$ be a knot, $L$ a link.  Define invariants
\begin{align*}
\genus{}K\isdefinedas \min\{\genus{}{S} &
     \Suchthat S \text{ is a Seifert surface for $K$}\},\\
\genus{r}{K}\isdefinedas \min\{\genus{}{S} &
     \Suchthat S \text{ is a ribbon surface for $K$}\}, \\
\genus{s}{K}\isdefinedas \min\{\genus{}{S} &
     \Suchthat S \text{ is a slice surface for $K$}\},\\
\eulermax{}{L}\isdefinedas \max\{\euler{S} &
     \Suchthat S \text{ is a Seifert surface for $L$}\},\\
\eulermax{r}{L}\isdefinedas \max\{\euler{S} &
     \Suchthat S \text{ is a ribbon surface for $L$}\}, \\
\eulermax{s}{L}\isdefinedas \max\{\euler{S} &
     \Suchthat S \text{ is a slice surface for $L$}\},\\
\clasp{L}\isdefinedas \min\{\clasp{S} &
     \Suchthat S=f(D^2\times\compsof{L}) \text{ is a clasp surface for $L$ }\},\\
\node{L}\isdefinedas \min\{\node{S} &
     \Suchthat S=f(D^2\times\compsof{L}) \text{ is a nodal surface for $L$}\}.
\end{align*}%
Call $\genus{}{K}$ \resp{$\genus{r}K$; $\genus{s}K$} 
the \moredef{genus} \resp{\moredef{ribbon genus}; \moredef{slice genus}}%
     \index{genus!of a knot ($\genus{}{\phold}$)}%
     \index{genus!ribbon, of a knot ($\genus{r}{\phold}$)}%
     \index{genus!slice, of a knot ($\genus{s}{\phold}$)}%
     \index{knot(s)!genus ($\genus{}{\phold}$)}%
     \index{ribbon!genus of a knot ($\genus{r}{\phold}$)}%
     \index{slice!genus of a knot ($\genus{s}{\phold}$)}%
     \index{knot(s)!genus ($\genus{}{\phold}$)}%
     \index{knot(s)!ribbon genus ($\genus{r}{\phold}$)}%
     \index{knot(s)!slice or Murasugi genus ($\genus{s}{\phold}$)}%
     \index{g@$\genus{}{\phold}$ (genus)!of a knot}%
     \index{gr@$\genus{r}{\phold}$ (ribbon genus of a knot)}%
     \index{gs@$\genus{s}{\phold}$ (slice, or Murasugi, genus of a knot)} of 
$K$, and say $K$ is a \moredef{slice} \resp{\moredef{ribbon}}%
    \index{slice!knot}%
    \index{knot(s)!slice}%
    \index{ribbon!knot}%
    \index{knot(s)!ribbon} knot 
in case $\genus{s}{K}=0$ \resp{$\genus{s}{K}=0$}.
Another name for $\genus{s}{K}$ is the ``Murasugi genus'' of $K$.  
\end{definitions}

\begin{definition}
Let $\diagram{L}$ be a standard link diagram, 
$n\isdefinedas\card{\compsof{L}}$.  It is easy to prove
and well known \citep[see \citeauthor{Hoste2004} or][]{Kauffman2004} 
that there is a standard link diagram $\diagram{O^{(n)}}$ for an 
unlink $O^{(n)}$ 
such that $\pict{O^{(n)}}=\pict{L}$ and 
$\info{O^{(n)}}$ differs from
$\info{L}$ precisely to the extent that some number 
$u\ge 0^{(n)}$ of
crossings of $\pict{O^{(n)}}=\pict{L}$ have opposite signs in 
$\info{O^{(n)}}$ and $\info{L}$.  The \moredef{unknotting number}%
    \index{unknotting number!of a standard link diagram} of
$\diagram{L}$ is the least such $u$.  The \moredef{unknotting number}%
    \index{unknotting number!of a link ($\usz(\phold)$)} of $L$
is the least unknotting number of all standard link diagrams for $L$.
(The unknotting number is also called the \moredef{\"Uberschneidungszahl}%
    \index{Uberschneidungszahl@\"Uberschneidungszahl|see {unknotting number}} %
\citet{Wendt1937,Milnor1968} and the \moredef{Gordian number}%
    \index{Gordian!number|see {unknotting number}} %
\citet{BoileauWeber1982,Bennequin1993,A'Campo1998}.)
More generally, the \moredef{Gordian distance}%
   \index{Gordian!distance} $d_G(L,L')$ between two links
$L$, $L'$ with $\card{\compsof{L}}=\card{\compsof{L'}}$ is the minimum
number of sign changes at crossings needed to transform some standard
link diagram $\diagram{L}$ to some standard link diagram 
$\diagram{L'}$ \citep{Murakami1985}; thus $\usz(L)=d_G(L,O^{(n)})$.
\end{definition}

Various more or less obvious inequalities 
relate $\usz$ and the several invariants named in 
\eqref{genus and other invariants} 
\citep[see][]{Shibuya1974,Rudolph1983b}.

\subsection{Framed links; Seifert forms}\label{framed links; Seifert forms}

Let $L$ be a link.  A \bydef{framing} of $L$ is 
a function $f\from \compsof{L}\to\Z$; 
the pair $(L,f)$ is a \moredef{framed link}%
    \index{framed!link}.  A framing of a knot is identified with the 
integer which is its value.  For framings $f,g$ of $L$, 
write $f\preccurlyeq g$%
    \index{$less$@$\preccurlyeq$ (less twisted framing)}, and
say $f$ is \moredef{less twisted}%
    \index{framing!less twisted than ($\preccurlyeq$)} than $g$
provided that $f(K)\leq g(K)$ for every $K\in\compsof{L}$.

\begin{unproved}{proposition}\label{normal bundle of a link}
The normal bundle of $L$ is trivial.  Given 
a Seifert surface $S$ for $L$, there exists a trivialization 
$n\from L\times\C\to\normbun{L}$ such that,
under the identification of $\Nb{L}{S^3}$ 
with $n(L\times D^2)$ \textup(as in 
\eqref{normal bundle stuff}.\eqref{adapted maps to the circle}\textup),
$S\cap\Nb{L}{S^3}=\Nb{L}{S}$ is identified 
with $n(L\times\clint{0}{1})$.
The homotopy class of $n$ is well-defined, independent of $S$.
\end{unproved}

Let $(K,k)$ be a framed knot.
A $k$-\bydef{twisted annulus of type $K$} 
is any annulus $\AKn{K}{k}\sub S^3$ such that 
$K\sub\Bd\AKn{K}{f}$ and $\link{K}{\diff{\Bd\AKn{K}{k}}{K}}$ is $-k$; 
note that, since $\diff{\Bd\AKn{K}{k}}{K}$ is clearly isotopic 
to $-K$, all four of $\AKn{K}{k}$, $\AKn{-K}{k}$, $-\AKn{K}{k}$, 
and $-\AKn{-K}{k}$ are isotopic.  For a framed link $(L,f)$, 
$\AKn{L}{f}$ is defined componentwise.
Given a \nsubmanifold{2} $S\sub S^3$ and a link $L\sub S$,
the $S$-\moredef{framing of $L$}%
     \index{S-framing@$S$-framing} 
is the framing $\Seifframe{L}{S}$ such that 
$\collar{L}{S}=\AKn{-L}{\Seifframe{L}{S}}$. 
A framed link $(L,f)$ is \moredef{embedded}%
     \index{framed!link!embedded on a Seifert surface} on 
a Seifert surface $S$ in case $L\sub S$ and $f=\Seifframe{L}{S}$.

Let $S$ be a Seifert surface with collaring $\col{S}{S^3}$.
The \moredef{Seifert pairing}%
    \index{Seifert!pairing}%
     \index{pairing, Seifert} (on $S$) 
of an ordered pair of links $(L_0, L_1)$ with $L_0, L_1\sub S$ is 
$\Seifpair{L_0}{L_1}{S}\isdefinedas\link{L_0}{\pushoff{L_1}}$;
if $K\sub S$ is a knot, then $\Seifpair{K}{K}{S}=\Seifframe{K}{S}$.  
Given an ordered \ntuple{\mu}
$\ontuple{L}{\mu}$ of links on $S$ the homology classes of which 
form a basis for $H_1(S;\Z)$, the 
\moredef{Seifert matrix}%
     \index{Seifert!matrix}%
     \index{matrix, Seifert} of $S$ 
with respect to that basis is the $\mu\times\mu$ matrix
$\Seifmatrix{\Seifpair{L_i}{L_j}{S}}$, and the 
\moredef{Seifert form}%
     \index{Seifert!form} %
is the (typically non-symmetric) bilinear form on $H_1(S;\Z)$
represented by $\Seifmatrix{\Seifpair{L_i}{L_j}{S}}$.

\subsection{Fibered links, fiber surfaces, and open books}

Let $L$ be a link. Let $n\from L\times\C\to\normbun{L}$ be 
a trivialization, as in \eqref{normal bundle of a link},
in the homotopy class corresponding to any Seifert
surface $S$ for $L$.  Call $L$ \moredef{fibered}%
    \index{fibered link}%
    \index{link(s)!fibered} in case there exists a map 
$\phi\from\diff{S^3}{L}\to S^1$ (called a \moredef{fiber map}%
    \index{fiber!map}%
    \index{map!fiber} for $L$)
which is adapted to $n$, has $d(K)=1$ for all 
$K\in\compsof{L}$, and is a fibration (in particular, a Morse map).  
If $L$ is a fibered link with fiber map $\phi$,
then for each $e^{\imunit\theta}\in S^1$, 
$L\cap \phi^{-1}(e^{\imunit\theta})$ is a Seifert surface for 
$L$.  A \moredef{fiber surface}%
    \index{fiber!surface}%
    \index{surface!fiber} is any Seifert surface 
$\Milnorfiber{L}$ isotopic 
to $L\cap\phi^{-1}(e^{\imunit\theta})$ for any fibered link $L$ with 
fiber map $\phi$.  The \moredef{Milnor number}%
    \index{Milnor!number of a fibered link ($\Milnornumber{\phold}$)} of $L$
is $\Milnornumber{L}\isdefinedas \dim_\R H_1(\Milnorfiber;\R)$.

Let $S$ be a Seifert surface.
The \moredef{top}%
    \index{top, of a Seifert surface $S$ ($\top{S}$)} of 
$S$ is $\top{S}\isdefinedas \collar{S}{S^3}$.
A \ndisk{2} $D\sub\top{S}$ is a \moredef{top-}\bydef{compression disk} 
in case $\Bd D=D\cap S$ and $\Bd D$ bounds no disk on $S$.
Call $S$ \moredef{compressible}%
    \index{surface!Seifert!compressible}%
    \index{Seifert!surface!compressible}%
    \index{compresible Seifert surface} in case there exists 
a top-compression disk for at least one of $S$ and $-S$,
and \moredef{incompressible}%
    \index{surface!Seifert!incompressible}%
    \index{Seifert!surface!incompressible}%
    \index{incompresible Seifert surface} in case it is not compressible.
Call $S$ \moredef{least-genus}%
    \index{surface!Seifert!least-genus}%
    \index{Seifert!surface!least-genus}%
    \index{least-genus Seifert surface} provided that 
$\euler{S}=\eulermax{}{L}$.
The following facts are well known \citep[see][]%
{Stallings1978,Gabai1983a,Gabai1983b,Gabai1986}.

\begin{unproved}{proposition}\label{fiber facts}
\begin{inparaenum}
\item
$S$ is a fiber surface if and only if 
$S$ is connected and a push-off map induces an isomorphism 
$\fundgp{\Int{S}}{\basept}\to 
\fundgp{\diff{S^3}{S}}{\pushoff{\basept}})$.
\item
A least-genus surface $S$ is incompressible.
\item
A fiber surface is least-genus, and up to isotopy 
it is the unique incompressible surface with its boundary.
\item
$\AKn{K}{n}$ is least-genus if and only if $(K,n)\ne (O,0)$.
\item
$\AKn{K}{n}$ is a fiber surface if and only if $(K,n)=(O,\mp1)$.
\end{inparaenum}
\end{unproved}

The fiber surface $\AKn{O}{-1}$ \resp{$\AKn{O}{1}$} is called a
\moredef{positive} \resp{\moredef{negative}} 
\bydef{Hopf annulus}%
     \index{Hopf annulus!positive}%
     \index{Hopf annulus!negative}
(sometimes ``Hopf band''); the choice of the adjectives 
``positive'' and ``negative'' reflects the linking number 
of the components of $\Bd \AKn{O}{\mp 1}$.

Fibered links are slightly flaccid.  They may be rigidified
as follows.  An \moredef{open book}%
    \index{open!book} is a map 
$\ob{f} \from S^3 \to \C$ such that $0$ is a regular value and 
$\arg{\ob{f}}\isdefinedas\argument{\ob{f}}%
\from \diff{S^3}{\binding{\ob{f}}}\to S^1$ 
is a fibration.  The \bydef{binding} $\binding{f}$ 
of $\ob{f}$ clearly is a fibered link, and 
for each $e^{\imunit\theta}\in S^1$, 
the \moredef{$\theta$th} \bydef{page} 
$F_\theta\isdefinedas \ob{f}^{-1}(\{re^{\imunit\theta}\suchthat r\ge 0\})$ 
of $\ob{f}$ is a fiber surface. Every fibered link is the 
binding of various open books; any two fibered books with
the same binding are equivalent in a straightforward sense
\citep[cf.][]{KauffmanNeumann1977}. 

\citet{Milnor1968} discovered a rich source (now called
\moredef{Milnor fibrations}%
    \index{Milnor!fibration}; see 
\eqref{Milnor maps and fibrations}) of open books as part of
his investigation of the topology of singular points of
complex hypersurfaces.  The simplest special cases are
fundamental to the knot theory of complex plane curves
and easy to write down.  Let $m, n \in\N$, $(m,n)\ne(0,0)$.

\begin{unproved}{theorem}\label{torus link Milnor map}
$\ob{o}\{m,n\}\from S^3\to\C^2\suchthat (z,w)\mapsto z^m+w^n$
is an open book.  
\end{unproved}

The binding $\ob{o}\{m,n\}^{-1}(0)$ is a \bydef{torus link} 
\moredef{of type $\{m,n\}$}, sometimes 
\citep[as in][cf. \citeauthor{Litherland1979}, 
\citeyear{Litherland1979}]{Rudolph1982,Rudolph1988} denoted $O\{m,n\}$.
Call $\ob{o}\isdefinedas\ob{o}\{1,0\}$ 
\resp{$\ob{o}'\isdefinedas\ob{o}\{0,1\}$} 
the \moredef{vertical} \resp{\moredef{horizontal}}
\moredef{unbook}%
    \index{unbook!vertical ($\ob{o}$)}%
    \index{unbook!horizontal ($\ob{o}'$)}%
    \index{vertical unbook ($\ob{o}$)}%
    \index{horizontal unbook ($\ob{o}$)} and its binding
$\vunknot\isdefinedas O\{1,0\}$ 
\resp{$\hunknot\isdefinedas O\{0,1\}$} the \moredef{vertical} 
\resp{\moredef{horizontal}} \moredef{unknot}%
    \index{unknot!vertical ($\vunknot$)}%
    \index{unknot!horizontal ($\hunknot$)} \citep{Rudolph1988}.

\subsection{Polynomial invariants of knots and links}
\label{polynomial invariants of knots and links}

The intent of this section is to establish notations 
and conventions, and to state without proof two useful
theorems.  For a thorough treatment of polynomial link 
invariants, see \citet{Kauffman2004}.

For any ring $\Rs$, 
for any Laurent polynomial $H(s)\in \Rs[s^{\pm 1}]$, write 
$\ord_s H(s)\isdefinedas 
\sup\{ n\in\Z\Suchthat s^{-n}H(s)\in\Rs[s]\sub\Rs[s^{\pm 1}]\}$, 
$\deg_s H(s)\isdefinedas -\ord_s H(s^{-1})$.

\begin{definition} 
Let $K$ be a knot.  Let $S$ be a Seifert surface for $K$.
Let $A=\Seifmatrix{\Seifpair{L_i}{L_j}{S}}$
be a Seifert matrix of $S$, with transpose $\transpose{A}$.  
The \moredef{unnormalized Alexander polynomial}%
    \index{Alexander polynomial ($\Alex{\phold}$)!unnormalized}%
    \index{polynomial!Alexander ($\Alex{\phold}$)!unnormalized} of $K$ 
is $\det(\transpose{A}-tA)\in\Z[t]$.  It is easily shown that
$\mu\isdefinedas\deg_t \det(\transpose{A}-tA)$ is even.  
The \moredef{Alexander polynomial}%
    \index{Alexander polynomial ($\Alex{\phold}$)}%
    \index{polynomial!Alexander ($\Alex{\phold}$)}%
    \index{Delta@$\Alex{\phold}$ (Alexander polynomial)} of $K$
is $\Alex{K}(t)\isdefinedas t^{-\mu/2}\det(\transpose{A}-tA)
\in\Z[t,t^{-1}]$.
\end{definition}

\begin{unproved}{proposition}
$\Alex{K}$ depends only on $K$, not on the choice
of $S$ or $A$.
\end{unproved}

Let $L\sub\R^3\sub S^3$ be a link.  Let $\diagram{L}$ be 
a standard link diagram for $L$.  Let $\Seifx\in\SeifX{\diagram{L}}$.
Call $\Seifx$ \moredef{homogeneous}%
    \index{homogeneous crossing (of a link diagram)}%
    \index{crossing!of a standard link diagram!homogeneous} \resp{%
\moredef{heterogeneous}}%
    \index{heterogeneous crossing (of a link diagram)}%
    \index{crossing!of a standard link diagram!heterogeneous} 
in case $\card{\compsof{\smoothloc{L}{\Seifx}}}$
equals $\card{\compsof{L}}+1$ \resp{$\card{\compsof{L}}-1$}.  
In case $\Seifx\in\SeifX{\diagram{L}}^{+}$, let
$L_{+}\isdefinedas L$, and define standard link diagrams 
$\diagram{L_{-}}$ and $\diagram{L_0}$ (and thereby links
$L_{-}$ and $L_{0}$) as follows:
\begin{inparaenum}
\item
$\pict{L_{-}}=\pict{L_{+}}$, and $\info{L_{-}}$ differs from
$\info{L_{+}}$ precisely to the extent that
$\Seifx\in\SeifX{\diagram{L_{-}}}^{-}$;
\item
$\pict{L_{0}}$ is the local smoothing 
$\smoothloc{\pict{L_{\pm}}}{\Seifx}$ of $\pict{L_{\pm}}$ at $\Seifx$,
and $\info{L_{0}}$ differs from $\info{L_{\pm}}$ precisely to the extent 
that $\Seifx\notin\SeifX{\diagram{L_{0}}}$.
\end{inparaenum}
In case $\Seifx\in\SeifX{\diagram{L}}^{-}$, 
modify these definitions accordingly; the two sets of definitions are 
consistent. 

\begin{figure}
\begin{minipage}{.49\textwidth}
\setlength{\unitlength}{\textwidth}
\centering
\begin{picture}(.8175,0.2325)(0,-0.03)
\includegraphics[width=.8175\textwidth]{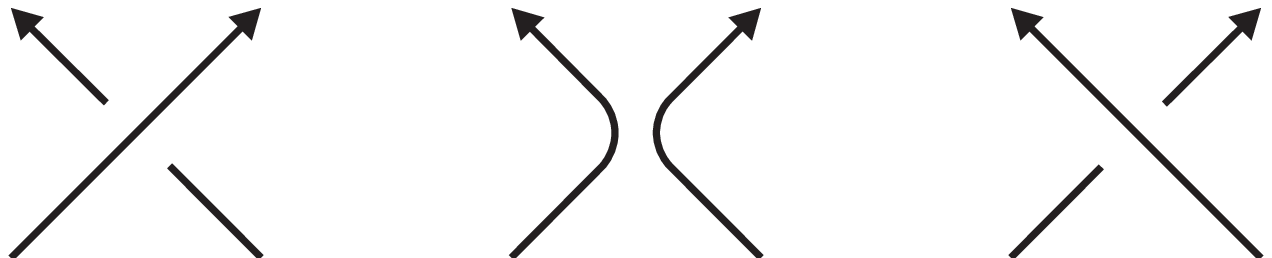}
\end{picture}
{\scriptsize
\put(-.75,0){$L_{+}$}\put(-0.4275,0){$L_0$}
\put(-.105,0){$L_{-}$}
}
\end{minipage}
\begin{minipage}{.49\textwidth}
\setlength{\unitlength}{\textwidth}
\centering
\begin{picture}(.49,.231)(0.0,-0.03)
\includegraphics[width=.49\textwidth]{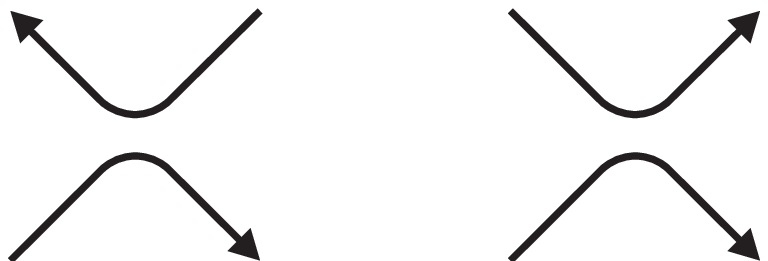}
{\scriptsize
\put(-.4275,-.03){$(1)$}
\put(-0.1095,-.03){$(2)$}
}
\end{picture}
\end{minipage}
\setlength{\unitlength}{\textwidth}
\put(-.755,-.1){$\text{\small{(A)}}$}
\put(-.255,-.1){$\text{\small{(B)}}$}
\caption{(A) $\protect\knot{L}{+}$, 
$\protect\knot{L}{0}$, and $\protect\knot{L}{-}$.
(B) $\protect\knot{L}{\infty}$ (homogeneous and 
heterogeneous cases).\label{local diagrams}}
\end{figure}

Up to isotopy, the local situation at $\Seifx$ 
is as in \Figref{local diagrams}(A).  
In case $\Seifx$ is homogeneous \resp{heterogeneous}, 
let $\diagram{L_\infty}$ be the standard link diagram 
differing from $\diagram{L_\pm}$ and $\diagram{L_0}$ only
as required by case $(1)$ \resp{case $(2)$} of 
\Figref{local diagrams}(B), 
let $p$ \resp{$q$} be the linking number of the right-hand visible 
component of $L_0$ with the rest of $L_0$ \resp{of the lower 
visible component of $L_{+}$ with the rest of $L_{+}$}, 
and define $r\isdefinedas 4p+1$ \resp{$r\isdefinedas 4q-1$}.%

\begin{definition}
The \moredef{oriented polynomial}%
    \index{oriented polynomial}%
    \index{polynomial!oriented} $P_L(v,z)\in\Z[v^{\pm1},z^{\pm1}]$
and \moredef{semi-oriented polynomial}%
    \index{semi-oriented polynomial}%
    \index{polynomial!semi-oriented} $F_L(a,x)\in\Z[a^{\pm1},x^{\pm1}]$
of $L$ are defined recursively as follows,
with the initial conditions $P_O(v,z)=1=F_O(a,x)$.%
\begin{align}
\label{oriented recursion}
P_{L_{+}}(v,z) & =vzP_{L_0}(v,z)+v^2P_{L_{-}}(v,z)\\
\label{semi-oriented recursion}
F_{L_{+}}(a,x) & =a^{-1}xF_{L_0}(a,x)-a^{-2}F_{L_{-}}(a,x)
                 +a^{-r}x F_{L_\infty}(a,x)
\end{align}
The nomenclature is that of \citet{Lickorish1986};
the choice of variables $v, z$ in 
\dispeqref{oriented recursion} follows \citet{Morton1986}.
The oriented \resp{semi-oriented} polynomial is often 
known, eponymously, as the \moredef{FLYPMOTH}%
    \index{FLYPMOTH polynomial}%
    \index{polynomial!FLYPMOTH} \citep{Freydetal1985,PrzytyckiTraczyk1988} %
\resp{\moredef{Kauffman} \citep{Kauffman1987}}%
    \index{Kauffman polynomial}%
    \index{polynomial!Kauffman} polynomial.
\end{definition}

\begin{definitions}\label{other link polynomials}
Several other polynomials, though mere adaptations of
the oriented or semi-oriented polynomials, nonetheless
have their uses.
\begin{asparaenum}[\bf \thedefinitions.1.]
\item\label{framed polynomial}
Let $(L,f)$ be a framed link.
The \moredef{framed polynomial}%
    \index{framed!polynomial}%
    \index{polynomial!framed} $\{L,f\}(v,z)\in\Z[v^{\pm1},z^{\pm1}]$ is
\begin{equation}
\smash[t]{
(-1)^{\card{\compsof{L}}}
(1+(v^{-1}-v)z^{-1} 
\sum_{\emptyset\ne C\sub\compsof{L}}
     (-1)^{\card{C}}P_{\Bd A(\cup C, f\restr{\cup C})})
}
\end{equation}
\citep[see][]{Rudolph1990}.
For any $(L,f)$, 
$\smash[t]{\{L,f\}=v^{-2\sum_{K\in\compsof{L}} f(K)}\{L,0\}}$.

\item\label{R-polynomial}
Let $R_L(v)\isdefinedas 
\left.(z^{\card{\compsof{L}}-1}P_L(v,z))\right|_{z=0}$. 
$R_L$ can be calculated from 
$R_O(v)=1$ and 
$R_{L_{+}}(v)=hR_{L_0}(v)+v^2R_{L_{-}}(v)$,
where $h$ is $0$ \resp{$1$} in case $\Seifx$ is 
heterogeneous \resp{homogeneous}.
\item\label{G-polynomial}
Let $F^*_L(a,x)\isdefinedas (F_L \pmod{2})\in(\Z/2\Z)[a^{\pm1},x^{\pm 1}]$.
For $k=0,1$, let $G^k_L(a)\isdefinedas 
\left.(x^{1-c(L)}F^*_L(a,x))\right|_{x=k}
\in(\Z/2\Z)[a^{\pm 1}]$,  
so $G^0_L(a)=R_L(a^{-1}) \pmod{2}$ 
and can be calculated using 
\eqref{other link polynomials}.\eqref{R-polynomial},
while $G^1_L$ can be calculated from 
$G^1_{O}(a)=1$ and 
$G^1_{L{_{+}}}(a)=a^{-2}G^1_{L{_{-}}}(a)+a^{-1}G^1_{L{_0}}(a)+
a^{-r}G^1_{L{_\infty}}(a)$.
\end{asparaenum}
\end{definitions}

The result underlying most applications of polynomial 
invariants to the knot theory of complex plane curves, 
due to \citet{Morton1986} and \citet{FranksWilliams1987},
is rephrased here to fit the expository order of this survey;
the usual statement, in terms of braids, is given
in \eqref{MFW--braids}.

\begin{theorem}\label{MFW--no braids}
If $\diagram{L}$ is a standard link diagram such that
\begin{inparaenum}
\item
$\diagram{L}$ is nested and
\item
$\SeifO{\diagram{L}}=\SeifO{\diagram{L}}^{+}$,
\end{inparaenum}
then 
\begin{equation*}
\quad
\ord_v P_{L}\geq\
\card{\SeifX{\diagram{L}}^{+}}-\card{\SeifX{\diagram{L}^{-}}}
-\card{\SeifO{\diagram{L}}}+1.
\quad
\qquad
\qedsymbol
\end{equation*}
\end{theorem}
The framed polynomial provides a bridge between 
the oriented and semi-oriented polynomials, as the 
following result \citep{Rudolph1990} makes plain.%
\begin{theorem}\label{congruence theorem}
$(1+(v^{-2}+v^2)z^{-2})F_L(v^{-2},z^2)
\equiv v^{4\tau(L)}\{L,0\}(v,z)\pmod{2}.\quad\qedsymbol$\linebreak[0]  
\end{theorem}

\subsection{Polynomial and analytic maps; algebraic and analytic sets}

This section recalls needed definitions from real and complex 
algebraic and analytic geometry, and establishes notations.
General background, and proofs of stated results, can be found in 
\citet{Whitney1957,Whitney1972,Milnor1968,%
Narasimhan1960,GunningRossi1965}, and,
for \eqref{Tarski-Seidenberg Theorem}(\ref{TS statement}), 
\citet[Appendix B]{AbrahamRobbin1967}.

Let $\F$ be one of the fields $\R$ or $\C$, with its metric
topology.  The algebra of polynomials 
\resp{somewhere-convergent power series} in $n$ 
variables with ground field $\F$ is denoted 
$\F[\phi_1,\dots,\phi_n]$ \resp{$\F\{\phi_1,\dots,\phi_n\}$},
where $\phi$ stands for $x$ \resp{$z$} in case $\F$ is $\R$ 
\resp{$\C$}.  As usual, $f\in\F[\phi_1,\dots,\phi_n]$ is
conflated with the \moredef{polynomial function}%
    \index{polynomial!function}%
    \index{function!polynomial} $f\from\F^n\to\F$ that it defines,
and $f\in\F\{\phi_1,\dots,\phi_n\}$ with both the $\F$-analytic 
function that it defines in a neighborhood of $0\in\F^n$ and
the germ of that function at $0\in\F^n$.  (In particular $\phi_s$ 
is conflated with the coordinate projection $\pr_s\from\F^n\to\F$ 
for every $s\in\intsto{n}$, no notational distinction being made 
between the functions $\phi_s$ on $F^m$ and $F^{n}$ so long as 
$s\in\intsto{m}\cap\intsto{n}$.)  Given a non-empty open set
$\Omega\sub\F^n$, a function $f\from\Omega\to\F$ is called
\moredef{$\F$-analytic}%
    \index{analytic!function}%
    \index{C-analytic@$\C$-analytic!function}%
    \index{F-analytic@$\F$-analytic!function}%
    \index{R-analytic@$\R$-analytic!function}%
    \index{function!R-analytic@$\R$-analytic}%
    \index{function!C-analytic@$\C$-analytic (holomorphic)} (simply 
\moredef{analytic} when $\F$ is clear from context; also
\moredef{holomorphic}%
    \index{holomorphic!function}%
    \index{function!holomorphic}%
for $\F=\C$, and \moredef{entire}%
    \index{entire function}%
    \index{function!entire}%
when $\Omega=\C^n$) in case,
for every $(\phi_1^{(0)},\dots,\phi_n^{(0)}\in\Omega$, 
the germ of 
$f(\phi_1-\phi_1^{(0)},\dots,\phi_n-\phi_n^{(0)})$ at $(0,\dots,0)$ 
belongs to $\F\{\phi_1,\dots,\phi_n\}$.  
The set $\holo{\Omega}$%
    \index{O()@$\holo{\phold}$ (analytic functions)} of 
all $\F$-analytic functions on $\Omega$ is an algebra containing 
(a natural isomorphic image of) $\F[\phi_1,\dots,\phi_n]$.

A \moredef{polynomial} \resp{\moredef{$\F$-analytic}}%
    \index{polynomial!map}%
    \index{analytic!map}%
    \index{holomorphic!map}%
    \index{map!polynomial}%
    \index{map!analytic}%
    \index{map!holomorphic}%
    \index{map!F-analytic@$\F$-analytic}%
    \index{C-analytic@$\C$-analytic!map}%
    \index{F-analytic@$\F$-analytic!map}%
    \index{R-analytic@$\R$-analytic!map} \moredef{map} %
$F=(f_1,\dots,f_m)$ from $\F^n$ \resp{$\Omega$} to $\F^m$ 
is one with components $f_s$ that are polynomial 
\resp{$\F$-analytic} functions.  
An \bydef{algebraic set} \resp{\bydef{global analytic set}}%
    \index{V()@$\variety{\phold}$ (algebraic or global analytic set)} is
a subset $\variety{F}\isdefinedas F^{-1}(0,\dots,0)$ of 
$\F^m$ \resp{$\Omega$}, where $F$ is a polynomial 
\resp{analytic} map.  An \moredef{analytic set}%
    \index{analytic!set} is a subset  
$X$ of $\Omega$ for which every point of $\Omega$ has an
open neighborhood $U$ such that $X\cap U$ is a global analytic
set $\variety{F}$ for some $F\in\holo{U}$.  Every global analytic 
set (in particular, every algebraic set) is an analytic set; the
converse fails for many $\Omega$ when $n>1$.

Let $F\from\Omega\to\F^m$ be an analytic map (allowing
the possibility that $\Omega=\F^n$ and $F$ is polynomial).
It may happen that
\[
\reg{F}\isdefinedas
\{\ontuple{y}m\in \variety{F}\Suchthat\rank_\F{DF\ontuple{y}m}=n-m\}
\]
is not dense in $\variety{F}$.  However, 
theorems of commutative algebra show that there is another
analytic map $F_0\from\Omega\to\F^m$ such that 
$\variety{F_0}$ and $\variety{F}$ are equal \shout{as sets} 
(that is, ignoring multiplicities), and 
$\reg{F_0}$ is dense in $\variety{F_0}=\variety{F}$.  
Call $\ontuple{y}m\in \variety{F}$ a 
\moredef{regular \resp{\bydef{singular} point}}%
    \index{regular!point}\index{singular!point} of 
the global analytic set $\variety{F}$ in case 
$\rank_\F{DF_0\ontuple{y}m}$ equals \resp{is less than} $n-m$;
these definitions are independent of the particular choice of
$F_0$.  Regularity in $\variety{F}$ is clearly a local property, 
and is therefore well-defined in any analytic set $X$.
The set $\reg{X}$%
    \index{reg()@$\reg{\phold}$ (regular locus)} of 
regular points of $X$ is called the \moredef{regular locus}%
     \index{regular!locus ($\reg{\phold}$)} of 
$X$; it is an $\F$-analytic manifold.  The \moredef{singular locus}%
     \index{singular!locus ($\sing{\phold}$)} 
$\diff{X}{\reg{X}}\defines\sing{X}$%
     \index{sing()@$\sing{\phold}$ (singular locus)} of $X$
is an algebraic, global analytic, or analytic set according as $X$ is.
Let $\singkth{0}{X}$ and $\regkth{1}{X}$ both mean 
$\reg{X}$; for $s\in\Nstrict$, let 
$\singkth{s+1}{X}\isdefinedas\sing{\singkth{s}{X}}$, 
$\regkth{s+1}{X}\isdefinedas\reg{\singkth{s}{X}}$.  
An $\F$-analytic set $X$ is partitioned 
by the (finitely many) non-empty sets in the sequence 
$\{\regkth{s}{X}\}_{s\in\Nstrict}$.  
In case $X$ is algebraic, the refinement of this partition 
obtained by separating each $\regkth{s}{X}$ into its connected 
components is a finite stratification 
\citep[\citealp{Whitney1957}; see][Theorems 2.3 and 2.4]{Milnor1968}; %
call it the \moredef{\naive{} stratification}%
     \index{naive@\naive{} stratification}
     \index{stratification!\naive{} (of an algebraic set)} of $X$.

A \moredef{basic semi-algebraic set}%
    \index{semi-algebraic set!basic} in $\R^n$ is the intersection 
of an $\R$-algebraic set $\variety{F}$ and finitely many sets 
of the form $G^{-1}(\Rg{0})$, with $G\in\R[x_1,\dots,x_n]$; 
a \bydef{semi-algebraic set} is the union of finitely many 
basic semi-algebraic sets.  An algebraic set is semi-algebraic.

\begin{unproved}{proposition}\label{Tarski-Seidenberg Theorem}
\begin{inparaenum}
\item\label{TS statement}
The image of a semi-algebraic set by a polynomial map is
semi-algebraic.
\item\label{Whitney stratification}
A semi-algebraic set has a finite \naive{} stratification.
\end{inparaenum}
\end{unproved}

\eqref{Tarski-Seidenberg Theorem}(\ref{TS statement}) is due 
to Tarski and to Seidenberg.  
\eqref{Tarski-Seidenberg Theorem}(\ref{Whitney stratification}) 
is a result of \citet{Whitney1957}.

Let $U\sub \C^2$ be an open set.  A 
\moredef{holomorphic curve}%
    \index{holomorphic!curve}%
    \index{curve(s)!holomorphic} in $U$ is a $\C$-analytic
set $\G$ such that the complex manifold $\reg{\G}$ is non-empty,
everywhere of real dimension $2$, and dense in $\G$.

\begin{unproved}{proposition}\label{resolution of curve}
Let $\G\sub U$ be a holomorphic curve in an open set in $\C^2$.
There exists a complex manifold $G$ of real dimension $2$,
and a holomorphic map $\resolve\from G\to U$, such that:
\begin{inparaenum}
\item
$\resolve(G)=\G$;
\item
$\critpts{\resolve}\sub\resolve^{-1}(\sing{\G})$, and 
$\resolve\restr{\resolve^{-1}(\sing{\G})}$ has finite fibers;
\item
$\resolve\restr\resolve^{-1}(\reg{\G})$ is a holomorphic diffeomorphism.
\end{inparaenum}
\end{unproved}

The map $\resolve$ is essentially unique, and is called 
\moredef{the resolution of $\G$}%
    \index{resolution ($\resolve$)!of a holomorphic curve}%
    \index{curve(s)!holomorphic!resolution of ($\resolve$)}.  A 
\moredef{branch}%
    \index{branch!of a holomorphic curve}%
    \index{curve(s)!holomorphic!branch} of $\G$ at 
$P\in\G$ is the image by $\resolve$ of a component of
$\Int{
      \Nb
         {\{Q\}}
         {\diff
               {G}
               {\critpts
                        {\resolve}
                }
          }
}$ 
for some $Q\in\resolve^{-1}(P)$(or the germ of the image of such a component).
If $P\in\reg\G$ then there is only one branch of $\G$ at $P$, but there can
also be singular points of $\G$ at which there is only one branch of $\G$.

\begin{examples}\label{examples of curve resolutions}
Classical algebraic geometers gave names to quite a few 
special cases of branches (and resolutions).  Two examples
are of particular importance in the knot theory of complex
plane curves.
\begin{asparaenum}[\bf \theexamples.1.]
\item\label{node example}
Define $f\in\holo{\Int{D^4}}$ by $f(z,w)=z^2+w^2$.  
The holomorphic curve $\variety{f}$ has two branches at $(0,0)$.  
Its resolution is 
$\resolve\from \Int{D^2}\times\{{+},{-}\}\to\variety{f}\suchthat
(\zeta,\pm)\mapsto 2^{-1/2}(\zeta,\pm\imunit \zeta)$.  A point
$P$ of a holomorphic curve $\G\sub\Omega$ such that there exist
an open neighborhood $U$ of $P$ in $\Omega$ and 
a diffeomorphism (which may in fact be required to be holomorphic) 
$h\from(U,U\cap\Gamma,P)\to(\Int{D^4},\variety{f},(0,0))$ 
is called a \moredef{node}%
    \index{holomorphic!curve!node}%
    \index{node!of a holomorphic curve}%
    \index{curve(s)!holomorphic!node of} of $\G$.
\item\label{cusp example}
Define $f\in\holo{\Int{D^4}}$ by $f(z,w)=z^2+w^3$.  The holomorphic curve 
$\variety{f}$ has one branch at $(0,0)$.  Its resolution is
$\resolve\from \Int{D^2}\to\variety{f}\suchthat
(\zeta,\pm)\mapsto 2^{-1/2}(\zeta^3,-\zeta^2)$.  A point
$P$ of a holomorphic curve $\G\sub\Omega$ such that 
there exist an open neighborhood $U$ of $P$ in $\Omega$ and 
a diffeomorphism (which may in fact be required to be holomorphic) 
$h\from(U,U\cap\Gamma,P)\to(\Int{D^4},\variety{f},(0,0))$ 
is called a \moredef{cusp}%
    \index{holomorphic!curve!cusp}%
    \index{cusp!of a holomorphic curve}%
    \index{curve(s)!holomorphic!cusp of} of $\G$.
\end{asparaenum}
\end{examples}

A holomorphic curve $\G$ such that every point of $\sing{\G}$
is a node \resp{either a node or a cusp} is called a 
\moredef{node} \resp{\moredef{cusp}} \moredef{curve}%
    \index{node!curve}%
    \index{curve(s)!node}%
    \index{curve(s)!cusp}. %

\subsection{Configuration spaces and spaces of monic polynomials}
\label{configuration space}

Let $X$ be a topological space.  For $n\in\N$, 
the sets $\multsets{X}n$, $\subsets{n}{X}$, and $\disc{n}{X}$
are endowed with topologies by the application of 
the bijection \dispeqref{bijection from sympower to multsets} 
to the quotient topology induced on $\sympower{X}n$
from the product topology on $X^n$; with these topologies,
they are called the $n$th \moredef{multipower space}%
     \index{multipower!space}%
\index{MPnX@\protect$\multsets{\phold}n$ ($n$th multipower)!space}, %
the $n$th \moredef{configuration space}%
     \index{configuration!space}, and 
the $n$th \moredef{discriminant space}%
     \index{discriminant!space}%
\index{DeltanX@\protect$\disc{n}{\phold}$ ($n$th discriminant)!space} of $X$, 
respectively. 

If $M$ is a manifold, then clearly each equivalence class 
of the partition by type $\multsets{M}{n}/\type$ of $\multsets{M}{n}$ 
is a manifold.  However, even for connected $M$ it often happens 
that not every fiber of $\type$ is connected.
The stratification $\multsets{M}{n}/\ctype$ of $\multsets{M}{n}$ 
by the connected components of the fibers of $\type$ will
be called the \moredef{standard stratification}%
     \index{stratification!standard!of a multipower space} of 
$\multsets{M}{n}$; the standard stratification of $\multsets{M}{n}$
induces a \moredef{standard stratification}%
     \index{standard!stratification}%
     \index{stratification!standard!of a discriminant space}%
     \index{stratification!standard!of a configuration space} on 
each of $\multsets{M}{n}/\ctype$ and $\disc{n}{M}$,
since they are evidently unions of strata of $\multsets{M}{n}/\ctype$.

If $\F$ is a field, then of course $\F^n$ is an algebraic set
over $\F$, and the standard action of $\symgroup{n}$ on $\F^n$ is 
algebraic.  In general this is not enough to ensure that the set
$\sympower{\F}n$ can be endowed with as much structure as
might be desirable to full-fledged algebraic geometers.  
However, for $\F=\C$ (or any algebraically closed field), 
on general principles $\sympower{\C}n$ does have a natural 
structure as an algebraic set 
\citep[more or less naturally embedded in an affine space $\C^N$; 
see, e.g.,][]{Cartan1957} with respect to which 
the unordering map $\Unorder{\C}{n}\from\C^n\to\sympower{\C}n$ 
is a polynomial map.  By contrast, if $n>1$, then $\sympower{\R}n$ 
not an algebraic set (in any natural way), although it is semi-algebraic.

Denote by $\monicpoly{n}\isdefinedas
\{p(w)\in\C[w]\Suchthat p(w)=w^n+c_1 w^{n-1}+\dots+c_{n-1}w+c_n\}$
the \ndimensional{n} complex affine space of 
\moredef{monic polynomials}%
    \index{MPn@$\monicpoly{n}$ (monic polynomials of degree $n$)}%
    \index{monic polynomials}%
    \index{polynomial!monic} of degree $n\in\N$.
Define the 
\moredef{roots map}%
    \index{roots map ($\roots$)}%
    \index{r (roots map)@$\roots$ (roots map)} %
$\roots\from\monicpoly{n}\to\multsets{\C}{n}$ 
by $\roots(p)\isdefinedas(p^{-1}(0),\mult{p}\restr{p^{-1}(0)})$.
where $\mult{p}(z)$ is the usual
\moredef{multiplicity of $p(w)\in\monicpoly{n}$ at $z\in\C$}%
\index{multiplicity!of a polynomial $p$ at a point ($\mult{p}(\phold)$)}. Let 
$V$ be the polynomial map 
$\C^n\to\monicpoly{n}\suchthat\ontuple{z}n\mapsto(w-z_1)\dotsm(w-z_n)$.
The diagram
\begin{equation}
\begin{CD}
\R^n @<{\textstyle{\Re}}<< \C^n @>{\textstyle{V}}>> \monicpoly{n}\\
@VV{\textstyle{U_{\R,n}}}V  @VV{\textstyle{U_{\C,n}}}V 
@VV{\textstyle{\roots_{\phantom{\C}}}}V    \\
\sympower{\R}n @<{\textstyle{\Re}}<< \sympower{\C}n 
   @>{\textstyle{\iso}}>{\dispeqref{bijection from sympower to multsets}}>
   \multsets{\C}n\\ 
@VV{\textstyle{\ctype_{\mathstrut}}}V 
@VV{\textstyle{\type^{\phantom{c}}_{\mathstrut}}}V @.  \\
(\sympower{\R}n)/\ctype @. (\sympower{\C}n)/\type@. {}
\end{CD}
\end{equation}
is commutative, and defines
$\resolve\from\monicpoly{n}\to\sympower{\C}n$%
\index{resolution ($\resolve$)!of $\sympower{\C}n$}
the following results are standard; some go back, in essence, 
to \citet{Viete1593}\footnote{%
    $V$ stands for \moredef{Vi\`ete map}%
       \index{V@$V$ (Vi\`ete map)}%
       \index{Vi\`ete map ($V$)}, a coinage due (apparently) 
    to \citeauthor{Arnold1968}, now widely used.} 
and \citet{Descartes1637}.

\begin{proposition}\label{complex stratification}
The \naive{} stratification of $\sympower{\C}n$ as 
an algebraic set coincides with its stratification
by type $(\sympower{\C}n)/\type$.  In particular:
\begin{enumerate}
\item
there are no strata of odd \textup(real\textup) codimension;
\item
the only \ncodimension{0} stratum is
$\smash[t]{\reg{\sympower{\C}n}=\config{n}=\type^{-1}(\{1,\dots,1\})}$,
so $\sing{\sympower{\C}n}=\disc{n}{\C}$ is the union of 
the strata of codimension $>0$;  
\item
the only \ncodimension{2} stratum is
$\regkth{2}{\sympower{\C}n}=\type^{-1}(\{1,\dots,1,2\})$.
\hfill\qedsymbol 
\end{enumerate}
\end{proposition}

Call the stratification in \eqref{complex stratification} the
\moredef{complex stratification of $\sympower{\C}{n}$}%
    \index{stratification!complex!of $\sympower{\C}n$}%
    \index{$complex stratification$@$\phold/{\equiv_\C}$ (complex stratification)}, %
and denote it by $(\sympower{\C}{n})/{\equiv_\C}$.

\begin{unproved}{proposition}
$\resolve$ is a homeomorphism and 
$\resolve\restr\roots^{-1}(\disc{n}{\C})$ is a diffeomorphism; 
in fact, $\resolve$ is a normalization and minimal resolution 
of the algebraic set $\sympower{\C}n$.
\end{unproved}

The resolution $\resolve$ and the complex stratification 
of $\sympower{\C}{n}$ together impose a 
\moredef{complex stratification} $\monicpoly{n}/{\equiv_\C}$ on 
$\monicpoly{n}$\index{stratification!complex!of $\monicpoly{n}$}. 

\begin{unproved}{proposition}
\label{contractibility of real configuration space}
$(\sympower{\R}n)/\ctype$ is a cellulation; in particular, 
$\smash[t]{\subsets{n}{\R}}$, its unique \ncodimension{0} stratum,
is an \ncell{n}.
\end{unproved}

\begin{proposition}\label{real cellulation}
\begin{inparaenum}
\item\label{Fox-Neuwirth}
If $S$ is a stratum of $(\sympower{\C}n)/{\equiv_\C}$, then
$S/(\ctype\after\Re)$ is a stratification of $S$, and in fact
a cellulation of $S$; thus the partition $(\sympower{\C}n)/{\equiv_\R}$,
such that each ${\equiv_\R}$-class is a 
$\ctype\after\Re$\nobreakdash-class of some stratum $S$
of $(\sympower{\C}n)/{\equiv_\C}$, is a cellulation.
\item\label{semi-algebraic cells}
Each cell of $(\sympower{\C}n)/{\equiv_\R}$ is a real semi-algebraic set.
\end{inparaenum}
\end{proposition}
\begin{proof}
(\ref{Fox-Neuwirth}) is apparently originally due to 
\citet{FoxNeuwirth1962} 
\citep[see also][]{Fuchs1970,Vainstein1978,Napolitano1998}.
(\ref{semi-algebraic cells}) follows from 
\eqref{Tarski-Seidenberg Theorem}(1).
\end{proof}

Call the cellulation in \eqref{real cellulation} the
\moredef{real cellulation of $\sympower{\C}{n}$}%
    \index{cellulation!real!of $\sympower{\C}n$},
and denote it by $(\sympower{\C}{n})/{\equiv_\R}$.
The resolution $\resolve$ and the real cellulation of $\sympower{\C}n$ 
together impose a \moredef{real cellulation} $\monicpoly{n}/{\equiv_\R}$
on $\monicpoly{n}$%
    \index{cellulation!real!of $\monicpoly{n}$}%
    \index{$real cellulation$@$\phold/{\equiv_\R}$ (real cellulation)}. The
real cellulations of $\sympower{\C}{n}$ and $\monicpoly{n}$
in turn define real cellulations%
    \index{cellulation!real} of $\multsets{\C}{n}$, 
$\disc{n}{\C}$, $\resolve^{-1}(\multsets{\C}{n})$, 
and $\resolve^{-1}(\disc{n}{\C})$.

\begin{examples}\label{explicit C-strata}
For small $n$, very explicit descriptions of the 
complex stratifications are easily given.
\begin{asparaenum}[\bf \theexamples.1.]
\item\label{strata:n=1}
$\monicpoly{1}/{\equiv_\C}$ consists of a single stratum,
necessarily of \ncodimension{0}.
\item\label{strata:n=2}
$\monicpoly{2}/{\equiv_\C}$ consists of two incident strata:
$\resolve^{-1}(\config{2})$, of codimension $0$, is 
diffeomorphic to $\C\times(\diff{\C}{\{0\}})$;
$\resolve^{-1}(\disc{2}\C)$, of codimension $2$, is 
diffeomorphic to $\C$.  Explicitly, 
$\sympower{\C}2\to\C^2\suchthat%
\{w_1,w_2\}\mapsto(w_1+w_2,(w_1-w_2)^2)$ is a homeomorphism 
that maps
$\smash[b]{\config{2}}$ \resp{$\disc{2}\C$} diffeomorphically
onto $\diff{\C}{\{0\}}$ \resp{$\C\times{0}$}.
\item\label{strata:n=3}
$\monicpoly{3}/{\equiv_\C}$ consists of three mutually incident
strata: 
$\smash[t]{\config{3}}$, of codimension $0$, is diffeomorphic 
to $\C\times(\diff{\C^2}{\variety{f}})$, where 
$f(z_1,z_2)=4z_1^3+9z_2^2$ and so $\variety{f}$ is a 
cuspidal cubic curve, homeomorphic to $\C$ and having a 
single singular point; 
$\reg{\disc{3}{\C}}$, of codimension $2$, 
is diffeomorphic to $\C\times\reg{\variety{f}}$ and therefore
to $\C\times\diff{\C}{\{0\}}$; and 
$\sing{\disc{3}{\C}}$, of codimension $4$, 
is diffeomorphic to $\C\times\sing{\variety{f}}$ and therefore
to $\C$.  It is easy to write down an explicit polynomial 
homeomorphism $\sympower{\C}3\to\C^3$ giving an isomorphic 
stratification.
%
\end{asparaenum}
\end{examples}

The real cellulations of $\monicpoly{n}$, and thus 
of $\config{n}$ and $\disc{n}{\C}$, can be described
very explicitly, in all dimensions
\citep[see][]{FoxNeuwirth1962,Napolitano2000}. 
For the purposes of this survey, it is sufficient to 
describe the cells of dimension $2n$, $2n-1$, and $2n-2$ only,
along with their incidence relations; this can be done in a uniform 
manner for all $n$.  

\begin{example}\label{low-codimension cells of MPn}  
In $\smash[t]{\config{n}/{\equiv_\R}}$ there is exactly one
cell of dimension $2n$, exactly $n-1$ cells of dimension $2n-1$, 
and exactly $(n-1)(n-2)/2$ cells of dimension $2n-2$.  
\begin{enumerate*}
\item
The cell $C_0$ of dimension $2n$ consists of all $\untuple{z}{n}$ with
$\Re z_1<\dots<\Re z_n$.  
\item
For $k=1,\dots,n-1$, there is a cell $C_k$
of dimension $2n-1$ consisting of all $\untuple{z}{n}$ with
$\Re z_1<\dots<\Re z_k=\Re z_{k+1}<\dots <\Re z_n$ and 
$\Im z_k\ne\Im z_{k+1}$; $C_k$ is transversely oriented
by the complex orientation of $\monicpoly{n}$, and 
$C_0$ is incident on $C_k$ from both sides---more
precisely, there is a simple closed curve in $\smash{\config{n}}$
that intersects $C_k$ in a single point, transversely, and is
otherwise contained in $C_0$.
\item
For $1\le k\le n-2$, there is a cell $C_{k,k+1}$ of dimension $2n-2$ 
consisting of all $\untuple{z}{n}$ with
$\Re z_1<\dots<\Re z_k=\Re z_{k+1}=\Re z_{k+2}<\dots <\Re z_n$ and 
$\card\{\Im z_k,\Im z_{k+1},\Im z_{k+2}\}=3$.  The two cells $C_k$
and $C_{k+1}$ of dimension $2n-1$ are each triply incident on 
$C_{k,k+1}$---more precisely, the stratification induced on a
small \ndisk{2} in $\monicpoly{n}$ intersecting $C_{k,k+1}$
in a single point, transversely, is as pictured to the left
of \Figref{incidence figure}.
\item
For $1\le i<j-1\le n-2$, there is a cell $C_{i,j}$ of dimension $2n-2$ 
consisting of all $\untuple{z}{n}$ with
$\Re z_1<\dots<\Re z_i=\Re z_{i+1}<\dots<\Re z_j=\Re z_{j+1}<\dots<\Re z_n$,
$\Im z_i\ne\Im z_{i+1}$, and $\Im z_j\ne\Im z_{j+1}$.
The two cells $C_i$
and $C_j$ of dimension $2n-1$ are each doubly incident on 
$C_{i,j}$---more precisely, the stratification induced on a
small \ndisk{2} in $\monicpoly{n}$ intersecting $C_{i,j}$
in a single point, transversely, is as pictured in the middle 
of \Figref{incidence figure}.
\end{enumerate*}
\begin{figure}
\centering
\includegraphics[width=.9\textwidth]{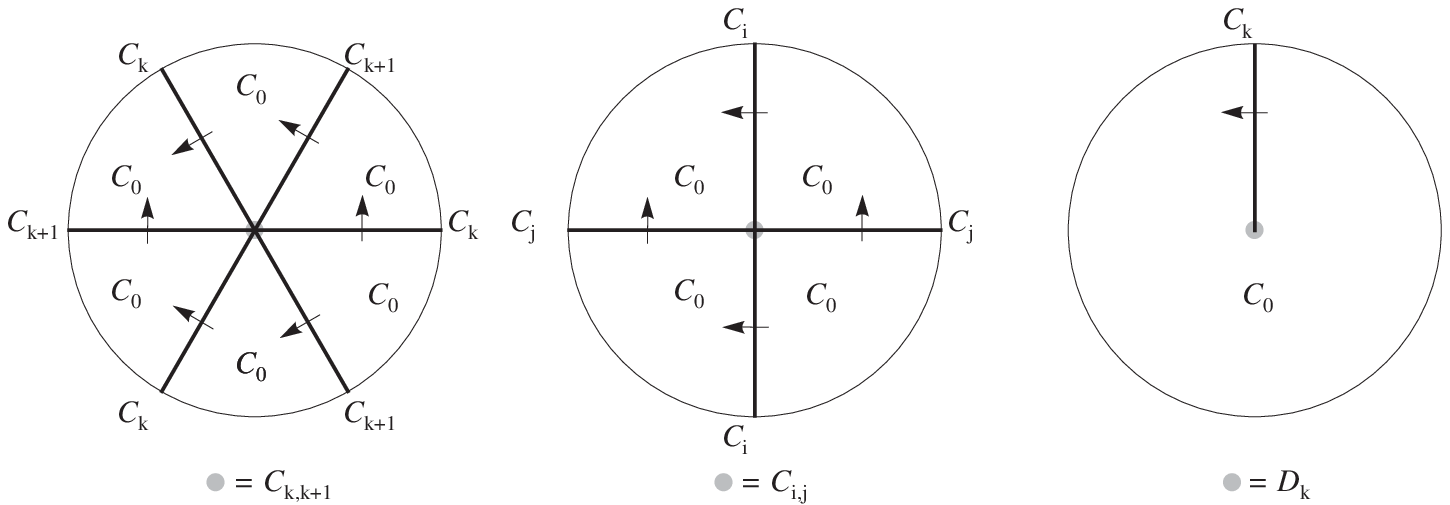}
\caption{The transverse structure of $\monicpoly{n}/{\equiv_\R}$
along its codimension-$2$ cells.\label{incidence figure}}
\end{figure}
In $\disc{n}{\C}/{\equiv_\R}$, there are no cells
of dimension $2n$ or $2n-1$.  For $1\le k\le n-2$, there 
is one cell $D_k$ of dimension $2n-2$ 
consisting of all $\untuple{z}{n}$ with
$\Re z_1<\dots<\Re z_k=\Re z_{k+1}<\dots <\Re z_n$ and 
$\Im z_k=\Im z_{k+1}$.  The stratification induced on a
small \ndisk{2} in $\monicpoly{n}$ intersecting $D_k$
transversely at a single point is as pictured at the right 
of \Figref{incidence figure}.
\end{example}

\subsection{Contact 3-manifolds, Stein domains, and Stein surfaces}
\label{contact, Stein}

This section simply resumes basic definitions and needed results.
For details on contact structures and contact \nmanifold{3}s, 
see \citet{Etnyre2004}.  
For the topology of Stein domains and Stein surfaces, 
see \citet{Gompf1998}.  
For the complex function theory of Stein domains,
Stein manifolds, and Stein spaces (in general dimensions), see
\citet{GunningRossi1965}.

\begin{definitions}\label{contact definitions}
Let $M$ be a \nmanifold{3} without boundary.  
A \moredef{contact structure}%
    \index{contact!structure} on $M$ is a completely
non-integrable field $\xi$ of tangent $2$-planes on $M$;  
for instance, on any round \nsphere{3} $S^3\sub\C^2$, 
the field $\xi$ of tangent $2$-planes that are actually complex 
lines is a contact structure, called the \moredef{standard} 
contact structure%
    \index{standard!contact structure}%
    \index{contact!structure!standard} on that sphere.
A \nmanifold{3} with a contact structure is called a 
\moredef{contact manifold}%
    \index{contact!manifold}%
    \index{manifold!contact}.  A 
(closed) \nsubmanifold{1} $L$ in a contact manifold 
$M$ with contact structure $\xi$ is \moredef{Legendrian}%
    \index{Legendrian!$1$-manifold} in case 
$\tanspace{x}{L}\sub\xi_x\sub\tanspace{x}{M}$ for all $x\in L$.
Of course $L$ is Legendrian if and only if $-L$ is.\footnote{In 
   particular, in the context of links in $S^3$ (equipped with 
   its standard contact structure) it is common practice to refer 
   to either $L$ or $\unorient{L}$ as a \moredef{Legendrian link}%
    \index{link(s)!Legendrian}%
    \index{Legendrian!link} (or knot, as the case may be), 
   and---contrary to the conventions established earlier, which
   require that a link (or knot) be oriented---this practice will 
   be followed here.}  In any contact manifold $M$, a Legendrian 
\nsubmanifold{1} $L$ is naturally endowed with a normal line field 
$\xi^\perp_L$ (unique up to isotopy), which determines an annular 
surface $A\sub M$ (also unique up to isotopy) 
containing $L$ as a retract and such that 
$\normbun{\inclusion{L}{A}}=\xi^\perp_L$.
In particular, a Legendrian link $L$ in $S^3$ with its standard 
contact structure has a \moredef{natural framing}%
    \index{natural!framing, of a Legendrian link ($\Legf{\phold}$)}%
    \index{framing!natural, of a Legendrian link ($\Legf{\phold}$)}%
    \index{Legendrian!link!natural framing ($\Legf{\phold}$)} $\Legf{L}$ 
for which $\AKn{L}{\Legf{L}}$ is such an annular surface.  The 
\moredef{Thurston--Bennequin number}%
    \index{Thurston--Bennequin number ($\tb(\phold)$)} of
a Legendrian knot $K\sub S^3$ is $\tb(K)\isdefinedas \Legf{K}(K)$.
For an arbitrary knot $K\sub S^3$, denote by $\TB(K)$
the \moredef{maximal Thurston--Bennequin number}%
    \index{maximal Thurston--Bennequin number ($\TB(\phold)$)}%
    \index{Thurston--Bennequin number ($\tb(\phold)$)!maximal ($\TB(\phold)$)} %
$\max\{\tb(K')\Suchthat \text{$K'$ is a Legendrian knot isotopic to $K$}\}$ 
of $K$; $\TB(K)$ is an integer \citep[i.e., neither $-\infty$ nor 
$\infty$; see][]{Bennequin1983}.
\end{definitions}

\begin{definitions}\label{Stein definitions}
An \moredef{open Stein manifold}%
    \index{manifold!Stein}%
    \index{manifold!Stein!open}%
    \index{Stein!manifold} is a complex manifold that is
holomorphically diffeomorphic to a topologically closed complex
submanifold of some complex affine space $\C^N$ 
\citep[equivalently, to a non-singular global analytic 
set $\variety{f}$ with $f\in\holo{\C^N}$; see][]{GunningRossi1965}.
For instance, $\C^N$ itself is an open Stein manifold, as is
(non-obviously) any open subset of $\C$.
An \moredef{open Stein surface}%
    \index{surface!Stein!open}%
    \index{open!Stein surface}%
    \index{Stein!surface!open} is an open Stein manifold of
real dimension $4$.
Let $M$ be an open Stein surface.
An \bydef{exhausting strictly plurisubharmonic function}%
    \index{function!exhausting strictly plurisubharmonic} on 
a non-empty open set $U$ of $M$ is a smooth function 
$\rho\from U\to\R$ that is bounded below, 
\moredef{proper}%
    \index{proper!function}%
    \index{function!proper} (in the sense that  
$\rho^{-1}(\clint{a}{b})$ is compact for all $a,b\in\R$),
and such that for each $c\in\R$, the field of tangent complex 
lines on the \nmanifold{3} $\diff{\rho^{-1}(c)}{\critpts{\rho}}$ 
is a contact structure,
called the \moredef{natural}%
    \index{natural!contact structure}%
    \index{contact!structure!natural} contact structure on 
that \nmanifold{3}.  (Thus, the standard structure on $S^3$
is the natural structure for the standard embedding $S^3\embto\C^2$.)
A \moredef{Stein domain}%
    \index{Stein!domain}%
    \index{domain, Stein} in $M$ is a compact \ncodimension{0}
submanifold $X\sub M$ such that 
$X$ is a sublevel set $\rho^{-1}(\Rle{c})$
($c\in\diff{\R}{\critvals{\rho}}$) of an 
exhausting strictly plurisubharmonic Morse function 
$\rho\from U\to\R$ on an open set $U\sub M$.
A \moredef{Stein surface with boundary}%
  \index{Stein!surface!with boundary}%
  \index{surface!Stein!with boundary} is a compact \nmanifold{4} $X$ 
with a complex structure on $\Int{Y}$ such that $X$ is diffeomorphic
to a Stein domain in some open Stein surface by a diffeomorphism that
is holomorphic on $\Int{X}$.
A \moredef{Stein disk}%
    \index{Stein!disk}%
    \index{disk!Stein} is a Stein domain $D$ in $\C^2$ 
diffeomorphic to $D^4$ (for instance, $D^4$ itself, 
where $\rho(z,w)$ can be taken 
to be $\norm{z,w}^2$).  Any non-singular level \nmanifold{3} of 
an exhausting strictly plurisubharmonic function on a Stein surface
is called a \moredef{strictly pseudoconvex}%
    \index{manifold!strictly pseudoconvex}%
    \index{strictly!pseudoconvex $3$-manifold} \nmanifold{3}.
A closed contact manifold is \moredef{Stein-fillable}%
    \index{Stein!-fillable $3$-manifold}%
    \index{manifold!Stein-fillable} in case it is
diffeomorphic to a strictly pseudoconvex \nmanifold{3} $N$ by
a diffeomorphism carrying its contact structure to the natural
contact structure on $N$.
\end{definitions}

Let $X\sub\C^2$ be a Stein domain.  It is convenient 
to establish notation for several subsets of 
$\cont{X}\isdefinedas\{f\from X\to\C\Suchthat \text{$f$ is continuous}\}$;
although $\cont{X}$, equipped with its $\sup$ norm, is 
a well known Banach algebra, here it (and its subsets) will
not be endowed with any topology.  
\begin{equation}\label{function algebras}
\begin{aligned}
\intholo{X}\isdefinedas \{f\in\cont{X}\Suchthat &
\text{$f\restr{\Int{X}}$ is holomorphic}\} 
\\	
\holo{X}\isdefinedas \{f\in\cont{X}\Suchthat &
f=F\restr X \text{ for some open neighborhood }
U \text{ of} 
\\ 
& \text{$X$ in $\C^2$ and some holomorphic $F\from U\to\C$}\}
\\
\nzholo{X}\isdefinedas 
\{f\in\holo{X}\Suchthat &
0\notin f(X)\}
\\
\nizholo{X}\isdefinedas 
\{f\in\holo{X}\Suchthat &
0\notin f(\Int{X})\}.
\end{aligned}
\end{equation}
An element of $\intholo{X}$ is called a \bydef{germ} of a holomorphic 
function on $X$.  Both $\intholo{X}$ and $\holo{X}$ are algebras. 
By a standard argument, $\nzholo{X}$ is the group of units 
of $\holo{X}$.  Clearly $\nizholo{X}$ is a multiplicative subsemigroup 
of $\holo{X}$ properly containing $\nzholo{X}$.

Major reasons for complex analysts' interest in Stein manifolds
include general theorems of which the following are special cases.

\begin{unproved}{theorem}\label{Cousin problem}
If $X\sub \C^2$ is a Stein disk and $\G\sub U$ is a holomorphic 
curve in an open neighborhood $U$ of $X$, then 
$\G\cap X=f^{-1}(0)$ for some $f\in\holo{X}$.
\end{unproved}

\begin{unproved}{theorem}\label{polynomial convexity}
Any holomorphic function on a Stein disk $X$ can be 
arbitrarily closely uniformly approximated, along with
any finite number of its derivatives, by the restriction 
to $X$ of a polynomial function.
\end{unproved}

The following results are especially useful for topological
applications.

\begin{unproved}{theorem}\label{ribbon property of cpx curves}
If $X\sub \C^2$ is a Stein disk with exhausting plurisubharmonic
Morse function $\rho$, and $f\in\holo{X}$ is such that
$\sing{\variety{f}}=\emptyset$, then 
$\variety{f}$ is $\rho$-ribbon.
\end{unproved}

\begin{unproved}{theorem}\label{covering of Stein is Stein}
Every covering space of an open Stein manifold is an open
Stein manifold.  A finite-sheeted branched 
covering space of a Stein disk branched along a non-singular
holomorphic curve is a Stein surface with boundary.
\end{unproved}
%
%

\section{Braids and braided surfaces}
\label{braids and braided surfaces}

Much of the following material is treated (usually more
generally and often from a different perspective) by 
\citet{Birman1974} and \citet{BirmanBrendle2004}, references to
which should be assumed throughout.

\subsection{Braid groups}\label{braid groups}

\begin{definition}\label{braid group definition}
For any $n\in\Nstrict$ and $X\in \smash[t]{\config{n}}$, let 
$\smash[t]{\brgp{X}\isdefinedas\fundgp{\config{n}}{X}}$%
    \index{B@$\brgp{\phold}$ (braid group)!$\brgp{n}$ (standard, on $n$ strings)}%
    \index{B@$\brgp{\phold}$ (braid group)!$\brgp{X}$ (based at \protect$X\in\smash[t]{\config{n}}$)}, %
and call $\brgp{X}$ \moredef{an \nstring{n}} 
\moredefdanger{braid group}%
    \index{braid!group ($\brgp{\phold}$)}%
    \index{braid!group ($\brgp{\phold}$)!n-string@\nstring{n}}.\footnote{%
    As \citet{Magnus1974,Magnus1976} points out, 
    in effect $\fundgpnbp{\config{n}}$
    was investigated, and recognized as a ``braid group'', by 
    \citeauthor{Hurwitz1891} as early as \citeyear{Hurwitz1891}.   
    Apparently this insight had been long forgotten when 
    \citet[p. 119]{FoxNeuwirth1962} described as ``previously 
    unnoted'' their ``remark that $\brgp{n}$ may be considered as the 
    fundamental group of the space~\dots~of configurations of $n$ 
    undifferentiated points in the plane.''  It is interesting to 
    speculate as to possible reasons for this instance of what
    \citeauthor{Epple1995} calls ``elimination of contexts'' 
    \citep[see especially][p. 386, n. 32]{Epple1995}.} The 
\moredef{standard}%
    \index{standard!n-string@\nstring{n} braid group} %
\nstring{n} braid group is $\brgp{n}\isdefinedas \brgp{\intsto{n}}$.
By convention, the (unique) \nstring{0} braid group 
is $\brgp{0}\isdefinedas \{\nbrid0\}$.  Write $o_{X}$%
    \index{oX@$o_{X}$ (identity of $\brgp{X}$)} for 
the identity of $\brgp{X}$, and let $\nbrid{n}\isdefinedas 
o_{\intsto{n}}$.  
\end{definition}

Of course, since $\smash[t]{\config{n}}$ is connected, 
every \nstring{n} braid group is isomorphic to $\brgp{n}$, 
but it is very convenient to allow more general basepoints.  
With the conventions in \eqref{braid group definition}, 
$\brgp{0}$ is isomorphic but not identical to 
$\brgp{1}=\{\nbrid{1}\}$; 
this is consistent with the obvious fact that the groups 
$\brgp{n}$, $n\in\Nstrict$, being fundamental groups of pairwise 
distinct spaces, are pairwise disjoint.

\begin{theorem}\label{presentation theorem}
\begin{equation}\label{standard presentation of Bn}
\smash{%
\quad
\brgp{n}=
   \gp{
      \brgen{s}, s\in\{1,\dots,n-1\}
   }
   {%
   \genfrac{}{}{0pt}{}
   {\commute{\brgen{s}}{\brgen{t}}~(|s-t|>1),}
   {\YBax{\brgen{s}}{\brgen{t}}~(|s-t|=1)}
   }
}
\qquad\quad\qedsymbol
\end{equation}
\end{theorem}

It is usual to call 
\dispeqref{standard presentation of Bn} the
\moredef{standard presentation}%
    \index{standard!presentation of $\brgp{n}$}%
    \index{presentation!standard, of $\brgp{n}$} of $\brgp{n}$,
the generators $\brgen{s}$ of 
\dispeqref{standard presentation of Bn} the
\moredef{standard generators}%
    \index{standard!generator(s) of $\brgp{n}$ ($\brgen{i}$)}%
    \index{generator(s)!standard, of $\brgp{n}$ ($\brgen{i}$)} of $\brgp{n}$,
and the relators of \dispeqref{standard presentation of Bn} the
\moredef{standard relators}%
    \index{standard!relator(s) of $\brgp{n}$)}%
    \index{relator(s)!standard, of $\brgp{n}$} of $\brgp{n}$.
\citep[It is also usual, and perhaps regrettable, to conflate 
$\brgen{s}\in\brgp{n}$ with $\brgen{s}\in\brgp{n'}$ for all 
$n,n'>s$. A more precise notation was proposed by]%
[see \protect\eqref{canonical injections of standard braid groups}.]%
{Rudolph1985} The 
detailed proof of \eqref{presentation theorem} 
given by \citet{FoxNeuwirth1962} is, more or less exactly,
an application of the usual algorithm 
\citep*[as in][]{Magnusetal1976}, which produces a
presentation of the fundamental group of a 
\ndimensional{2} cell complex with one \ncell{0}
(the basepoint) having a generator for each \ncell{1}
and a relator for each \ncell{2}, to the $2$-skeleton 
of the cellulation of $\config{n}$ that is dual to 
$\config{n}/{\equiv_\R}$.  

\subsection{Geometric braids and closed braids}
Let $I\sub\R$ be a closed interval.
Let $p\from I\to\smash[t]{\config{n}}$ be a closed path.  
The multigraph $\gr{p}\sub I\times\C$ of $p$ is called a 
\nstring{n} \moredef{geometric braid}%
    \index{geometric braid}%
    \index{braid!geometric} for the (algebraic) braid
in $\brgp{p(\Bd I)}$ represented by $p$.

Let $q\from\hunknot\to\smash[t]{\config{n}}$ be a loop
with domain the horizontal unknot $\hunknot\sub S^3$.
The multigraph $\gr{q}\sub S^1\times\C$ of $q$ is a
simple collection of closed curves in the open solid torus
$\hunknot\times\C$.  Using the vertical unbook $\ob{o}$,
it is easy to create a (nearly standard) identification
of $\hunknot\times\C$ with $\diff{S^3}{\vunknot}$,
and thus a (nearly standard) embedding of $\gr{q}$
into $\diff{S^3}{\vunknot}$.  This embedding has the property 
that $\arg{\ob{o}}\restr\gr{q}\from\gr{q}\to S^1$ is a covering 
map of degree $n$, and every simple collection of closed curves 
$L\sub\hunknot\times\C$ such that 
$\arg{\ob{o}}\restr L\from L\to S^1$ is a covering map
of degree $n$ arises in this way from some $q$.  Call any 
such $L$ an \nstring{n} \moredef{closed $\ob{o}$-braid}%
    \index{o-braid@$\ob{o}$-braid, closed}%
    \index{closed!braid}%
    \index{braid!closed} \citep{Rudolph1988}.  In general,
if $\ob{a}$ is an \bydef{unbook} (that is, an open book on $S^3$ 
with unknotted binding $A$), then $\ob{a}$ is equivalent to $\ob{o}$, 
and the class of \moredef{closed $\ob{a}$-braids} (sometimes called, 
slightly abusively, simply ``closed braids with axis $A$'') is defined 
by any such equivalence.  Given a basepoint $\basept\in \hunknot$, 
$q$ naturally represents an element of $\beta_q\in\brgp{q(\basept)}$.  
The notation $\closedbraid{\beta_q}$ is often used for 
the closed $\ob{o}$-braid $\gr{q}$.

\subsection{Bands and espaliers}\label{bands and espaliers}

For $n=2$, \dispeqref{standard presentation of Bn} says that 
$\brgp{2}$ is infinite cyclic.  More specifically and directly,
\eqref{explicit C-strata}(\ref{strata:n=2}) shows that, 
for any edge $\eedge\sub\C$, the \nstring{2} braid group 
$\brgp{\Bd\eedge}$ is infinite cyclic with a preferred generator, 
say $\brgen\eedge$; in fact $\brgen\eedge$ depends only on $\Bd\eedge$.
In particular, any two \nstring{2} braid groups are canonically 
isomorphic.  For $n>2$, typically no isomorphism between distinct 
\nstring{n} braid groups has much claim to be called canonical;
however, the following is an immediate consequence of 
\eqref{contractibility of real configuration space}.

\begin{unproved}{proposition}
\label{canonical isomorphism for real basepoints}
If $\smash[t]{X\in\subsets{n}{\R}}$, then $\brgp{X}\iso\brgp{n}$,
and any path in $\smash[t]{\subsets{n}{\R}}$ from $X$ to 
$\intsto{n}$ induces this canonical isomorphism. 
\end{unproved}

For $n\ge 2$, if $X\cap\eedge=\Bd\eedge$, 
then there is a \moredef{natural injection}%
    \index{injection!natural, of braid groups ($\brinj{\phold}$)}%
    \index{natural!injection, of braid groups ($\brinj{\phold}$)}%
    \index{iota@$\brinj{\phold}$ (natural or canonical injection)} %
$\brinj{\eedge;X}\from\brgp{\Bd\eedge}\to\brgp{X}$,
and $\brinj{\Bd\eedge;X}$ depends 
only on the isotopy class of $\eedge$ (rel.\ $\Bd\eedge$) 
in $\diff{\C}{(\diff{X}{\eedge})}$.

\begin{definitions}
A \moredef{positive $X$-band}\label{def:positive X-band}%
    \index{positive!X-band@$X$-band}%
    \index{X-band@$X$-band!positive}%
    \index{band!positive} is any 
$\brgen{\eedge;X}\isdefinedas \brinj{\eedge;X}(\brgen\eedge)\in \brgp{X}$.
(When $X$ is understood, or irrelevant, $\brgen{\eedge;X}$ may
be abusively abbreviated to $\brgen\eedge$.)
A \moredef{negative $X$-band}%
    \index{negative!X-band@$X$-band}%
    \index{X-band@$X$-band!negative}%
    \index{band!negative} is the 
inverse of a positive \nband{X}.  An \moredef{$X$-band}%
    \index{X-band@$X$-band}%
    \index{band} is a positive or negative \nband{X}.
Let $\absband{\brgen{\eedge;X}^{\pm 1}}\isdefinedas \brgen{\eedge;X}$%
    \index{$absband$@$\absband{\phold}$ (``absolute value'' of a band)} denote
the \moredef{absolute value}%
    \index{absolute value of a band} of 
the band $\brgen{\eedge;X}^{\pm 1}$,
$\sign{\brgen{\eedge;X}^{\pm 1}}\isdefinedas \pm$%
    \index{epsilon@$\sign{\phold}$ (sign)} the
\moredef{sign}%
    \index{sign ($\sign{\phold}$)!of a band} (positive or negative)
of $\brgen{\eedge;X}^{\pm 1}$,
and $\edgeband{\brgen{\eedge;X}^{\pm 1}}$%
    \index{e@$\edgeband{\phold}$ (edge-class of a band)} the 
\moredef{edge-class}%
    \index{edge!-class of a band ($\edgeband{\phold}$)} of 
$\brgen{\eedge;X}^{\pm 1}$, that is,
the isotopy class 
(rel.\ $\Bd\eedge$) of $\eedge$ in $\diff{\C}{(\diff{X}{\eedge})}$.  An
\moredef{\nband{X}word}%
    \index{X-bandword@$X$-bandword}%
    \index{bandword} of \moredef{length}%
    \index{length (of a bandword)} $k$ is a \ntuple{k} 
$\brep{b}\defines\brepvec{b}{k}$ such that each $b(i)$
is an \nband{X}.\footnote{%
    In particular, \nband{\intsto{n}}words are just 
    what (since \citeyear{Rudolph1983a}) I have long called
    ``band representations''---a coinage which I would like to 
    suppress, due to its misleading suggestion of a relation to 
    group representation theory.} An 
\nband{X}word $\brep{b}$ is \moredefdanger{quasipositive}%
    \index{quasipositive!bandword}%
    \index{X-bandword@$X$-bandword!quasipositive}%
    \index{bandword!quasipositive} in case each $b(i)$ is
a positive \nband{X}.  The \moredef{braid of}%
    \index{braid!of a bandword ($\braidof{\phold}$)}%
    \index{bandword!braid of ($\braidof{\phold}$)}%
    \index{beta@$\braidof{\phold}$ (braid of a bandword)} $\brep{b}$
is $\braidof{\brep{b}}\isdefinedas b(1)\dotsm b(k)\in\brgp{X}$.
Every braid in $\brgp{X}$ is the braid of some \nband{X}word.
A braid in $\brgp{X}$ is \moredef{quasipositive}%
    \index{quasipositive!braid}%
    \index{braid!quasipositive} in case it is the braid
of a quasipositive \nband{X}word.
\end{definitions}

\citet{Orevkov2004} gives an algorithm to determine 
whether or not a given braid in $\brgp{3}$ is quasipositive.
\citet{Bentalha2004} generalizes this to all $\brgp{n}$.

\begin{unproved}{proposition}
Any two positive \nband{X}s are conjugate in $\brgp{X}$.
\end{unproved}

An \nstring{n} braid group is the knotgroup of $\disc{n}{\C}$ 
in $\subsets{n}{\C}$ \citep[see][]{Rudolph1983b}; a positive 
band is a meridian.

\begin{proposition}\label{basic band relations}
Let $\smash[t]{X\in \config{n}}$, $n\ge 2$.  
Let $\eedge, \fedge\sub\C$ be two edges with 
$X\cap\eedge=\Bd\eedge, X\cap\fedge =\Bd\fedge$.  
If $\eedge\cap\fedge=\emptyset$ 
\resp{$\eedge\cap\fedge=\{x\}\sub \Bd\eedge\cap\Bd\fedge$},
then 
$\commute{\brgen{\eedge;X}}{\brgen{\fedge;X}}=o_X$ 
\resp{$\YBax{\brgen{\eedge;X}}{\brgen{\fedge;X}}=o_X$}.
\end{proposition}

\begin{proof}
This is readily proved directly.  Alternatively, note
that: 
\begin{inparaenum}
\item
the general case is isotopic to a special case 
in which $X=\intsto{n}$, $\eedge=\clint{1}{2}$, and 
$\fedge=\clint{s}{s+1}$ with $s\in\{2,\dots,n\}$; 
\item
in such a special case, $\brgp{X}=\brgp{n}$,
$\brgen{\eedge;X}=\brgen{1}$ and $\brgen{\fedge;X}=\brgen{s}$ 
are two standard generators 
of \dispeqref{standard presentation of Bn}, 
and the claimed commutator and yangbaxter relators
are two standard relators 
of \dispeqref{standard presentation of Bn}.
\end{inparaenum}
\end{proof}

\begin{unproved}{proposition}
\label{group from tree} Let $\Ts\sub\C$ be a planar tree.
The positive \nband{\Verts{\Ts}}s $\brgen{\eedge;\Verts{\Ts}}$, 
$\eedge\in\Edges{\Ts}$, generate 
$\brgp{\Verts{\Ts}}$; no proper subset of them does so.
If $\eedge\cap\fedge=\emptyset$ 
\resp{$\eedge\cap\fedge=\{z\}, z\in\Verts{\Ts}$}, then 
$\commute{\brgen{\eedge;\Verts{\Ts}}}{\brgen{\fedge;\Verts{\Ts}}}=%
o_{\Verts{\Ts}}$ %
\resp{$\YBax{\brgen{\eedge;\Verts{\Ts}}}{\brgen{\fedge;\Verts{\Ts}}}=%
o_{\Verts{\Ts}}$}.
\end{unproved}

Call the $\Verts{\Ts}$-bands $\brgen{\eedge;\Verts{\Ts}}$, 
$\eedge\in\Edges{\Ts}$, 
the $\Ts$-\moredef{generators}%
    \index{generator(s)!$\Ts$- (of $\brgp{\Verts{\Ts}}$)}%
    \index{T-generators@$\Ts$-generators (of $\brgp{\Verts{\Ts}}$)} of 
$\brgp{\Verts{\Ts}}$.  
Note that \eqref{group from tree} asserts that the 
braid group $\brgp{\Verts{\Ts}}$ is a quotient of 
\begin{equation}\label{T-presentation of non-braid group}
\gp{
      \brgen{\eedge}, \eedge\in\Edges{\Ts}
   }
   {%
   \genfrac{}{}{0pt}{}
   {\commute{\brgen{\eedge}}
              {\brgen{\fedge}}~(\card{\eedge\cap\fedge}=0),}
   {\YBax{\brgen{\eedge}}
           {\brgen{\fedge}}~(\card{\eedge\cap\fedge}=1)}
   }
\end{equation}
but not that these groups are identical.
In fact they are easily seen to be so if and only if
$\Ts$ has no intrinsic vertices.

An \bydef{espalier} is a planar tree $\Ts$ 
such that each $\eedge\in\Edges{\Ts}$ is a proper edge in $\LHP$ 
and $\Re\restr({\diff{\eedge}{\collar{\Bd\eedge}{\eedge}}})$ is 
injective.  \citet{Rudolph2001a} gives proofs of the following facts 
about espaliers.

\begin{unproved}{proposition}\label{espaliers}
\begin{inparaenum}%
\item\label{every tree can be espaliered}
Every planar tree is isotopic, in $\C$, to an espalier.
\item\label{combinatorics determine espalier}
The embedding of an espalier $\Ts$ in $\LHP$ is 
determined, up to isotopy of $(\Ts,\Verts{\Ts})$ 
in $(\LHP,\R)$, by the combinatorial structure 
of its cellulation together with the order induced 
on $\Verts{\Ts}$ by its embedding in $\R=\Bd{\LHP}$. 
In particular, given 
$X=\utupsep{x}{1}{n}{<}{\dotsm}\sub\R$
and $n-1$ pairs $\{x_{i(p)}<x_{j(p)}\}\sub X$,
the following are equivalent.
\begin{inparaenum}
\item
There is an espalier $\Ts$ with 
$\Verts{\Ts}=\{x_1,\dots,x_n\}$, 
$\Edges{\Ts}=\{\eedge_1,\dots,\eedge_{n-1}\}$, and 
$\Bd\eedge_p=\{x_{i(p)},x_{j(p)}\}$.
\item
For $1\le p<q\le n-1$ the
pairs $\{x_{i(p)}<x_{j(p)}\}$ and $\{x_{i(q)}<x_{j(q)}\}$ are not
linked \textup(i.e., they are either in touch or unlinked\textup).
\end{inparaenum}
\end{inparaenum}
\end{unproved}

\begin{definition}
Let $X\in \smash[t]{\subsets{n}{\R}}$.  An \nband{X} 
$\brgen{}$ is \moredef{embedded}%
    \index{embedded!$X$-band}%
    \index{X-band@$X$-band!embedded}%
    \index{band!embedded} provided that some 
$\eedge\in\edgeband{\brgen{}}$ is a proper edge in $\LHP$;
when it is given that $\brgen{}=\brgen{\eedge;X}^{\pm1}$ is 
embedded, it is assumed without further comment that 
$\eedge\in\edgeband{\brgen{}}$ is such an edge.
Clearly there are exactly $n(n-1)/2$ embedded positive 
\nband{X}s, one---which may be denoted $\brgen{x_i,x_j;X}$---for 
each $\{x_i<x_j\}\sub X$.  The 
embedded positive \nband{X}s generate $\brgp{X}$.
\citep[A presentation \protect\dispeqref{presentation of a group}
of $\brgp{n}$ with the positive \nband{\intsto{n}}s as generators 
is given by][and many interesting conclusions drawn therefrom.]%
{BirmanKoLee1998}  An \nband{X}word $\brep{b}$ is \moredef{embedded}%
    \index{embedded!$X$-bandword}%
    \index{X-bandword@$X$-bandword!embedded}%
    \index{bandword!embedded} in case each $b(i)$ is embedded.
An \nband{X}word $\brep{b}$ is called \moredef{positive}%
    \index{bandword!positive}%
    \index{positive!bandword@\nband{X}word}
in case it is both quasipositive and embedded.
\end{definition}

Let $X_0\sub X\in\smash[t]{\config{n}}$.
Unless $n_0\isdefinedas\card{X_0}$ equals $n$, 
typically no non-trivial homomorphism $\brgp{X_0}\to\brgp{X}$ 
has much claim to be called canonical.  For $X\sub\R$, however,
the situation is much better (cf.\ \eqref{canonical 
isomorphism for real basepoints}).  In fact, given 
a positive embedded \nband{X_0} $\brgen{x_i,x_j;X_0}\in\brgp{X_0}$,
let $\brinj{X_0;X}(\brgen{x_i,x_j;X_0})\isdefinedas 
\brgen{x_i,x_j;X}\in\brgp{X}$.

\begin{unproved}{proposition}
\label{canonical injection exists}
There is a unique homomorphism $\brinj{X_0;X}\from\brgp{X_0}\to\brgp{X}$
extending $\brinj{X_0;X}$ as defined on the embedded \nband{X_0}s, and 
it is injective.
\end{unproved}

Call $\brinj{X_0;X}$ the \moredef{canonical injection}
    \index{injection!canonical, of braid groups ($\brinj{\phold}$)}%
    \index{injection!canonical, of braid groups ($\brinj{\phold}$)}%
    \index{canonical!injection, of braid groups ($\brinj{\phold}$)}%
    \index{iota@$\brinj{\phold}$ (natural or canonical injection)} of
$\brgp{X_0}$ into $\brgp{X}$.  (The collision of notation with
$\brinj{\eedge;X}$ is unproblematic; if $X\sub\R$ and $\eedge$ is 
embedded with $\Bd{\eedge}=\{x_i,x_j\}$, 
then $\brinj{x_i,x_j;X}=\brinj{\Bd\eedge;X}$.)
\eqref{canonical injection exists} excuses the conflation,
under the single name $\brgen{x_i,x_j}$, of
all the positive embedded bands $\brgen{x_i,x_j;X}$
for ${\x_i<\x_j}\sub X\sub \R$, $X$ finite.

\begin{example}
\label{canonical injections of standard braid groups}
For $n\ge m$, the canonical injection $\brinj{\intsto{m};\intsto{n}}$ 
is implicit in the identification of $\brgp{m}$ as a subgroup 
of $\brgp{n}$, discussed following \eqref{presentation theorem}.
\citet{Rudolph1985} proposed the notation $\beta^{(n-m)}$
for $\brinj{\intsto{m};\intsto{n}}(\beta)$, and the 
convention (extending the notations $\nbrid0$ and
$\nbrid{n}$ introduced in \eqref{braid group definition})
that for each $m\in\N$, $\brgen{m-1}$ denote only an element
of $\brgp{m}$, the other standard generators of $\brgp{m}$
being $\brgen{1}^{(m-2)},\dots,\brgen{m-2}^{(1)}$.
This notation and convention have been widely unadopted.
\end{example}

\begin{definition}
Let $\Ts$ be an espalier.  A \moredef{\nband{\Ts}word} is
a \nband{\Verts{\Ts}}word $\brep{b}$ such that every $\absband{b(i)}$
is a $\Ts$\nobreakdash-generator (so, in particular, a
\nband{\Ts}word is an embedded \nband{\Verts{\Ts}}word).  
A positive \nband{\Ts}word $\brep{b}$ is called 
\moredef{strictly $\Ts$-positive}%
    \index{bandword!positive!strictly}%
    \index{strictly!$\Ts$-positive \protect\nband{\Ts}word} in case
every $\Ts$\nobreakdash-generator appears among the bands $b(s)$.
\end{definition}

For $X=\utupsep{x}1n<\dotsm\in\multsets{\R}{n}$, let
$\Is_X$ \resp{$\Ys_X$} denote any espalier 
$\Ts$ with $\Verts{\Ts}=X$ and 
$\{\Bd\eedge\Suchthat \eedge\in\Edges{\Ts}\}=
\{\{x_p,x_{p+1}\}\Suchthat1\le p<n\}$ 
\resp{$\{\{x_1,x_{p}\}\Suchthat 1< p\le n\}$}.
Among the combinatorial types of trees $\Ts$
with $\Verts{\Ts}=X$, $\Is_X$ and $\Ys_X$ 
represent two extreme types, which may be 
called \moredef{linear}%
   \index{espalier!linear}%
   \index{linear espalier} (minimal number of endpoints) 
and \moredef{star-like}%
   \index{espalier!star-like}%
     \index{star!-like espalier} (maximal number of endpoints), 
respectively; further, among the linear \resp{star-like} espaliers, 
$\Is_X$ \resp{$\Ys_X$} is again extreme, in a sense that 
is obvious and easily formalized.

\begin{example}\label{I-generators}
The $\Is_n$-generators of $\brgp{n}$ are the
standard generators of \dispeqref{standard presentation of Bn}, 
an \nband{\Is_{\intsto{n}}}word 
is just a ``braid word'' in the usual sense 
\citep[][p. 70 ff.]{Birman1974}, and the braid of an
$\Is_n$-positive (or, for some authors, strictly $\Is_n$-positive) 
\nband{\Is_n}word is a ``positive braid'' in the usual sense 
\citep[][etc., etc.]{Birman1974,Rudolph1982b,FranksWilliams1987}.
\end{example}

\begin{question}\label{tree cellulation question}
As noted after \eqref{presentation theorem}, the 
standard generators and standard relators of 
\dispeqref{standard presentation of Bn} 
correspond naturally to the \ncodimension{1} and \ncodimension{2} 
cells of $\config{n}/{\equiv_\R}$.  As noted in 
\eqref{I-generators}, the standard generators also correspond
naturally to the edges of $\Is_n$.  Clearly, given any linear
planar tree $\Ts$ with $\card{\Verts{\Ts}}=n$, from an isotopy
of $\Ts$ to $\Is_n$ may be contrived a cellulation 
$\config{n}/{\equiv_\Ts}$ with a unique \ncodimension{0} cell,
such that 
\begin{inparaenum}
\item\label{q:generators}
the \ncodimension{1} cells of $\config{n}/{\equiv_\Ts}$
correspond to the $\Ts$-generators of $\brgp{\Verts{\Ts}}$, and
\item\label{q:relators}
the \ncodimension{2} cells of $\config{n}/{\equiv_\Ts}$ 
correspond to the relators of \eqref{standard presentation of Bn}.
Now suppose that $\Ts$ is a planar tree for which
there exists a (natural or contrived) cellulation 
$\config{n}/{\equiv_\Ts}$, with a unique \ncodimension{0} cell, 
that has property (\ref{q:generators}) but, 
rather than property (\ref{q:relators}), satisfies both 
\item\label{q:old relators}
some of the \ncodimension{2} cells 
of $\config{n}/{\equiv_\Ts}$ correspond to 
the relators of \dispeqref{T-presentation of non-braid group}, and
\item\label{q:new relators}
the remaining \ncodimension{2} cells
of $\config{n}/{\equiv_\Ts}$ correspond to some family of
extra relators exactly sufficient to convert 
\dispeqref{T-presentation of non-braid group} into a 
presentation of $\brgp{\Verts{\Ts}}$.
\end{inparaenum}
\ask{Does this imply that, in fact, $\Ts$ is linear?}
Specifically, 
\ask{does there exist a cellulation 
$\config{4}/{\equiv_{\Ys_{\intsto{4}}}}$
with properties (\ref{q:generators}),
(\ref{q:old relators}), and (\ref{q:new relators})?}
\end{question}

The methods of \citet*{BirmanKoLee1998} might help answer
\eqref{tree cellulation question} (although the presentation 
they give is not obviously associated to a cellulation).

\subsection{Embedded bandwords and braided Seifert surfaces}
\label{embedded bandwords and braided Seifert surfaces}

Early versions of the construction of braided Seifert surfaces described 
in this section appeared in \cite{Rudolph1983b,Rudolph1983c}.

\begin{construction}\label{braided Seifert surface construction}
Let $X\in \smash[t]{\subsets{n}{\R}}$, $T\in\subsets{k}{\R}$,
say $X\defines\ultuple{x}{n}$, $T\defines\ultuple{t}{k}$. 
Let $0<\e<\min\{|t-t'|/2\Suchthat t, t'\in T,\, t\ne t'\}$;
in case $k>0$, let $I=\clint{\min{T}-\e}{\max{T}+\e}$.
Let $\brep{b}$ be an embedded \nband{X}word of length $k$.
To implement this construction of braided Seifert surfaces,
choose a proper arc $\eedge_s\sub\LHP$ in 
the edge-class $\edgeband{b(s)}$ for each $s\in\intsto{k}$, 
and embeddings $\hparam{x}\from\h0{}\embto\UHP\times\R$ ($x\in X$)
and $\hparam{t}\from\h1{}\embto\LHP\times\opint{t-\e}{t+\e}$ ($t\in T$), 
subject to the following conditions.
For each $x\in X$,
{\parskip6pt
\begin{enumerate*}
    \item\label{brsurf: 0-handles proper} 
    $\hparam{x}$ is proper along a boundary arc of $\h0{}$, 
    \item\label{brsurf: 0-handles stacked}
    $\hparam{x}(\h0{})\sub \{x+\imunit{y}\Suchthat y\in\Rge{0}\}\times \R$,
    \item\label{brsurf: 0-handles tall enough}
    and $I\sub\Bd{\hparam{x}(\h0{})}$.
\suspend{enumerate*}
For each $t\in T$, say $t=t_p$, $p\in\intsto{k}$,
\resume{enumerate*}
    \item\label{brsurf: 1-handles proper}
    $\hparam{t_p}$ is proper along the attaching arcs
    of $\h1{}$, 
    \item\label{brsurf: 1-handles have correct cores}
    $\pr_1\after\hparam{t_p}\restr\core{\h1{}}\from \core{\h1{}}\to\LHP$
    is a diffeomorphism onto $\eedge_p$,
    \item\label{brsurf: 1-handles half-twisted}
    $(\Re\after\pr_1,\pr_2)\after\hparam{t}\from\h1{}\to\R\times\R$
    is a bowtie, and
    \item\label{brsurf: sign condition}
    the sign of the crossing of 
    $(\Re\after\pr_1,\pr_2)\after\hparam{t}\restr\Bd{\h1{}}$
    is equal to $\sign{b(t)}$.
\end{enumerate*}
}
It follows that 
\begin{equation}\label{equation:braided Seifert surface construction}
\smash{%
\Sigma(\brep{b})\isdefinedas
\bigcup_{x\in X}\hparam{x}(\h0{})
\cup 
\bigcup_{t\in T}\hparam{t}(\h1{})
\sub \C\times\R
}
\end{equation}
is a \nhandle{(0,1)} decomposition 
\dispeqref{handle decomposition of a surface} of a surface.  
Any Seifert surface $\bSeif{\brep{b}}\isdefinedas
\d(\Sigma(\brep{b}))\sub \R^3\sub S^3$, where 
$\d\from\C\times\R\to\R^3\suchthat(z,t)\mapsto(\Re{z},\Im{z},t)$
and $\Sigma(\brep{b})$ is constructed as in 
\eqref{braided Seifert surface construction},
is called a \moredef{braided Seifert surface of}%
    \index{braided!Seifert surface}%
    \index{Seifert!surface!braided}%
    \index{surface!Seifert!braided} the 
embedded \nband{X}word $\brep{b}$.  It is easy to see that
$\gbraid{\brep{b}}\isdefinedas
\Sigma(\brep{b})\cap(\LHP\times I)$ is a geometric
braid for $\braidof{\brep{b}}$, and that the link 
$\Bd{\bSeif{\brep{b}}}$ is isotopic in $S^3\supset\R^3$ 
to $\closedbraid{\brep{b}}$.

The \moredef{braid diagram of $\brep{b}$}%
    \index{braid!diagram}%
    \index{diagram!braid ($\brdiagram{\phold}$)}%
    \index{BD(L)@$\brdiagram{\phold}$ (braid diagram)}%
    \index{diagram!braid ($\brdiagram{\phold}$)!information of ($\brinfo{\phold}$)}%
    \index{diagram!braid ($\brdiagram{\phold}$)!picture of ($\brpict{\phold}$)}%
    \index{I@information!about a braid ($\brinfo{\phold}$)}%
    \index{P@picture!of a braid ($\brpict{\phold}$)} is the pair
$\brdiagram{\brep{b}}\defines
(\brpict{\brep{b}},\brinfo{\brep{b}})$, where
$\brpict{\brep{b}}\isdefinedas%
(\Re,\pr_2)(\gbraid{\brep{b}})\sub\R\times I$
and $\brinfo{\brep{b}}$ is the information 
about the signs of the crossings of $\brpict{\brep{b}}$
(indicated graphically in the style of 
\eqref{apparatus of standard diagram}(\ref{crossings of standard diagram});
see \Figref{braided surface construction figure}).
A \moredef{braid diagram for $\beta\in \brgp{X}$} is a braid diagram 
of any bandword $\brep{b}$ with $\braidof{\brep{b}}=\beta$.
A \moredef{standard braid diagram for $\beta\in \brgp{n}$}%
    \index{standard!braid diagram}%
    \index{diagram!braid ($\brdiagram{\phold}$)!standard} is a braid diagram 
of any \nband{\Is_{\intsto{n}}}word $\brep{b}$ with 
$\braidof{\brep{b}}=\beta$.
\begin{figure}
\centering
\includegraphics[width=.7\textwidth]{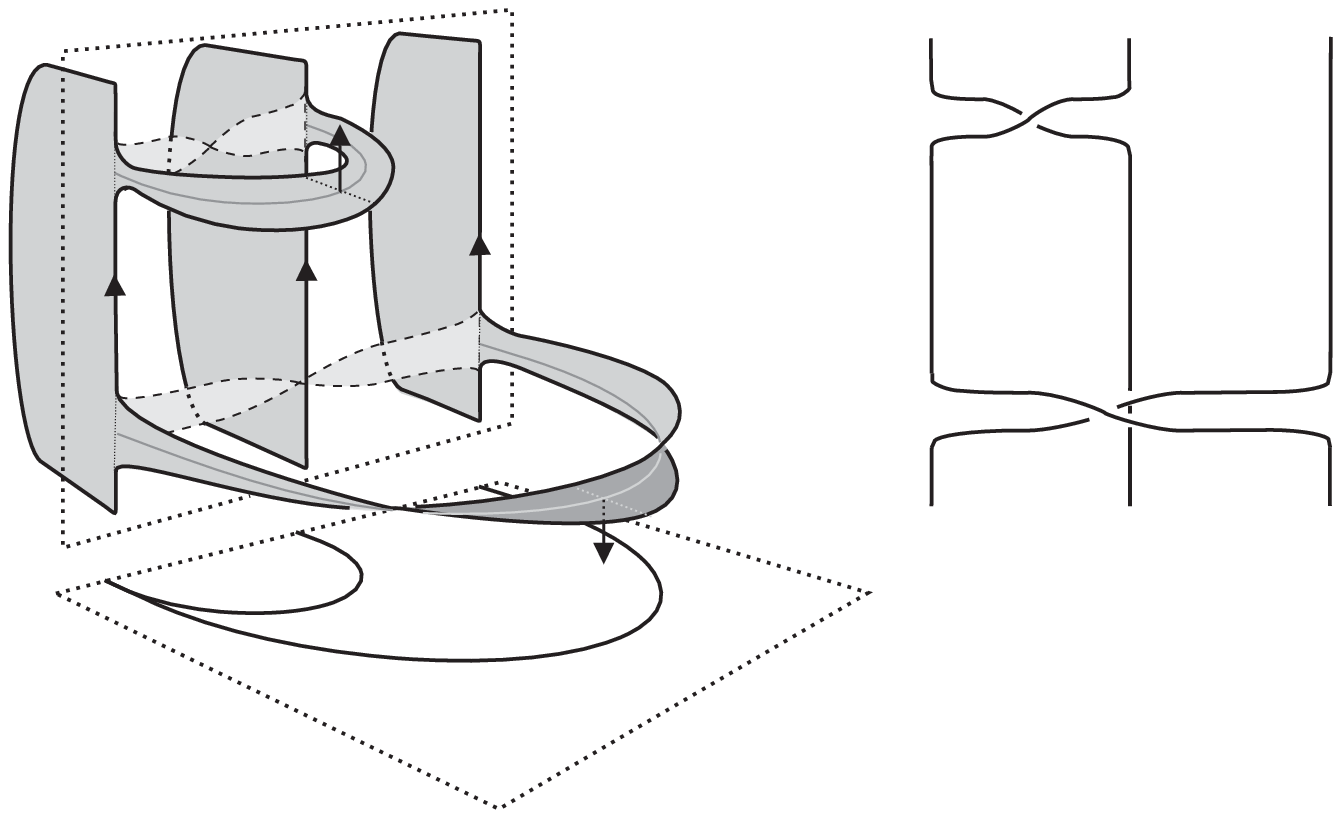}
\caption{Left: a surface $\Sigma(\brep{b})$; 
$\sign{b(1)}={-}$, $\sign{b(1)}={+}$.
Right: a braid diagram for $\braidof{\brep{b}}$.
\label{braided surface construction figure}}
\end{figure}
\end{construction}

\begin{figure}
\centering
\includegraphics[width=.30\textwidth]{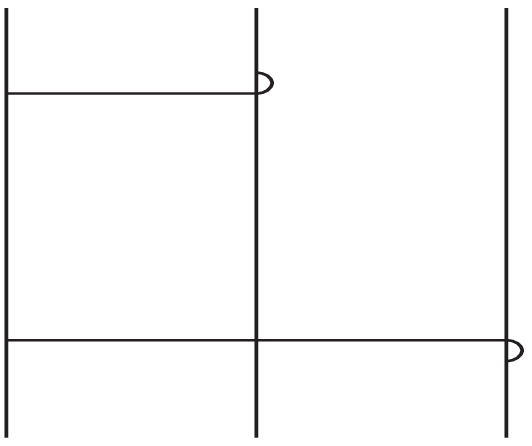}
\caption{A charged fence diagram for the surface
in \protect\Figref{braided surface construction figure}.
\label{fence figure}}
\end{figure}

Of course $\brep{b}$ determines $\bSeif{\brep{b}}$ 
and $\closedbraid{\brep{b}}$ up to isotopy; 
and clearly, up to isotopy (even isotopy through braided 
Seifert surfaces), none of $X$, $T$, or the collection of 
specific edge-class representatives $\eedge_s$ is
necessary \foreign{per se} to the construction of 
$\bSeif{\brep{b}}$---all that is needed is a modicum of 
combinatorial information extracted from $\brep{b}$.  
That information can be encoded as the \ntuple{k} of triples
$((i_\brep{b}(1),j_\brep{b}(1),\varepsilon_\brep{b}(1)),\dots,
(i_\brep{b}(k),j_\brep{b}(k),\varepsilon_\brep{b}(k)))$
with 
$\{x_{i_\brep{b}(s)}<x_{j_\brep{b}(s)}\}\isdefinedas\Bd{\eedge_s}\sub X$
and $\varepsilon_\brep{b}(s)\isdefinedas\sign{b(s)}$.
Another convenient encoding is graphical, using 
\moredef{charged fence diagrams} 
\citep[see][]{Rudolph1992b,Rudolph1998}.  \Figref{fence figure} 
pictures a (very simple) charged fence diagram.\footnote{%
   In \citet{Rudolph1992b}, I described fence diagrams 
   as ``my synthesis of some diagrams that H.~Morton used 
   to describe certain Hopf-plumbed fiber surfaces in 1982~\dots~and 
   `square bridge projections' as described by H.~Lyon'' 
   \citep{Lyon1980}.  By \citeyear{Rudolph2001a}, I had recalled
   that I first saw fences (and braid groups!) ``c.\ 1959, in 
   one of Martin Gardner's `Mathematical Games' columns'' 
   \citep[reprinted in][as Chapter 2, ``Group Theory and Braids'']%
   {Gardner1966}.  Since 2001 I have become aware of 
   ``amida-diagrams'' 
   \citep[introduced by][precisely to construct certain special
   Seifert surfaces]{Yamamoto1978} and ``wiring diagrams'' 
   \citep[cf.][and other literature from the theory of line
   arrangements]{CordovilFachada1995}, both closely related to
   fences, as are (via their essential identity with square
   bridge projections) the ``barber-pole projections'' attributed
   in \citet{Rudolph1992b} to Thurston (unpublished), 
   \citet{Erlandsson1981}, and \citet{Kuhn1984}.  
   Is this another instance of \citeauthor{Epple1995}'s
   ``elimination of contexts'', or a mere multiplication of contexts?}%

There is a more or less canonical way to turn an 
\nstring{n} braid diagram $\brdiagram{\brep{b}}$ 
into a standard link diagram $\diagram{\closedbraid{\brep{b}}}$:
identify $\C$ with $\R^2\supset\brpict{\brep{b}}$ by
$(\Re,\Im)$; in $\C$, attach $n$ arcs to 
$\brpict{\brep{b}}\sub\C$, preserving orientation,
so as to create a normal collection $\pict{\closedbraid{\brep{b}}}$ 
of closed curves with no new crossings;
and let $\info{\closedbraid{\brep{b}}}$ 
be $\brinfo{\brep{b}}$.

\begin{unproved}{proposition}
\label{braid diagrams and Seifert surfaces}
\begin{inparaenum}
\item\label{In-word braided surfaces are nested diagrammatic}
Let $\brep{b}$ be an \nband{\Is_{\intsto{n}}}word.  
Let $L=\closedbraid{\brep{b}}$.  The arcs
$\diff{\pict{L}}{\Int{\brpict{\brep{b}}}}$
can be chosen so that 
$\SeifO{\diagram{L}}=\SeifO{\diagram{L}}^{+}$;
if they are, then $\diagram{L}$ is nested
and the diagrammatic Seifert surface 
$\Seifsurf{L}$ is isotopic to $\bSeif{\brep{b}}$.
\item\label{nested diagrammatic surfaces are In-word braided}
Conversely, if $\diagram{L}$ is nested and 
$\SeifO{\diagram{L}}=\SeifO{\diagram{L}}^{+}$, then
\textup(up to isotopy of $\pict{L}$\textup)
there exists an \nband{\Is_{\intsto{n}}}word $\brep{b}$,
$n\isdefinedas\card{\SeifO{\diagram{L}}}$, 
such that $L=\closedbraid{\brep{b}}$ and $\Seifsurf{L}$ 
is isotopic to $\bSeif{\brep{b}}$.
\item\label{non-In-word braided surfaces are not diagrammatic}
If $\brep{b}$ is an embedded \nband{\intsto{n}}word but not 
an \nband{\Is_{\intsto{n}}}word, and 
$L=\closedbraid{\brep{b}}$, 
then $\diagram{L}$ is nested
but $\euler{\Seifsurf{L}}<\euler{\bSeif{\brep{b}}}$,
so $\Seifsurf{L}$ and $\bSeif{\brep{b}}$ are not 
diffeomorphic, let alone isotopic.
\end{inparaenum}
\end{unproved}

\eqref{braid diagrams and Seifert surfaces} shows that 
\eqref{MFW--no braids} is equivalent to the following.

\begin{unproved}{theorem}
\label{MFW--braids}
For all $n$, for all $\beta\in\brgp{n}$, 
\begin{equation}
\ord_v P_{\widehat\beta}\geq e(\beta)-n+1,
\end{equation}
where $e\from\brgp{n}\to\Z$ is abelianization \textup(exponent
sum\textup).
\end{unproved}

Proofs of the following theorem apppear in 
\citet{Rudolph1983b,Rudolph1983c}; it also 
follows from results of \citet{Bennequin1983} on ``Markov surfaces''
\citep[see][]{Rudolph1985}. 

\begin{unproved}{theorem}
\label{every Seifert surface can be braided}
If $S\sub \R^3\sub S^3$ is a Seifert surface,
then there exist $X$ and an embedded \nband{X}word $\brep{b}$
such that $S$ is isotopic to $\bSeif{\brep{b}}$.
\end{unproved}

\begin{question}\label{estimates for braided surfaces}
The \moredef{braid index}%
    \index{braid!index ($\bi{\phold}$)} $\bi{L}$ of a link $L$ 
is the minimum $n\in\Nstrict$ such that $L$ is isotopic to 
$\closedbraid{\brep{b}}$ for some \nband{n}word $\brep{b}$
(which can, since its length is not an issue, be taken to
be embedded); see \citet{BirmanBrendle2004}.
Let the \moredef{braided Seifert surface index}%
    \index{braided!Seifert surface!index} of a Seifert surface
$S$ be the minimum $n\in\Nstrict$ such that $S$ is 
isotopic to $\bSeif{\brep{b}}$ for some
embedded \nband{\intsto{n}}word $\brep{b}$;
by \eqref{every Seifert surface can be braided}, 
this invariant of $S$ is an integer.
\ask{How can it be calculated or estimated?}
This question is the analogue for Seifert surfaces of 
``Open Problem~1'' of \citet{BirmanBrendle2004},
which asks how to calculate $\bi{L}$.
Clearly the braided Seifert surface index of $S$ is at least as great
as $\bi{\Bd{S}}$.  Examples of \cite{HirasawaStoimenow2003} show that 
this inequality can be strict.
\end{question}

\begin{definition}\label{def:strongly quasipositive}
A Seifert surface $S$ is \moredefdanger{quasipositive}%
    \index{quasipositive!Seifert surface}%
    \index{Seifert!surface!quasipositive}%
    \index{surface!Seifert!quasipositive} in case $S$ is 
isotopic to a braided Seifert surface $\bSeif{\brep{b}}$ 
for some positive \nband{X}word $\brep{b}$.  A 
link is \bydef{strongly quasipositive}%
    \index{quasipositive!strongly} in case it has a 
quasipositive Seifert surface.
\end{definition}

\subsection{Plumbing and braided Seifert surfaces}
\label{plumbing and braided Seifert surfaces}

Let $X\in\subsets{n}{\R}$.  For $x\in X$, write
$\Xle{X}x\isdefinedas X\cap\Rle{x}$, 
$\Xge{X}x\isdefinedas X\cap\Rge{x}$.

\begin{definitions}
Let $x\in X$.
Let $\brep{b}$ be an embedded \nband{X}word of length $k$.  
Let $T_{\brep{b};\le x}\isdefinedas%
\{s\in\intsto{k}\Suchthat b(s)\in\brinj{\Xle{X}x,X}\}$ and
$T_{\brep{b};\ge x}\isdefinedas%
\{s\in\intsto{k}\Suchthat b(s)\in\brinj{\Xge{X}x,X}\}$; 
let $k_{\ssle}\isdefinedas\card{T_{\brep{b};\le x}}$,
$k_{\ssge}\isdefinedas\card{T_{\brep{b};\ge x}}$, so
in every case $k_{\ssle}+k_{\ssge}\le k$.  In case 
$k_{\ssle}+k_{\ssge}=k$, say that $\brep{b}$ is 
\moredef{deplumbed by $x\in X$}%
    \index{deplumbed}.  Writing 
$T_{\brep{b};\le x}\defines\ultuple{u}{k_{\ssle}}$,
$T_{\brep{b};\ge x}\defines\ultuple{v}{k_{\ssge}}$,
let $\brep{b}_{\le x}\isdefinedas\brepvec{b_{\le x}}{k_{\le}}$ 
\resp{$\brep{b}_{\ge x}\isdefinedas\brepvec{b_{\ge x}}{k_{\ge}}$}
be the embedded \nband{\Xle{X}x}word 
\resp{\nband{\Xge{X}x}word} 
for which $\brinj{\Xle{X}x,X})b_{\le x}(s)=b(u_s)$
\resp{$\brinj{\Xge{X}x,X})b_{\ge x}(t)=b(v_t)$}.
By construction \eqref{braided Seifert surface construction}, 
if $\brep{b}$ is deplumbed by $x\in X$, then 
$\Sigma(\brep{b})=\Sigma(\brep{b}_{\le x})\cup\Sigma(\brep{b}_{\ge x})$
and $\Sigma(\brep{b}_{\le x})\cap\Sigma(\brep{b}_{\ge x})=%
\Sigma(\brep{b}\cap\C\times\{x\}$ is a disk (in fact, a
\nhandle{0} of the constructed \nhandle{(0,1)} decompositions
of $\Sigma(\brep{b}_{\le x})\cup\Sigma(\brep{b}_{\ge x})$);
the situation is illustrated using fence diagrams in 
\Figref{deplumbing figure}.
\begin{figure}
\centering
\includegraphics[width=.7\textwidth]{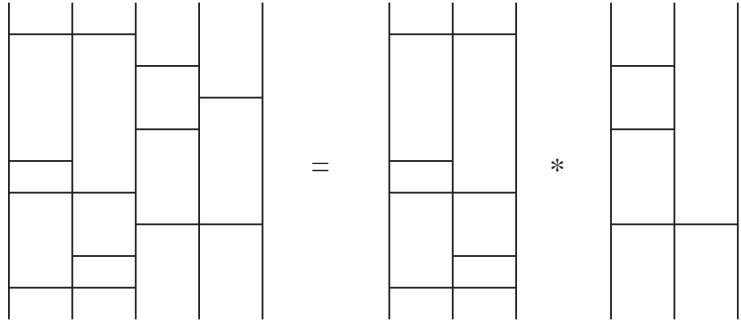}
\caption{An example of deplumbing.
\label{deplumbing figure}}
\end{figure}
The inverse of this operation,
which combines two braided Seifert surfaces into one,
\moredef{braided Stallings plumbing}%
    \index{braided!Stallings plumbing}%
    \index{Stallings plumbing!braided}%
    \index{plumbing!braided}. 
\end{definitions}

Although it appears to be a rather special case 
of \bydef{Stallings plumbing}%
    \index{plumbing!Stallings (Murasugi sum)}, or 
    \bydef{Murasugi sum}%
    \index{sum!Murasugi}, a geometric operation on Seifert 
surfaces that has been extensively studied 
\citep[by][and many others; notations, definitions, and 
some history are given by \citealp{Rudolph1998}]%
{Murasugi1963,
Conway1970,
Siebenmann1975,
Stallings1978,
Gabai1983a,
Gabai1983b,
Gabai1986},
braided Stallings plumbing is in fact perfectly general,
according to \citet{Rudolph1998}.

\begin{unproved}{theorem}
Up to isotopy, every Murasugi sum $S_0\plumb{}S_1$
of Seifert surfaces $S_0$, $S_1$ is a braided Stallings 
plumbing of braided Seifert surfaces.  
\end{unproved}

\subsection{Labyrinths, braided surfaces in bidisks, and braided ribbons}
\label{labyrinths, braided surfaces in bidisks, and braided ribbons}

Let $X$ be a topological space. 
Given $f\from X\to\monicpoly{n}$, let 
$\Weier{f}\from X\times\C\to\C\suchthat(x,w)\mapsto f(x)(w)$,
so that $\gr{\resolve\after f}=\Weier{f}^{-1}(0)$%
    \index{$weier$@$\Weier{\phold}$ (Weierstrassification)}
    \index{$unWeier$@$\unWeier{\phold}$ (unWeierstrassification)}. Generally 
it can be a subtle matter to determine whether a map
$F\from X\times\C\to\C$ (with fibers of cardinality bounded
by $n\in\N$) is $\Weier{f}$ for some (continuous)
$f\from X\to\monicpoly{n}$, but in case $F$ is known to be such, 
write $f=\unWeier{F}$.  

\begin{definition}
A map $f\from X\to\monicpoly{n}$ is \moredef{amazing}%
    \index{amazing map} provided that the partition of $M$ 
into connected components of inverse images $f^{-1}(C)$ of 
the various cells $C$ of $\monicpoly{n}/{\equiv_\R}$ 
has some finite refinement that is 
a stratification of $X$.  A \bydef{labyrinth} for $f$ 
is any such stratification.
\end{definition}

\begin{examples}
Amazing maps have proved useful in 
several applications.
\begin{asparaenum}[\bf \theexamples.1.]
\item\label{LR:algebraic functions}
If $F(z,w)\isdefinedas 
f_0(z)w^n + f_1(z)w^{n-1} + \dotsm + f_n(z)\in \C[z,w]$
is a polynomial in Weierstrass form, then 
$\unWeier{F}\from \C\to\monicpoly\C$ is (in effect)
$F$ considered as an element of $\C[z][w]\iso\C[z,w]$,
and is amazing.  Classically, in this situation 
$w$ is said to be the ``$n$-valued algebraic function 
of $z$ without poles'' determined by the algebraic curve 
$\variety{F}$ \citep{Bliss1933,Hansen1988}, and 
a more or less canonical set of ``branch cuts'' for
$w$ can be extracted systematically from the labyrinth
of $\unWeier{F}$.
\citet{Rudolph1983a} determined the topology 
of the labyrinth of $\unWeier{F}$ in the case 
$F_\e(z,w)=(w-1)(w-2)\dotsm(w-n+1)(w-z)+\e$,
first for $\e=0$, then for small $\e\ne 0$, 
and used it to give the first proof of 
\eqref{qp implies transverse-C}.
\item\label{Dung:labyrinth}
The term labyrinth was introduced by \citet{DungHa1995},
in their study of line arrangements in $\C^2$.  In a 
linear coordinate system chosen so that no line in the 
arrangement is vertical, such an arrangement is $\variety{F}$, 
where $F(z,w)\in\C[z,w]$ is a product of (unrepeated)
factors $w+az+b$ (as in (\ref{LR:algebraic functions}), where
$F_0$ corresponds to an arrangement of $n-1$ horizontal lines
crossed by a single diagonal line).  Again, $\unWeier{F}$ exists
and is amazing.  
(See \Secref{Alexander invariants} for further references on arrangements.)%
\item\label{Orevkov:Zariski labyrinth}
$\mathstrut$
\citet{Orevkov1988} made excellent use of 
the observation that, if $V\sub \C^2$ is both 
a $\C$-algebraic curve and a nodal surface, then,
after an arbitrarily small linear change of 
coordinates in $\C^2$, $V=\variety{F}$ where 
$F\in\C[z,w]$ is in Weierstrass form and 
the labyrinth of the amazing map $\unWeier{F}$ 
has a very special form.
(See \Secref{Zariski conjecture} for further discussion of this 
application.)
\item\label{Orevkov:Rudolph diagrams}
$\mathstrut$
\citet{Orevkov1996} developed the theory of labyrinths
in careful detail, and applied it to the Jacobian Conjecture. 
(See \Secref{Keller's Jacobian Problem}).
\end{asparaenum}
\end{examples}

Amazing maps, and the facts about the 
low-codimension cells of $\monicpoly{n}/{\equiv_\R}$
laid out in \eqref{low-codimension cells of MPn}, 
together allow two closely related constructions---of 
braided surfaces in bidisks, and braided ribbons in $D^4$---to 
be described and carried out more precisely than in the 
original sources \citep{Rudolph1983b,Rudolph1985}.

\begin{proposition}
\label{transversality results}
Let $M$ be a surface.  If $f\from M\to\monicpoly{n}$ is
transverse to all the cells of $\monicpoly{n}/{\equiv_\R}$,
then $f$ is amazing.  More specifically, $f$ has a 
labyrinth $M/{\equiv}$ with the following properties.
{\parskip6pt
\begin{enumerate*}
\item
The union of the vertices and edges of $M/{\equiv}$ is 
a graph $\Lambda(f)$ with no vertex of degree $d\notin\{1,4,6\}$.
\item
The association by $f$ to each edge $\eedge$ of $\Lambda(f)$
of a \ncodimension{1} cell of $\monicpoly{n}/{\equiv_\R}$ 
endows $\eedge$ with 
\begin{inparaenum}
\item
a natural transverse orientation
\textup(and therefore a natural orientation, so that
$\eedge$ is naturally either a simple closed curve or an arc\textup)
and 
\item\label{labyrinth edges}
a \bydefdanger{clew} 
$\sigma(\eedge)$ in the set $\{\brgen{1},\dots,\brgen{n-1}\}$
of standard generators of the presentation
\dispeqref{standard presentation of Bn} of $\brgp{n}$.
\end{inparaenum}
\item\label{labyrinth vertices}
The association by $f$ to each vertex $x$ of $\Lambda(f)$
of a \ncodimension{2} cell of $\monicpoly{n}/{\equiv_\R}$ 
endows $x$ with either
\begin{inparaenum}
\item\label{valence 4 or 6 vertices}
a standard relator of \dispeqref{standard
presentation of Bn},
in case $\valfun{\Lambda(f)}(x)\in\{4,6\}$, or
\item\label{endpoints}
a trivial relator $\brgen{s}^{\pm1}\brgen{s}^{\mp1}$ 
in case $\valfun{\Lambda(f)}(x)=1$.
\end{inparaenum}
\item
The clews \textup{(2\ref{labyrinth edges})} are consistent with the
relators \textup{(\ref{labyrinth vertices})} as pictured in 
\Figref{labyrinth figure}.  \textup(In each local picture,
the mirror image of the illustrated situation is also allowed, 
as is simultaneous reversal of all edge orientations.\textup)
\qquad\qquad\qedsymbol
\end{enumerate*}}
\begin{figure}
\centering
\includegraphics[width=.9\textwidth]{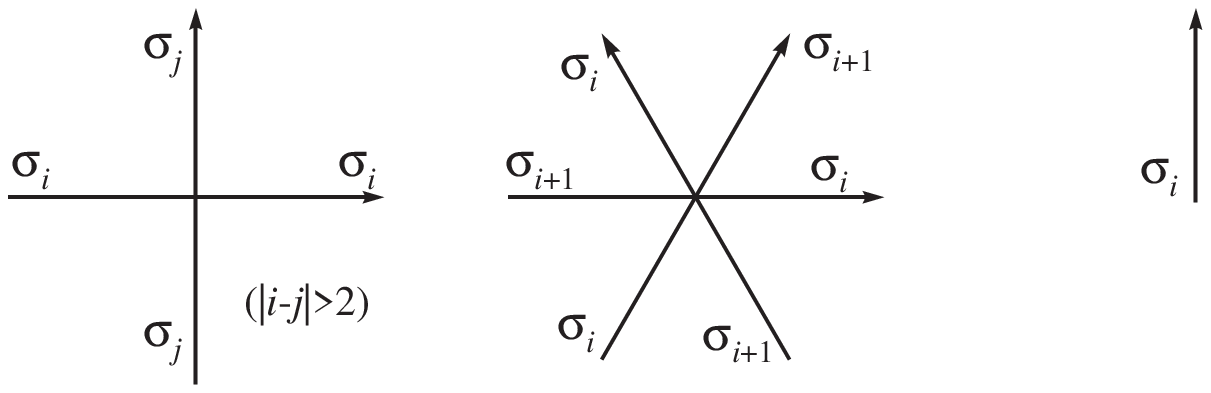}
\caption{Clews in the labyrinth of a 
$\monicpoly{n}/{\equiv_\R}$-transverse map on a
surface.\label{labyrinth figure}}
\end{figure}
\end{proposition}


\begin{construction}\label{braided surfaces in DxD and D4}
A \moredef{braided surface}%
    \index{braided!surface} of \moredef{degree $n$}%
    \index{degree (of a braided surface)} in $D^2\times\C$ is 
$S_f\isdefinedas\Weier{f}^{-1}(0)$, where 
$f\from D^2\to\monicpoly{n}$ is transverse to all the cells 
of $\monicpoly{n}/{\equiv_\R}$ 
as in \eqref{transversality results}. 
\end{construction}

\begin{unproved}{proposition}
A braided surface $S_f$ in $D^2\times\C$ is a surface, and 
$\pr_1\restr S\from S\to D^2$ is a simple
branched covering map branched over $\Verts{\Lambda(f)}$.
\end{unproved}

In these terms, various facts about
the constructions in \citep{Rudolph1983b,Rudolph1985}
can be phrased as follows.
Since any braided surface $S_f$ is compact, it is actually
contained in $D^2\times rD^2$ for some $r$; clearly only braided
surfaces in $D^2\times D^2$ need be considered.  The obvious
identification of $S^1\times D^2\sub\Bd(D^2\times D^2)$ with
$\Nb{\hunknot}{S^3}$ identifies the boundary of a braided
surface $S_f\in D^2\times D^2$ with a \nsubmanifold{1} of 
$\Nb{\hunknot}{S^3}$, which is a closed 
\nstring{n} $\ob{o}$-braid isotopic to some $\closedbraid{\brep{b}}$
of length $\card{\Verts{\Lambda(f)}}$.  The bandword $\brep{b}$ can 
be read off from the clews of $\Lambda(f)$.  In particular, 
$\braidof{\brep{b}}$ is the ordered product of the clews on
the (oriented) edges of $\Lambda(f)$ that intersect $S^1$.
Every bandword determines a braided surface which is unique
up to isotopy (through braided surfaces). 
\citet{Rudolph1983a} characterizes quasipositive braided 
surfaces as follows.

\begin{unproved}{theorem}\label{quasipositive ribbons are holomorphic}
A braided surface $S_f$, where $f$ is holomorphic, determines
a quasipositive bandword.  Any braided surface that determines
a quasipositive bandword is isotopic 
\textup(through braided surfaces\textup)
to $S_f$ with $f$ a polynomial.
\end{unproved}

There is a straightforward way to ``round the corners''
of $D^2\times D^2$ 
by a \bydef{smoothing} $D^2\times D^2\to D^4$---that is, a
homeomorphism which is a diffeomorphism off $S^1\times S^1$---so 
that a braided surface
in $D^2\times D^2$ is carried to a ribbon surface in $D^4$.  
The converse is also true \citep{Rudolph1983b,Rudolph1985}.

\begin{unproved}{theorem}\label{ribbons are braidable}
Up to isotopy, every ribbon surface in $D^4$ is the image 
of a braided surface in $D^2\times D^2$ by a smoothing 
$D^2\times D^2\to D^4$.
\end{unproved}

\begin{unproved}{theorem}\label{ribbons are ribbons}
Let $\brep{b}$ be an embedded bandword.  The ribbon
surface $S$ in $D^4$ with $\Bd S=\Bd\bSeif{\brep{b}}\sub S^3$,
produced by pushing $\Int{\bSeif{\brep{b}}}$ into $\Int{D^4}$,
is isotopic to the ribbon surface which is the image of the
braided surface in $D^2\times D^2$ determined by the bandword
$\brep{b}$.
\end{unproved}

\begin{definitions}
\label{nodes and cusps in braid groups}
A \moredef{node}%
    \index{node!in a braid group} \resp{\moredef{cusp}}%
    \index{cusp!in a braid group} in an \nstring{n} braid group
$\brgp{X}$ is the square \resp{cube} of a \nband{X}.  A \bydef{nodeword} 
\resp{\bydef{cuspword}} in $\brgp{X}$ is a \ntuple{k} 
$\brep{b}\defines\brepvec{b}{k}$ such that each $b(i)$
is a band or a node \resp{a band, a node, or a cusp} in $\brgp{X}$.
\end{definitions}

The proof of \eqref{ribbons are braidable} can be generalized 
\citep[cf.][]{Rudolph1983b} to show that every ``nodal ribbon surface''
in $D^4$ can be realized, up to smoothing and isotopy, by 
a ``nodal braided surface'' in $D^2\times D^2$.
\citet{Orevkov1998} proved an immense generalization of 
\eqref{quasipositive ribbons are holomorphic}; the following
special case of \citeauthor{Orevkov1998}'s result is easy
to state and more than adequate for this review.

\begin{unproved}{theorem}\label{Orevkov's monodromy theorem}
If $\brep{b}$ is a quasipositive cuspword in $\brgp{n}$, 
then there is a holomorphic cusp curve in $D^2\times D^2$
that corresponds to $\brep{b}$ via an appropriate holomorphic map
$f\from D^2\to\monicpoly{n}$.
\end{unproved}

The map $f$ in \eqref{Orevkov's monodromy theorem} is amazing;
it is transverse to all the cells of $\monicpoly{n}/{\equiv_\R}$ 
if and only if the cuspword $\brep{b}$ is a bandword.


%
%
\section{Transverse \texorpdfstring{$\C$}{C}-links}\label{transver}

\begin{definition}\label{def:transverse C-link} 
A link $L\sub S^3$ is a \moredef{transverse {\Clink}link}%
    \index{transverse!C-link@{\Clink}link}%
    \index{C-link@{\Clink}link!transverse} provided that 
$(S^3,L)$ is diffeomorphic to $(\sph,\variety{f}\cap\sph)$,
where:
\begin{inparaenum}%
\item\label{tcl:strictly pseudoconvex}
$D\sub\C^2$ is a Stein disk bounded by
a strictly pseudoconvex \nsphere{3} $\sph$, and 
\item\label{tcl:singularities allowed off sphere}
$f\in\holo{D}$ is a holomorphic function with
$\sing{\variety{f}}\cap\sph=\emptyset$, such that 
\item
the complex manifold $\reg{\variety{f}}$ intersects $\sph$ 
transversely and
\item
$\variety{f}\cap\sph$ is non-empty (and therefore a smooth,
naturally oriented closed \nsubmanifold{1} of $\sph$).
\end{inparaenum}
\end{definition}

Various modifications may be made to the specifications in
\eqref{def:transverse C-link} without changing the class of
links so specified.

\begin{proposition}\label{tcl:modifications}
\begin{inparaenum}
\item\label{tcl:relax strictly pseudoconvex}
The requirement that $\sph$ be strictly pseudoconvex 
can be somewhat relaxed. 
\item\label{tcl:strengthen strictly pseudoconvex}
The requirement that $\sph$ be strictly pseudoconvex 
can be considerably strengthened: $\sph$ can be required
to be convex, or even to be a round sphere.
\item\label{tcl:restrict variety}
$\variety{f}$ can be required to be non-singular,
to be algebraic, or to be both at once.
\end{inparaenum}
\end{proposition}

\begin{proof}[Proof \textup(sketches\textup)]
\begin{inparaenum}
\item\label{tcl:relax strictly pseudoconvex proof}
Here, no attempt will be made to state the most general theorem 
(see \eqref{transverse C-links questions}(\ref{fancy balls question})).  
A sufficient
example for present purposes is the following.  Let $\g\sub\C$
be an arbitrary simple closed curve and $G\sub\C$ the 
\ndisk{2} that it bounds, so that for any $r>0$ the 
\bydef{bicylinder} $D\isdefinedas G\times rD^2$ is 
diffeomorphic \textup(indeed, biholomorphic\textup) 
to the bidisk $D_{1,1}$.  Although
$D$ is not a Stein domain, $\Int{D}$
is an open Stein manifold, 
and $(D,\Bd{D})$ can be arbitrarily closely
approximated by Stein domains.  The result follows from
\eqref{polynomial convexity} and transversality.
\item
\label{strictly pseudoconvex proof}
As pointed out by \citet{BoileauOrevkov2001}, theorems of
\citet{Eliashberg1990b} show that for this purpose all Stein 
disks are equally good.  
\item\label{tcl:restrict variety, proof}
A nearby level set of $f$ will have no singularities and give 
the same link type; by \eqref{polynomial convexity}, $f$ is 
arbitrarily closely approximated by its sufficiently high-degree 
Taylor polynomials, and again we get the same link type.
\end{inparaenum}
\end{proof}

\begin{questions}\label{transverse C-links questions}
\eqref{def:transverse C-link} and \eqref{tcl:modifications}
immediately suggest a number of questions.
\begin{asparaenum}[\bf \thequestions.1.]
\item
\ask{Is there an algorithm for determining whether or not
a given link is a transverse {\Clink}link?}  (It follows
from \eqref{qp implies transverse-C} that the set of
isotopy classes of transverse {\Clink}links is recursive:
there exists an algorithm that produces a link in every
such isotopy class.  Thus the question becomes: Is
the set of isotopy classes of transverse {\Clink}links 
recursively enumerable?)
\item\label{fancy balls question}
\ask{How much can the requirement that $\sph$ be strictly pseudoconvex be 
relaxed (as in \eqref{tcl:modifications}%
(\ref{tcl:relax strictly pseudoconvex}))?}
See \citet{Rudolph1985a} for a cautionary result.
\item\label{tcl degree question}
Let $L$ be a transverse {\Clink}link.
By \eqref{tcl:modifications}(\ref{tcl:restrict variety}), 
up to isotopy $L=S^3\cap\variety{f}$ for some $f\in\C[z,w]$.
\ask{What is the minimum degree of such a polynomial~$f$?}
In light of 
\eqref{tcl:modifications}(\ref{tcl:strengthen strictly pseudoconvex}),
the same question can be asked with the roundness of $S^3$ 
weakened to convexity or strict pseudoconvexity of 
a \nsphere{3} $\sph$.  \ask{Can either of these weakenings 
strictly decrease the minimum?}
Calculations of, or even good upper bounds for, this (or
these) invariant(s) of $L$ would have applications 
to finding embeddings of certain Stein
domains (namely, cyclic branched covers of $D^4$ branched
over holomorphic curves) into algebraic surfaces of
unexpectedly low degree in $\C^3$ \citep{BoileauRudolph1995x}.
\end{asparaenum}
\end{questions}

\subsection{Transverse \texorpdfstring{$\C$}{C}-links 
are the same as quasipositive links}

\begin{unproved}{theorem}
\label{qp implies transverse-C}
Every quasipositive link is a transverse {\Clink}link.
\end{unproved}

Proofs of \eqref{qp implies transverse-C} (with some variation 
in the details) are presented in 
\citep{Rudolph1983a,Rudolph1984b,Rudolph1985}.
A remarkable theorem of \citet{BoileauOrevkov2001} asserts
the converse.

\begin{unproved}{theorem}\label{Boileau-Orevkov Theorem}
Every transverse {\Clink}link is a quasipositive link.
\end{unproved}

\begin{unproved}{corollary}
$L$ is a transverse {\Clink}link if and only if $L$ is quasipositive.
\end{unproved}

\begin{unproved}{corollary}\label{polynomials are enough}
Every isotopy class of transverse {\Clink}links is represented by a 
transverse intersection 
$\variety{f}\cap \{(z,w)\in\C^2\Suchthat \norm{(z,w)}=1\}$, where
$f\in\C[z,w]$ and $\sing{\variety{f}}=\emptyset$.
\end{unproved}

As \citeauthor{BoileauOrevkov2001} note, their proof (the only
one known to date) is completely non-constructive, relying 
strongly as it does on the theory of pseudoholomorphic curves
\citep{Gromov1985}.  

\begin{question}\label{constructive B-O question}
\ask{Is there a constructive proof of \eqref{Boileau-Orevkov Theorem}?}
Biding such a proof, \ask{can an upper bound on 
$\{ n \Suchthat \text{$L$ is isotopic to a quasipositive \nstring{n} 
braid}\}$ be deduced from other invariants of a 
transverse {\Clink}link $L$?}  
\end{question}

\subsection{Slice genus and unknotting number of
transverse \texorpdfstring{$\C$}{C}-links}

\begin{unproved}{theorem}\label{K-M theorem}
If $f\in\holo{D^4}$ is such that $\sing{\variety{f}}=\emptyset$,
and $L$ is the transverse {\Clink}link $\variety{f}\cap S^3$,
then $\eulermax{s}{L}=\euler{\variety{f}\cap D^4}$.
\end{unproved}

\begin{unproved}{corollary}\label{slice-B equality for QP knots}
If $\brep{b}$ is a quasipositive \nstring{n} bandword of length $k$,
then $\eulermax{s}{\closedbraid{\brep{b}}}=n-k$.
\end{unproved}

It had been known at least as early as 1982 
\citep{BoileauWeber1982,Rudolph1983d} that these results 
would follow from a local version of the Thom Conjecture
(see \eqref{Thom conjecture}): either 
(in the case of \eqref{K-M theorem}) that $L$ 
be a torus link $O\{n,n\}$ for some $n\in\Nstrict$,
or (in the case of \eqref{slice-B equality for QP knots})
that $\brep{b}$ be the important \nband{\Is_{\intsto{n}}}word\footnote{%
   There is a long and worthy tradition of designating 
   (the braid of) this \nband{\Is_{\intsto{n}}}word 
   as $\D^2$ \citep[][etc., etc.]{Birman1974,BirmanBrendle2004}, 
   a usage apparently introduced by \citet{Fadell1962} in homage 
   to P.~A.~M.~Dirac \citep[who wrote that he ``first thought of the string 
   problem about 1929'' in a letter quoted by \citeauthor{Gardner1966}, 
   \citeyear{Gardner1966}, addendum to Chapter 2; see also][]{Newman1942}.  
   Nonetheless, in the context of the knot theory of complex plane
   curves, where discriminants (which have their own, algebro-geometric,
   traditional claim on the notation $\Delta$) abound, I believe
   that the makeshift of using $\nabla$ to stand in for $\D^2$ is 
   preferable to the further overloading of the symbol $\D$.}
\begin{equation}
(\overbrace{%
\brgen{1},\brgen{2},\dots,\brgen{n-1},\dots,
\brgen{1},\brgen{2},\dots,\brgen{n-1}}
^{\text{$n$ repetitions}})\defines\nabla_n
\end{equation}
with closed braid $\closedbraid{\nabla_n}=O\{n,n\}$.
The first proof of the local Thom Conjecture was given by
\citet{KronheimerMrowka1993}, using gauge theory for embedded 
surfaces, a \ndimensional{4} technique.  The full Thom 
Conjecture (and more) was proved using Seiberg--Witten 
theory---also a \ndimensional{4} technique---by 
\citet*{KronheimerMrowka1994,MorganSzaboTaubes1996}.
Recently \citet{Rasmussen2004} announced a purely \ndimensional{3} 
combinatorial proof based on Khovanov homology.

\begin{unproved}{corollary}
If $L$ is a transverse \Clink{link} then 
$\eulermax{r}{L}=\eulermax{s}{L}$.
\end{unproved}

It was also known \citep{BoileauWeber1982,Rudolph1983d} that
the truth of the local Thom Conjecture implies
an affirmative answer to ``Milnor's Question'' 
(see \eqref{Milnor's question}) on unknotting numbers
of links of singularities (\Secref{links of singularities})
or---more generally \citep[see][]{Rudolph1983b}---closed positive braids 
in the sense of \eqref{I-generators}.

\begin{unproved}{corollary}
\label{unknotting number of link of singularity}
If $L$ is the link of a singularity \textup(e.g.,
a torus link $O\{m,n\}$ with $0\le m<n$\textup), then 
$\usz(L)=\node{L}=(\card{\compsof{L}}-\eulermax{}{L})/2$.
\end{unproved}

\subsection{Strongly quasipositive links}
\label{strongly quasipositive links}

\eqref{K-M theorem}\thru\eqref{unknotting number of link of singularity}
immediately imply the following.

\begin{unproved}{theorem}\label{various invariants of sqp links}
If $L$ is a strongly quasipositive link with quasipositive
Seifert surface $S$, then 
$\euler{S}=\eulermax{}{L}=\eulermax{r}{L}=\eulermax{s}{L}$.
If $K$ is a strongly quasipositive knot with quasipositive
Seifert surface $S$, then 
$\genus{}{S}=\genus{}K=\genus{r}{K}=\genus{s}{K}\le\node{K}\le\usz(K)$.
\end{unproved}

To fully exploit \eqref{various invariants of sqp links} it is 
necessary to have a good supply of transverse {\Clink}links
or, equivalently, quasipositive (closed) braids.  
Clearly every strongly quasipositive link 
(\eqref{def:strongly quasipositive}) is quasipositive.  
The converse is false; for instance, the closed braid
of $\conjugate{\brgen{2}^3}{\brgen{1}}\brgen{1}\in\brgp{3}$
is a quasipositive knot, but not strongly quasipositive
(by \eqref{slice-B equality for QP knots}).  
Nonetheless, the class of strongly quasipositive 
links is very varied.  

\subsubsection{\texorpdfstring{$S$}{S}-equivalence and strong quasipositivity}
\label{S-equivalence and strong quasipositivity}

\begin{proposition}\label{qp is Seifert-undistinguished}
No link invariant calculable from a Seifert matrix 
\textup(for instance, the Alexander polynomial of a knot,
the signatures of a knot or link, etc., etc.\textup)
can tell whether or not a link is quasipositive. 
\end{proposition}

\begin{proof}
Let $\Seifmatrix{\Seifpair{L_i}{L_j}{S}}$ be a Seifert
matrix for a link $L$ with Seifert surface $S$.
By \eqref{every Seifert surface can be braided},
we may assume that $S=\bSeif{\brep{b}}$, for some 
embedded \nband{X}word $\brep{b}$ of length $k$.
Let there be $m\le k$ negative bands in $\brep{b}$.
If $m=0$, we are done.  Otherwise, let $b(s)$ be negative,
say $b(s)=\brgen{x_p,x_q}^{-1}$, with $\{x_p<x_q\}\sub X$.  
Let $X'\isdefinedas X\cup\{x_1',x_2',x_3',x_4'\}$,
where $x_p<x_1'<x_2'<x_3'<x_4'<x_q$.
Let $b'(t)=\brinj{X,X'}b(t)$ for $1\le t<s$,
$b'(s)=\brgen{x_3',x_q}$,
$b'(s+1)=\brgen{x_2',x_4'}$,
$b'(s+2)=\brgen{x_1',x_3'}$,
$b'(s+3)=\brgen{x_p,x_2}$,
$b'(s+4)=\brgen{x_1',x_4'}$, and
$b'(t)=\brinj{X,X'}b(t+4)$ for $s+4< t\le k+4$
(see \Figref{trefins} for a rendition of this 
operation in terms of fence diagrams).
\begin{figure}
\centering
\includegraphics[width=.6\textwidth]{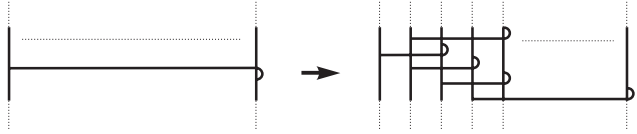}
\caption{``Trefoil insertion'' on a fence diagram
preserves the Seifert form while eliminating 
a negative band.\label{trefins}}
\end{figure}
Manifestly $\bSeif{\brep{b}}$ is diffeomorphic to 
$\bSeif{\brep{b}'}$ by a diffeomorphism that is the
identity of a single \nhandle{1}, 
and an easy calculation shows that
they have identical Seifert matrices.  There 
are only $m-1$ negative bands in $\brep{b}'$,
so this result of \citet{Rudolph1983c} is true
by induction on $m$.
\end{proof}

\begin{corollary}\label{S-equivalence and qp}
\begin{inparaenum}
\item\label{qp in every S-equivalence class}
Every $S$-equivalence class of knots contains 
strongly quasipositive knots.
\item\label{doubled-delta and qp}
Every knot can be converted to a strongly quasipositive knot
by a sequence of ``doubled-delta moves''.
\end{inparaenum}
\end{corollary}

\begin{figure}
\centering
\includegraphics[width=.6\textwidth]{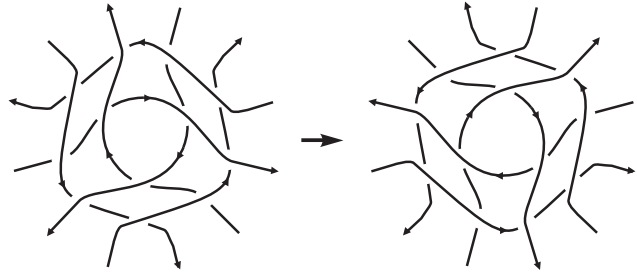}
\caption{The ``doubled-delta'' move.\label{doubdelt}}
\end{figure}

\begin{proof}
Knots $K_0$ and $K_1$ are \moredef{$S$-equivalent}%
    \index{S-equivalence@$S$-equivalence}%
\index{h-equivalence@$h$-equivalence|see {$S$-equivalence}} %
\citep[in the first instance,][``$h$-equivalent'']{Trotter1962} just 
in case $K_i$ has a Seifert surface $S_i$ such that
some homology bases of $S_0$ and $S_1$ produce identical
Seifert matrices; so (\ref{qp in every S-equivalence class})
is immediate from (the proof of) \eqref{qp is Seifert-undistinguished}.
\citet{NaikStanford2003} introduced the ``doubled-delta move'',
an operation on standard link diagrams illustrated in \Figref{doubdelt},
and proved that the equivalence relation on knots generated by 
applying such moves is precisely $S$-equivalence; 
(\ref{doubled-delta and qp}) follows directly.
\end{proof}

\begin{question}
\label{rotant/mutant question}
\citet{Conway1970} introduced an operation on knots and links
called \bydef{mutation}%
    \index{link(s)!mutation}.  The
result $L'$ of applying mutation to a link $L$
is called a \moredef{mutant}%
    \index{link(s)!mutant} of $L$.  Although mutant links
can be non-isotopic, they are indistinguishable by a wide 
variety of link invariants.  
Mutation was generalized by \citet*{AnsteePrzytyckiRolfsen1989}.
The result $L'$ of applying the operation defined
by \citeauthor*{AnsteePrzytyckiRolfsen1989} to $L$ is called a 
\bydef{rotant}%
    \index{link(s)!rotant} of $L$.  
The doubled-delta move is a special case of 
\citeauthor*{AnsteePrzytyckiRolfsen1989}'s operation.  
Obviously \eqref{S-equivalence and qp}(\ref{doubled-delta and qp}) 
can be rephrased as ``the doubled-delta move need not preserve
strong quasipositivity''.  
\ask{In what circumstances is a rotant
of a strongly quasipositive link strongly quasipositive?}  
In particular, \ask{is a mutant of a strongly quasipositive link 
necessarily strongly quasipositive?}
\end{question}

\subsubsection{Characterization of strongly quasipositive links}
\label{constructions III}

There is a simple characterization of quasipositive 
Seifert surfaces, and thus---in some sense---of strongly
quasipositive links.

\begin{unproved}{theorem}\label{charthm}
A Seifert surface $S$ is quasipositive if and only 
if, for some $n\in\Nstrict$, $S$ is ambient isotopic to a full 
subsurface of the fiber surface of the torus link
$O\{n,n\}$.
\end{unproved}

\begin{questions}\label{qp algorithms questions}
\ask{Does there exist an algorithm to determine whether 
or not a given Seifert surface is quasipositive?} 
\ask{Does there exist an algorithm to determine whether 
a given link is strongly quasipositive?}
\end{questions}

\subsubsection{Quasipositive annuli}
\label{constructions IV}

\begin{lemma}
Let $(L,g)$ be a framed link such that the annular 
Seifert surface $\AKn{L}{g}$ is quasipositive.
If the framing $f$ of $L$ is less twisted than $g$,
then $\AKn{L}{f}$ is quasipositive.
\end{lemma}

\begin{proof}
This is the ``Twist Insertion Lemma'' of \citet{Rudolph1983a}
(with the sign conventions in force here, it should be called
the ``Twist Reduction Lemma'').  A proof using fence diagrams 
is shown in \Figref{twistins}.
\begin{figure}
\centering
\includegraphics[width=.6\textwidth]{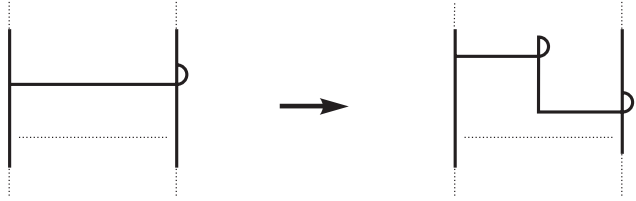}
\caption{``Twist reduction'' on a fence diagram.\label{twistins}}
\end{figure}
\end{proof}

\begin{definition} The \moredef{modulus of quasipositivity}%
    \index{modulus of quasipositivity ($\modqp{\phold}$)}%
    \index{quasipositivity!modulus of $\modqp{\phold}$}%
    \index{q()@$\modqp{\phold}$ (modulus of quasipositivity)} of a knot
$K$ is 
\begin{equation}\label{modqp defeq}
\modqp{K}\isdefinedas\sup\{t\in\Z\Suchthat \AKn{K}{t} 
\text{ is a quasipositive Seifert surface}\}.
\end{equation}
\end{definition}

\citet{Rudolph1990} applied \eqref{MFW--braids} 
(via \eqref{congruence theorem}) and \eqref{charthm}
to deduce the following.

\begin{unproved}{proposition}
If $\AKn{L}{f}$ is quasipositive, 
then $\ord_v\{L,f\}\ge 0$. 
\end{unproved}

\begin{unproved}{corollary}
If $K$ is a knot, then
$\modqp{K}\le \frac12\ord_v\{K\}\le -1-\deg_aF^*_K(a,x)$.
\end{unproved}

\begin{unproved}{corollary}\label{q(K) is an integer}
If $K$ is a knot, then $\modqp{K}\in\Z$.
\end{unproved}
Since $\modqp{K}=\TB(K)$ \citep{Rudolph1995}, 
\eqref{q(K) is an integer} also
follows from \citet{Bennequin1983}.

\begin{unproved}{corollary}
$\modqp{O}=-1$.
\end{unproved}

\begin{proposition}\label{modqp altdef}
For $K\ne O$, $\modqp{K}=\sup\{f\in\Z\Suchthat \text{ the framed knot }
(K,f) \linebreak[3]\text{embeds on the fiber surface of } O\{n,n\}
\text{ for some }n\in\Nstrict\}$.
\end{proposition}
\begin{proof}
If a framed knot $(K,f)$ embedded on a Seifert surface $S$
and $K$ is not full on $S$, then $(K,f)=(O,0)$; so \eqref{modqp altdef}
follows immediately from the characterization theorem \eqref{charthm}.
\end{proof}

\begin{unproved}{corollary}\label{qpSs is incompressible}
A quasipositive Seifert surface is incompressible.
\end{unproved}

Of course \eqref{qpSs is incompressible} also follows from 
\eqref{unknotting number of link of singularity}.  Yet another
proof can be extracted from \citet{Bennequin1983}.

\subsubsection{Quasipositive plumbing}
\label{constructions V}

A somewhat lengthy but straightforward combinatorial proof
of the following theorem appears in \citet{Rudolph1998}.

\begin{unproved}{theorem}\label{plumbing preserves qp} 
A Murasugi sum is quasipositive 
if and only if its summands are quasipositive.
\end{unproved}
 
\begin{unproved}{corollary}
\label{annular plumbings}
Iterated Murasugi sums of quasipositive annuli are quasipositive
Seifert surfaces.  
\end{unproved}
 
\begin{examples}
\eqref{annular plumbings} covers a surprising amount of ground.
\begin{asparaenum}[\bf \theexamples.1.]
\item\label{stably positive Hopf-plumbed stuff}
A \moredef{positive Hopf-plumbed fiber surface} is an iterated Murasugi
sum of positive Hopf annuli $\AKn{O}{-1}$ starting from the trivial
Murasugi sum $D^2$.  A fiber surface $S$
is \moredef{stably positive Hopf-plumbed} in case some 
positive Hopf-plumbed fiber surface $F$ can also be constructed
as an iterated Murasugi sum of positive Hopf annuli $\AKn{O}{-1}$ 
starting from $S$.  A remarkable theorem of \citet{Giroux2002}
states that a fibered link $L$ has a stably positive Hopf-plumbed
fiber surface if and only if a certain contact structure on $S^3$
(constructed from the fibration in a way he describes) is 
the standard contact structure.  In combination with 
\eqref{charthm} and \eqref{plumbing preserves qp}, 
\citeauthor{Giroux2002}'s 
theorem implies that a fibered link is strongly quasipositive
if and only if it is stably positive Hopf-plumbed.
\item\label{quasipositive fibered baskets}
A \bydef{basket} is a Seifert surface produced by repeatedly plumbing
unknotted annuli $\AKn{O}{k_i}$ to a single fixed $D^2\sub S^3$.
\citep[See][for details and examples.]{Rudolph2001a}  A fundamental 
theorem of \citet{Gabai1983b,Gabai1985} implies that a basket is a fiber 
surface if and only if it is Hopf-plumbed (allowing both positive
and negative Hopf annuli as plumbands).  In particular, according to
\citet{Rudolph2001a}, a fiber surface is a quasipositive
basket if and only if it is isotopic to a braided Seifert
surface $\bSeif{\brep{b}}$ where $\brep{b}$ is a strictly 
$\Ts$-positive \nband{\Ts}word for some espalier $\Ts$.
\item\label{strongly quasipositive special arborescent links}
General \moredef{arborescent links}%
    \index{link(s)!arborescent}%
    \index{arborescent link} are constructed as boundaries
of unoriented (possibly non-orientable) \nmanifold{2}s formed
by \moredef{unoriented plumbing}%
   \index{plumbing!unoriented}, in which not just
unknotted annuli $\AKn{O}{k_i}$ but also unknotted 
M\"obius bands are used as building blocks, while their
mode of assembly is suitably restricted \citep[it is coded 
by a planar tree with integer weights; see]
[]{Conway1970,Siebenmann1975,Gabai1986}.  So-called
\moredef{special arborescent links}%
    \index{link(s)!arborescent!special}%
    \index{arborescent link!special} \citep{Sakuma1994} are
those obtained by disallowing M\"obius plumbands.
As noted by \citet{Rudolph2001a}, the surfaces defining
special arborescent links are baskets; in particular, when
all weights of an arborescent presentation of an arborescent
link $L$ are strictly negative even integers, then $L$ is
strongly quasipositive.
\item\label{Whitehead doubles}
The \moredef{positively-clasped, $k$-twisted} 
\bydef{Whitehead double} $\mathscr{D}(K;k;{+})$ 
of a knot $K$ is defined
to be the boundary of the non-trivial plumbing
(that is, not a boundary-connected sum) of $\AKn{K}k$ 
and $\AKn{O}{-1}$ (\Figref{double figure} illustrates
this operation with a fence diagram).  
By \eqref{annular plumbings} and \eqref{charthm}, 
$\mathscr{D}(K;k;{+})$ is strongly quasipositive
if and only if $k\le\TB(K)$ \citep[see][]{Rudolph1993}.
\begin{figure}
\centering
\includegraphics[width=.3\textwidth]{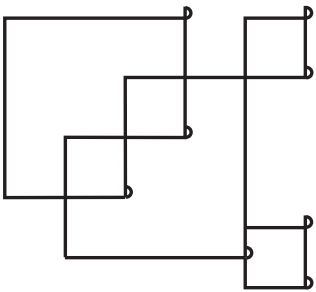}
\caption{A fence diagram for the plumbed Seifert surface
of $\mathscr{D}(O\{2,3\};0;{+})$.\label{double figure}}
\end{figure}
\item\label{charthm (bis)}
If a set of Seifert surfaces 
contains $\AKn{O}{-1}$ and
is closed under isotopy, Murasugi sum, and 
the operation of passing from a surface to 
a full subsurface,
then it contains all quasipositive Seifert surfaces; and 
the set of all quasipositive Seifert surfaces is the
smallest set with these properties.
\end{asparaenum}
\end{examples}

\subsubsection{Positive links are strongly quasipositive}
\label{positive links}
A standard link diagram $\diagram{L}$ is called \moredef{positive}%
    \index{link diagram!standard!positive}%
    \index{positive!standard link diagram} in case every crossing 
$\Seifx\in\SeifX{\diagram{L}}$ is positive; a link $L$ is 
\moredef{positive}%
    \index{positive!link}%
    \index{link(s)!positive} if $L$ has a positive standard link
diagram.  \citet{Nakamura1998,Nakamura2000} and \citet{Rudolph1999b}
independently proved the following.

\begin{unproved}{theorem}
If the standard link diagram $\diagram{L}$ is positive,
then the diagrammatic Seifert surface $\Seifsurf{L}$
is quasipositive.
\end{unproved}

\begin{unproved}{corollary}
A positive link is strongly quasipositive.
\end{unproved}

\begin{questions}
Two questions posed by \citet{Rudolph1999b} remain open and
are relevant here.
\begin{asparaenum}[\bf \thequestions.1.]
\item
Can positive links be characterized as strongly quasipositive links that 
satisfy some extra geometric conditions? 
\item
Does $K\ge O\{2,3\}$ imply that $K$ is strongly quasipositive?
(Here $K_1\ge K_2$ means that $K_1$ is concordant to $K_2$ 
inside a $4$-manifold with positive intersection form; this
partial order was defined and studied by \citet{CochranGompf1988}, 
who showed that if $K$ is a positive knot, then $K\ge O\{2,3\}$.)
\end{asparaenum}
\end{questions}

\subsubsection{Quasipositive pretzels}
\label{quasipositive pretzels}

\begin{figure}
\centering
\includegraphics[width=.9\textwidth]{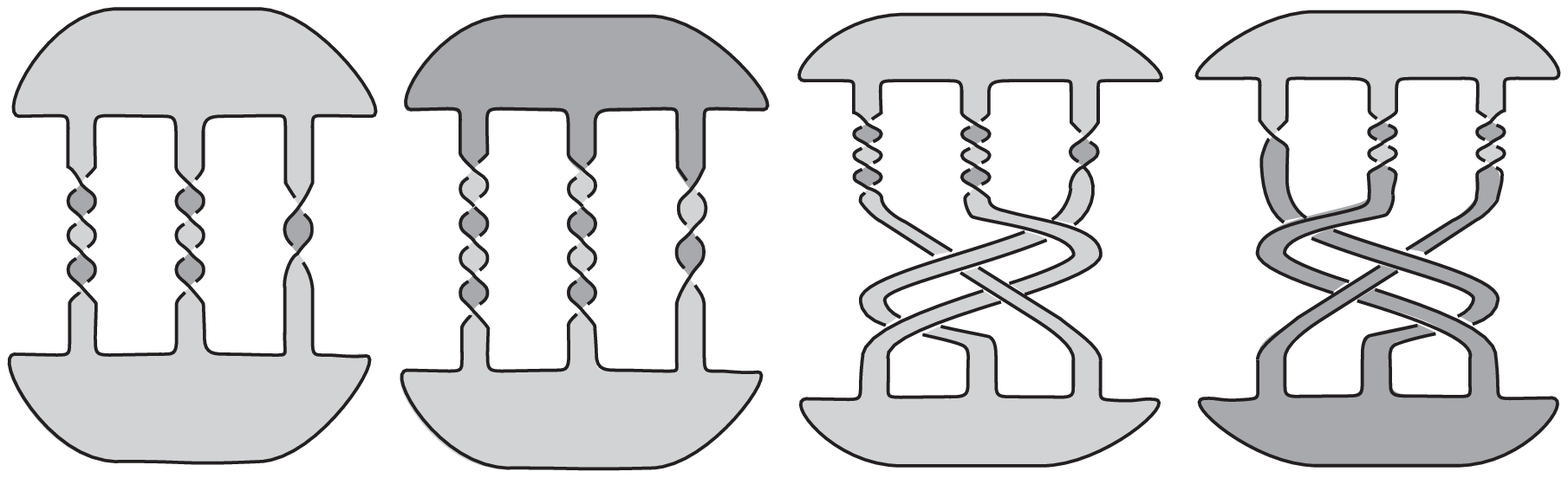}
\caption{The quasipositive pretzel surfaces $\pretzel{4,4,-2}$
and $\pretzel{5,5,-3}$; the braidzel surfaces
$\pretzel{\brgen{1},\brgen{2} \protect\mathstrut^{-1},
\brgen{1},\brgen{2} \protect\mathstrut^{-1};4,4,-2}$
and $\pretzel{\brgen{2} \protect\mathstrut^{-1},\brgen{1},
\brgen{2} \protect\mathstrut^{-1},\brgen{1};1,-3,-3}$.
\label{pretzel and braidzel figure}}
\end{figure}

\moredef{Braidzels}%
    \index{braidzel} are generalizations of the well-known
(oriented) pretzel surfaces%
    \index{pretzel!surfaces}%
    \index{link(s)!pretzel}, with braiding data supplementing 
the twisting data that specifies an ordinary pretzel surface 
(a braidzel with trivial braiding); two typical pretzel surfaces
and two typical braidzel surfaces are pictured in 
\Figref{pretzel and braidzel figure}.  Braidzels were defined
by \citet{Rudolph2001b} and have been further studied by 
\citet{Nakamura2004}, who showed that every link has a Seifert
surface which is a braidzel.
Using braidzels, \citet{Rudolph2001b} proved the following.

\begin{unproved}{proposition}
\label{quasipositive pretzel surface NASC}
The oriented pretzel surface 
$\pretzel{t_1,\dots,t_k}$ is quasipositive if and only if
the even integer $t_i+t_j$ is less than $0$ for 
$1\le i<j\le k$.
\end{unproved}

The question \eqref{rotant/mutant question} was somewhat
motivated by \eqref{quasipositive pretzel surface NASC},
since up to repeated mutation the pretzel link%
    \index{pretzel!link}%
    \index{link(s)!pretzel} $\Bd{\pretzel{t_1,\dots,t_k}}$
depends only on $\{t_1,\dots,t_k\}$.

\begin{question}
What are necessary and sufficient conditions on the defining data
of a braidzel that it be quasipositive?
\end{question}

\subsubsection{Links of divides}
\label{links of divides}

Identify $\C^2$ with the tangent bundle of $\C$.
\citet{A'Campo1998} calls a normal collection 
$P$ of edges in $D^2$ a \bydef{divide}, and constructs 
the \moredef{link of the divide $P$}%
    \index{link(s)!of a divide}%
    \index{divide!link of} as 
$L(P)\isdefinedas \{(z,w)\in S^3 \Suchthat z \in P, w \in T_z(P)\}$.  
This construction has been extended---first by allowing non-proper
edges \citep[``free divides'';][]{GibsonIshikawa2002a}, 
more generally by allowing immersed unoriented circle components 
\citep{Kawamura2002c}, and more generally yet 
(and most recently) by allowing immersions of graphs that are not
\nmanifold{1}s \citep[``graph divides'';][]{Kawamura2004}.  

\begin{unproved}{theorem}
If $P$ is a divide, then:
\begin{inparaenum}
\item
$L(P)$ is fibered, and every link $L$ of a singularity 
is of the form $L(P)$ for some divide $P$ \textup{\citep{A'Campo1998}};
\item
$L(P)$ is strongly quasipositive \textup{\citep{Hirasawa2000}};
\item
$L(P)$ is stably positive Hopf-plumbed \textup{\citep{Hirasawa2002}}; 
and indeed
\item
$L(P)$ is positive Hopf-plumbed \textup{\citep{Ishikawa2002}}.
\suspend{inparaenum}
If $P$ is a free divide, possibly with circle components, then
\resume{inparaenum}
\item
$L(P)$ is strongly quasipositive \textup{\citep{Kawamura2002c}}.
\suspend{inparaenum}
If $P$ is a graph divide, then
\resume{inparaenum}
\item
$L(P)$ is quasipositive \textup{\citep{Kawamura2004}}.
\end{inparaenum}
\end{unproved}

See also \citet*{%
CouturePerron2000,
QuachWeber2000,
Ishikawa2001,
GibsonIshikawa2002b,
Gibson2002,
GodaHirasawaYamada2002,
Chmutov2003}.

\subsection{Non-strongly quasipositive links}

Other than the tautologous construction by forming the closed
braid of a quasipositive braid, little is yet known about
systematic constructions of non-strongly quasipositive transverse
{\Clink}links.

\begin{example}\label{fibered qp but not strongly qp}
There exist fibered links which are quasipositive but
not strongly quasipositive (compare with 
\eqref{stably positive Hopf-plumbed stuff}).
Examples include all the links 
$H_{n+m,m}$ obtained from $O\{n+2m,n+2m\}$ ($m, n\ge 1$)
by reversing the orientation of (any) $m$ components.
\Figref{H_3_1 figure} pictures 
a ribbon surface $R=f(S^1\times \clint01 \coprod D^2)$ 
for $H_{1,1}$ and a fence diagram for the fiber surface
$\smooth{R}$ of $H_{1,1}$.
\begin{figure}
\centering
\includegraphics[width=.6\textwidth]{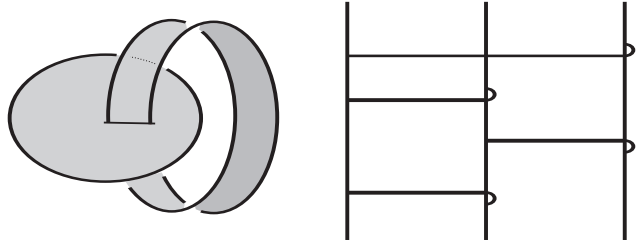}
\caption{A quasipositive fibered link which is not 
strongly quasipositive.\label{H_3_1 figure}}
\end{figure}
\end{example}

\begin{example}
\label{arborescent qp but not special arborescent}
There exist arborescent links which are quasipositive but
are not special arborescent and not strongly quasipositive 
(compare with \ref{strongly quasipositive special arborescent links}).
An example is the quasipositive arborescent ribbon knot pictured,
along with its defining weighted tree, in \Figref{loveknot figure}.
\begin{figure}
\centering
\includegraphics[width=.6\textwidth]{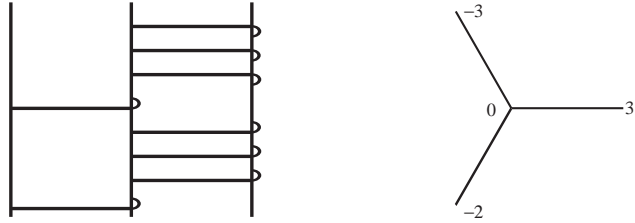}
\caption{Two representations of a quasipositive Lover's Knot.
\label{loveknot figure}}
\end{figure}
\end{example}

\begin{question}
Every pretzel link is arborescent.  (Caution: this does not mean 
that every pretzel link bounds an oriented pretzel surface, and in
fact there is a natural sense in which those which do are a small
minority.  Nor is a typical oriented pretzel surface itself an
arborescent Seifert surface.)  What are necessary and sufficient
conditions on an arborescent link that it be quasipositive, or
strongly quasipositive?  
\end{question}

Besides \eqref{quasipositive pretzel surface NASC},
some partial but systematic progress towards answering this 
question has been made by \citet{Tanaka1998b} (for the 
of rational links).

%
%
\section{Complex plane curves in the small and in the large}\label{alggeom}

As noted in \eqref{polynomials are enough}, every transverse
{\Clink}link $L$ is ``algebraic'' in the sense that, up to 
isotopy, $L$ is the transverse intersection of 
$rS^3=\{(z,w)\in\C^2\Suchthat \norm{(z,w)}=r\}$ and 
a complex algebraic curve $\variety{f}$;
clearly the same class of links arises whether $r>0$ is
fixed once for all, or allowed to vary throughout $\Rpos$,
so long as $f$ is allowed to vary throughout $\C[z,w]$.
By contrast, interesting additional restrictions on the 
isotopy type of $L$ arise in case, for each fixed $\variety{f}$,
$r$ is constrained to be either ``very small'' or ``very large''.
More precisely, the situation is as follows
\citep[essentially this is proved by][in general dimensions]{Milnor1968}.
Let $F\in\C[z,w]$ be non-constant and without repeated factors.
The \moredef{argument}
    \index{argument (of a polynomial)} of $F$ is 
$\argument{F}\from\diff{\C^2}{\variety{F}} \to S^1$.
The \moredef{Milnor map of $F$}
    \index{Milnor!map of a polynomial $F$ ($\Mmap{F}$)} is
    \index{Milnor!map of a polynomial $F$ ($\Mmap{F}$)!at radius $r$ ($\Mmapr{F}{r}$)} is
$\Mmap{F}\isdefinedas\arg(F)\restr(\diff{S^{2n-1}}{\divisor{F}})$; 
for $r>0$, the \moredef{Milnor map of $F$ at radius $r$}, 
denoted by $\Mmapr{F}{r}$, is the Milnor map of 
$(z,w)\mapsto F(z/r,w/r)$.
Let $m(F)\isdefinedas\inf\{\norm{(z,w)}^2\Suchthat F(z,w)=0\}
=\sup\{r\Suchthat rS^3\cap\variety{F}=\emptyset\}$.

\begin{unproved}{proposition}\label{Milnor maps proposition}
There is a finite set $\Xs(F)\sub\Rg{m(F)}$ of radii $r$
with the following properties.
\begin{inparaenum}
\item\label{almost all L(F,r) are links}
If $r\in \diff{\Rg{m(F)}}{\Xs(F)}$, 
then $\variety{F}$ intersects $rS^3$ transversally, 
so that $L(F,r)\isdefinedas (1/r)(\variety{F}\cap rS^3)$ 
is a link in $S^3$.
\item\label{isotopic Milnor links}
If $r$ and $r'$ are in the same component 
of $\diff{\Rg{m(F)}}{\Xs(F)}$ then 
$L(F,r)$ and $L(F,r')$ are isotopic. 
\item\label{Milnor maps are simple and adaptable} 
If $m(F)<r\notin \Xs(F)$, then there is 
a trivialization $n\from L(F,r)\times\C\to\normbun{L(F,r)}$,
as in \eqref{normal bundle of a link}, that lies 
in the homotopy class corresponding to any Seifert
surface for $L(F,r)$ and such that $\Mmapr{F}{r}$
is adapted to $\Mmapr{F}{r}$ with $d(K)=1$ for all 
$K\in\compsof{L(F,r)}$.  
\end{inparaenum}
\end{unproved}

\begin{definitions}
In case $F(0,0)=0$, for any $r\in\opint{0}{m(F)}$ the transverse
\Clink{link} $L(F,r)$ is called the 
\moredef{link of the singularity} of $F$ (or of $\variety{F}$)
at $(0,0)$.  In any case, for any $r>\max\Xs(F)$ the transverse
\Clink{link} $L(F,r)$ is called the 
\moredef{link at infinity} of $F$ (or of $\variety{F}$).
\end{definitions}

\subsection{Links of singularities as transverse 
\texorpdfstring{$\C$}{C}-links}
\label{links of singularities}
It can happen that $\Mmapr{F}{r}$ has 
degenerate critical points, and so is not a 
Morse map, or that $\Mmapr{F}{r}$ is a Morse map
but not a fibration (the latter in fact being common: 
as noted in the first paragraph of this section, every
transverse \Clink{link} can be realized as $L(F,r)$ 
for some $F$ and $r>0$, while according to the previous
section many, many transverse \Clink{link}s are not
fibered links).  However, \citet{Milnor1968} proved
the following (and its analogue in higher dimensions).

\begin{unproved}{theorem}
\label{Milnor maps and fibrations}
If $m(F)<r<\min\Xs(F)$, then $L(F,r)$ is a fibered link
and $\Mmapr{F}{r}$ is a fibration.
\end{unproved}

In other words, if $L$ is a link of a singularity, then
$L$ is fibered.

Let $L$ be the link of a singularity.
There is a huge literature devoted to studying the
algebraic and geometric topology $L$, and especially
of its fibration (which, at the geometric level,
completely determines $L$ and all its invariants---though
not necessarily in a perspicuous fashion).  Some starting 
points for the dedicated reader are 
\citet{EisenbudNeumann1985,Durfee1999,Neumann2001,LDT2003}.
Here it may simply be noted that many---though 
certainly not all---of the interesting topological
features of the fibration of $L$ can profitably be 
explored using one or another of the representations
of $L$ alluded to in earlier sections.  On one hand,
$L$ is the link of a divide \citep{A'Campo1998},
from which one representation of $L$ as a positive Hopf-plumbed
link can be derived and exploited \citep{A'Campo1998}.
On the other hand, it is easy to use labyrinths to
see that $L$ is a positive closed braid (and not overwhelmingly
difficult actually to derive a specific positive braidword for 
$L$), whence by \citet{Rudolph2001a} one may derive (and,
at least potentially, exploit) an apparently different---though
surely closely related---representation of $L$ as a positive
Hopf-plumbed link.

\subsection{Links at infinity as transverse \texorpdfstring{$\C$}{C}-links}
\label{links at infinity}

Let $L$ be a link at infinity.  In general, $L$ need not
be fibered (for instance, every unlink $O^{(n)}$ is a link
at infinity, but only $O=O^{(n)}$ is fibered), but as
\citet{Neumann1989b} puts it, $L$ is ``nearly'' fibered.
The literature on links at infinity is smaller than that
on links of singularities, but still too large and 
multifaceted to summarize
here.  The reader is referred to \citet{BoileauFourrier1998} 
for a review up to 1998.  Some more recent articles on 
various aspects of the subject are 
\citet*{Bodin1999,Nemethi1999,Neumann1999,
ArtalCassou-Nogues2000,NeumannNorbury2000,
PaunescuZaharia2000,GwozdziewiczPloski2001,
GwozdziewiczPloski2002,Rudolph2002,
NeumannNorbury2003,Cimasoni2004}.

%
%

\section{Totally tangential \texorpdfstring{$\C$}{C}-links}\label{totaltan}

\begin{definition}\label{def:totally tangential C-link} 
A link $L\sub S^3$ is a 
\moredef{totally tangential {\Clink}link}%
    \index{totally tangential C-link@totally tangential {\Clink}link}%
    \index{C-link@{\Clink}link!totally tangential} provided that
$L=\variety{g}\cap S^3$,
where
$g\in\holo{\C^2}$ is such that 
$\variety{g}\cap\Int{D^4}=\emptyset$
and $\variety{g}\cap S^3$ is a non-empty 
non-degenerate critical manifold of index $1$ of 
$\rho\restr\reg{\variety{g}}$,
with $\rho(z,w)\isdefinedas\norm{(z,w)}^2$.
\end{definition}

This section follows \citet{Rudolph1992c,Rudolph1995,Rudolph1997}.

\begin{unproved}{lemma}
\label{totally tangential iff Legendrian and analytic}
$L$ is a totally tangential {\Clink} if and only if 
\begin{inparaenum}
\item
$L$ is Legendrian with respect to the standard contact
structure on $S^3$, and
\item
$L$ is an $\R$-analytic submanifold of $S^3$.
\end{inparaenum}
\end{unproved}

By \ref{totally tangential iff Legendrian and analytic},
a totally tangential \Clink{link} has a natural framing $f=\Legf{L}$.

\begin{unproved}{proposition}\label{generalizing ttclinks}
Let $\rho\from U\to\R$ be an exhausting strictly plurisubharmonic 
function on an open set $U\sub C^2$ such that 
$D=\rho^{-1}(\Rle{1})$ is a Stein disk, so 
$\sph\isdefinedas\Bd{D}$ is a strictly pseudoconvex 
\nsphere{3} with a natural contact structure.
If $g\in\nizholo{D}\sub\holo{U}$ is such that 
$\variety{g}\cap\sph$ is a non-empty 
non-degenerate critical manifold $M$ of index $1$ of 
$\rho\restr\reg{\variety{g}}$,
then $M$ is a Legendrian \nsubmanifold{1} of 
the contact manifold $\sph$.  Further, there
exists a totally tangential \Clink{link} $L\sub S^3$ 
such that $(\sph,M)$ is diffeomorphic 
to $(S^3,L)$ by a diffeomorphism carrying 
the natural normal line field of $M$ in $\sph$
to the natural normal line field of $L$ in $S^3$.
\end{unproved}

\begin{question}\label{transverse C-links question}
\ask{Does every totally tangential {\Clink}link arise (up to isotopy)
as in \eqref{def:totally tangential C-link} with $g$ required
to be a polynomial?}  \citep[][shows that $g$ may be 
be required to be entire.]{Rudolph1995}
\end{question}

\begin{theorem}
\label{Legendrian isotopic etC}
If $L\sub S^3$ is Legendrian, then $L$ is isotopic through
Legendrian links to a totally tangential \Clink{link} $L'$
with $\tb(L')=\tb(L)$. 
\end{theorem}
\begin{proof}
The proof for knots \citep{Rudolph1995} extends directly
to links.
\end{proof}

\begin{unproved}{corollary}
\label{connection to transverse C-links}
If $L\sub S^3$ is Legendrian, then there is a polynomial map
$p\from\C^2\to\C$ such that
\begin{inparaenum}
\item
$\sing{\variety{p}}=\emptyset$,
\item
$\variety{p}$ intersects $S^3$ transversely,
\item
$\variety{p}\cap D^4$ is an annular surface, and 
\item
$\variety{p}\cap S^3 = \Bd\AKn{L}{\tb{L}}$.
\end{inparaenum}
\end{unproved}

Coupled with 
\eqref{K-M theorem}\thru\eqref{unknotting number of link of singularity},
the connection \eqref{connection to transverse C-links} between
totally tangential \Clink{links} and transverse \Clink{links},
leads to several interesting conclusions.

\begin{unproved}{theorem}
If $\TB(K)\ge 0$, then $K$ is not slice.
\end{unproved}

\begin{figure}
\centering
\includegraphics[width=2in]{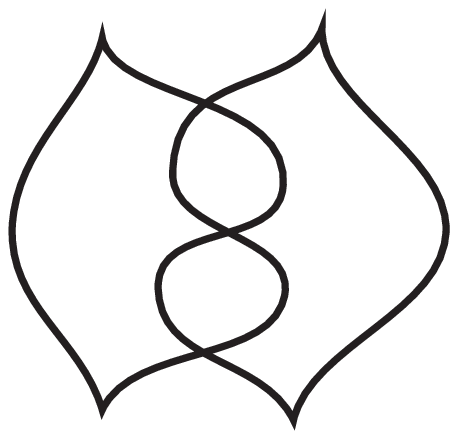}
\caption{A front.\label{front figure}}
\end{figure}
A \bydef{front} \citep[see][]{Eliashberg1993,Etnyre2004} is a 
piecewise-smoothly immersed curve $\gamma\sub\R^2$ 
of the sort illustrated in Figref{front figure}: 
$\gamma$ has finitely many cuspidal corners and
finitely many normal crossings; 
the corners are exactly the local extrema of $\pr_2\restr\gamma$; 
at a corner, the common tangent line of the two branches
is parallel to $\{0\}\times\R$; and elsewhere the tangent lines 
are parallel neither to $\{0\}\times\R$ nor to $\R\times\{0\}$.
A \moredef{front diagram}%
    \index{front!diagram}%
    \index{diagram!front}%
    \index{diagram!front ($\frdiagram{\phold}$)}%
    \index{FD(L)@$\frdiagram{\phold}$ (front diagram)}%
    \index{diagram!front ($\frdiagram{\phold}$)!information of ($\frinfo{\phold}$)}%
    \index{diagram!front ($\frdiagram{\phold}$)!picture of ($\frpict{\phold}$)}%
    \index{I@information!about a front ($\brinfo{\phold}$)}%
    \index{P@picture!of a front ($\brpict{\phold}$)} is a pair
$\frdiagram{L}\defines(\frpict{L},\frinfo{L})$, where
$\frpict{L}$ is a front, and $\info{L}$ is the information 
that $L$ is a closed Legendrian \nsubmanifold{1} of $\R^3$ 
equipped with a certain standard contact structure (consistent 
with its embedding as the complement of a point in $S^3$ equipped
with its standard contact structure).

\begin{unproved}{lemma}
\begin{inparaenum}
\item
Given any front $\frpict{L}$, the information $\info{L}$ is
true, and $\frdiagram{L}$ determines $\unorient{L}$ up to a 
vertical translation in $\R^3$.
\item
A standard link diagram $\diagram{L}$ for $L$ can be derived
from $\frpict{L}$ by 
\begin{inparaenum}
\item
letting $\pict{L}$ be $\frpict{L}$ with its cusps rounded off, and
\item
letting $\info{L}$ declare that at each crossing of $\pict{L}$,
the over-arc is the branch of which the tangent line has 
larger positive slope.
\end{inparaenum}
\end{inparaenum}
\end{unproved}

\begin{unproved}{proposition}\label{tb:Proposition 3} 
For any knot $K$, if $r\le \TB(K)$ then 
$\TB(D(K,r,{+}))\ge 1$.
\end{unproved}

\begin{unproved}{theorem}
For any knot $K\sub S^3$, $\genus{s}{K}\ge (\TB{K}+1)/2$.
\end{unproved}

%
%

\section{Relations to other research areas}
\label{applic}

This section is a bibliography of papers that relate
the knot theory of complex plane curves to other research
areas; in some of the subsections, one (or a few) paper(s)
are highlighted.  The choice of highlights, and indeed the
selections themselves, reflect both my bias and my ignorance; the 
latter, at least, is corrigible, and I welcome pointers to sources 
I have neglected.  The division into subsections is necessarily 
somewhat arbitrary.  

\subsection{Low-dimensional real algebraic geometry; 
Hilbert's 16th problem.}
An interesting series of recent papers 
\citep{Orevkov1999,
OrevkovPolotovskii1999,
Orevkov2000,
Orevkov2001a,
Orevkov2001b,
OrevkovShustin2002} has applied quasipositivity to 
Hilbert's $16$th problem
\hyperlink{Hilbert's $16$th problem}{Hilbert's $16$th problem}.
\subsection{The Zariski Conjecture; knotgroups of complex plane curves}
\label{knotgroups of plane curves}\label{Zariski conjecture}

\begin{unproved}{theorem}\label{Orevkov's Zariski theorem}
Let $f\in\C[z,w]$.  If 
\begin{inparaenum}
\item
$\variety{f}\sub\C^2$ is a node curve, and 
\item\label{positive closed braid at infinity}
the link at infinity of $\variety{f}$ is a positive 
\textup(and not merely quasipositive\textup) closed braid,
\end{inparaenum}
then the knotgroup $\fundgpnbp{\diff{\C^2}{\variety{f}}}$ is abelian.
\end{unproved}

\eqref{Orevkov's Zariski theorem} is the version of 
\hyperlink{Zariski's Conjecture}{``Zariski's Conjecture''}
proved by \citet{Orevkov1988}.  The 
original version of the conjecture (in which the 
link at infinity is required to be, not just any positive closed
braid, but actually $\closedbraid{\nabla_{\deg{f}}}$) had been
proved by \citet{Fulton1980,Deligne1981,Nori1983,Harris1986}.
The particular beauty of Orevkov's proof consists,
first, in his use of a 
labyrinth for $f$ (assumed to be in Weierstrass form) to 
construct a presentation of $\fundgpnbp{\diff{\C^2}{\variety{f}}}$ 
with very many generators ($\deg_w(f)$ generators for each 
\ndimensional{0} stratum of the labyrinth) and very many relations,
which are however all of very simple forms; and, second,
in his ingenious application of the clews of the labyrinth,
coupled with the positivity of the closed braid at infinity,
to inductively reduce that presentation to a smaller presentation
that visibly presents an abelian group.

The background and motivation for the conjecture, and
related material on knotgroups of complex plane curves, 
may be found in \citet*{Enriques1923,
Lefschetz1924,
Brauner1928,
Zariski1929,
vanKampen1933,
Burau1934,
Zariski1936,
Zariski1971,
Reeve1954,
Cheniot1973,
LDT1974b,
Oka1975,
Oka1974,
OkaSakamoto1978,
Chang1979,
Randell1980,
Kaliman1992,
MoishezonTeicher1996,
Kulikov1997,
Garber2003}.
The ideas of \citet{Orevkov1998} are extended by 
\citet*{Orevkov1990x,DethloffOrevkovZaidenberg1998,NeumannNorbury2003}. 

\subsection{Keller's Jacobian Problem; 
embeddings and injections of \texorpdfstring{$\C$}{C} 
in \texorpdfstring{$\C^2$}{C{2}}}
\label{Keller's Jacobian Problem}
For somewhat unclear reasons, 
\hyperlink{Keller's Jacobian Problem}{Keller's Jacobian Problem}
inspired \cite{Segre1956} to state---though not to prove 
correctly---a theorem classifying embeddings of $\C$ in $\C^2$,
which has now been proved and generalized to injections.

\begin{unproved}{theorem}\label{line embedding/injection}
Let $p,q\in\C[t]$ be such that 
$F\from\C\to\C^2\suchthat t\mapsto (p(t),q(t))$
is a \textup(set-theoretical\textup) injection.
\begin{inparaenum}
\item\label{A-M-S theorem}
If $F$ is a smooth embedding, then one of $\deg(p)$, $\deg(q)$ divides
the other; as a consequence, there is a polynomial change of coordinates
$Q\from\C^2\to\C^2$ such that $Q(F(\C))$ is a straight line.
\item\label{Z-L theorem}
If $F$ is not smooth, then either one of $\deg(p)$, $\deg(q)$ divides
the other or $\deg(p)$ and $\deg(q)$ are relatively prime; as
a consequence,  there is a polynomial change of coordinates
$Q\from\C^2\to\C^2$ such that $Q(F(\C))$ is one of the curves
$\{(z,w)\in\C^2\Suchthat z^m+w^n\}$ for some relatively prime
positive integers $m$ and $n$.
\end{inparaenum}
\end{unproved}

\eqref{line embedding/injection}(\ref{A-M-S theorem}) was first proved
by \citet{AbhyankarMoh1975} (algebraically) 
and by \citet{Suzuki1974} (analytically); a topological (knot-theoretical)
proof was given by \citet{Rudolph1982}.  
\eqref{line embedding/injection}(\ref{Z-L theorem}) was first 
proved by \citet{ZaidenbergLin1983} (analytically); a topological
(knot-theoretical) proof was given by 
\citet{NeumannRudolph1987}.  Many other proofs have appeared since.
Those proofs, related theorems, and some of the more knot-theoretical
work on the original Jacobian Problem may be found in
\citet*{AbhyankarSingh1978,
OkaSakamoto1978,
Orevkov1990a,
Orevkov1991,
Chang1991,
Kang1991,
ChangWang1993,
Kaliman1993,
Abhyankar1994,
Artal-Bartolo1995,
GurjarMiyanishi1996,
Orevkov1996,
A'CampoOka1996,
Essen1997,
Zoladek2003}.

\subsection{Chisini's statement; braid monodromy}
\label{braid monodromy}\label{Alexander invariants}\label{arrangements}
It was known early on \citep[see]{Zariski1929} by algebraic geometers 
that, for proposed applications to the study of algebraic surfaces,
the most important knotgroups of complex plane curves were those
of cusp curves.  Starting in \citeyear{Chisini1933} and continuing
through the 1950s, Oscar Chisini and his school
published a series of papers around that subject
\citep*{Dedo1950,
Chisini1952,
Tibiletti1952,
Chisini1954,
Tibiletti1955a,
Tibiletti1955b,
Tibiletti1955c,
Tibiletti1955d,
Chisini1955}.  The tools they
brought to bear were what would later be called ``braid monodromy''
and ``the arithmetic of braids''
\citep*{Moishezon1981,
Moishezon1983,
Moishezon1985,
MoishezonTeicher1988,
MoishezonTeicher1991,
Moishezon1994,
MoishezonTeicher1994a,
MoishezonTeicher1994b,
MoishezonTeicher1996}.

\begin{unproved}{theorem}\label{Chisini: yes and no}
Let $\brep{b}$ be a quasipositive cuspword in $\brgp{n}$ such
that $\braidof{\brep{b}}=\nabla_n$.
\begin{inparaenum}
\item\label{Chisini: local yes}
For any $n$ and any such $\brep{b}$, there exists
a holomorphic cusp curve in $D^2\times D^2$
that corresponds to $\brep{b}$ via an appropriate holomorphic map
$f\from D^2\to\monicpoly{n}$.
\item\label{Chisini: global no}
There exist $n$ and $\brep{b}$ such that no cusp curve in $\C^2$
corresponds to $\brep{b}$ via a polynomial map 
$f\from D^2\to\monicpoly{n}$.
\end{inparaenum}
\end{unproved}

\eqref{Chisini: yes and no}(\ref{Chisini: local yes}) is of course
a special case of \eqref{Orevkov's monodromy theorem}, as
proved (in much greater generality) by \citet{Orevkov1998}.
\eqref{Chisini: yes and no}(\ref{Chisini: global no}) was proved by
\citet{Moishezon1994}, in contradiction to 
``a statement of Chisini'' (\citeyear{Chisini1954}).

There is an enormous, interesting, and constantly growing literature
on braid monodromy and its many offsprings, including applications
to line arrangements, the Alexander invariants of complex plane curves,
and so on; see
\citet{%
Libgober1986,
LoeserVaquie1990,
Dung1994a,
Dung1994b,
DungHa1995,
CohenSuciu1997,
Dung1999,
ArtalCarmonaCogolludo2001,
KharlamovKulikov2001,
Dung2002,
KulikovKharlamov2003,
Yetter2003}.

\subsection{Stein surfaces}
\label{Stein surfaces as branched covers}

By \eqref{covering of Stein is Stein}, 
\eqref{quasipositive ribbons are holomorphic},
\eqref{ribbons are braidable}, and 
\eqref{ribbons are ribbons}, 
a branched covering of the bidisk branched over a 
quasipositive braided surface is diffeomorphic to a Stein surface.
\citet{LoiPiergallini2001} proved the following converse.

\begin{unproved}{theorem}
If $X$ is a Stein surface, then $X$ is diffeomorphic to 
a branched covering of $D^2\times D^2$ branched 
over a quasipositive braided surface.
\end{unproved}

%
%

\section{The future of the knot theory of complex plane curves}\label{future}
This section suggests some possible new directions for the
knot theory of complex plane curves.

\subsection{Transverse \texorpdfstring{$\C$}{C}-links and their Milnor maps.}

As observed before \eqref{Milnor maps and fibrations},
the Milnor map $\Mmapr{F}{r}$ at radius $r>m(F)$ of a 
non-constant reduced polynomial $F\in\C[z,w]$ can have 
degenerate critical points.  However, this can happen 
for only finitely many $r\in\Rg{m(F)}$.  Every
other $\Mmapr{F}{r}$ is a Morse map, and typically \shout{not}
a fibration \citep[even if $\variety{F}\cap rS^3$ is a fibered link;
see][]{HirasawaRudolph2003}.  These Morse maps deserve further
study.  

\begin{motto}\label{Milnor maps are interesting}
A Milnor map is interesting even if it is not a fibration.
\end{motto}

\begin{question} Here are a couple of specific questions about
Milnor maps that are Morse.  Is such a map necessarily topless?
Under what conditions is such a map minimal 
\citep*[that is, possessed
of exactly the minimum number necessary for any Morse map on
its domain; see][and references cited therein]%
{PajitnovRudolphWeber2001,HirasawaRudolph2003}?
\end{question}

Rational functions and non-reduced polynomials
also have Milnor maps.  \citeauthor{Cimasoni2004}'s
study (\citeyear{Cimasoni2004}) of Alexander invariants
of links of infinity of (reduced) polynomials leads him to 
investigate, in passing, the Milnor maps at infinity of 
non-reduced polynomials,
and \citet{Pichon2003} has begun the study of the Milnor map 
of rational functions at their points of indeterminacy, 
but a general study of these maps remains to be undertaken.

\subsection{Transverse \texorpdfstring{$\C$}{C}-links as links at infinity
in the complex hyperbolic plane.}

To date, with the exception of \citet{Rudolph2002}, all
investigations of links at infinity have implicitly or
explicitly concerned themslves with behavior at ``affine 
infinity'', that is, in affine space $\C^2$.  
Is that an oversight?

\begin{motto}\label{complex hyperbolic link at infinity}
A transverse {\Clink}link is a complex hyperbolic link at infinity.
\end{motto}

The meaning of \eqref{complex hyperbolic link at infinity} is
this.  As in the case of the real hyperbolic plane $\mathbb{H}^2(\R)$, 
there is a Klein model of the complex hyperbolic plane $\mathbb{H}^2(\C)$ 
which displays $\mathbb{H}^2(\C)$ as holomorphically equivalent to 
$\Int{D^4}\sub\C^2$.
From a purely topological---even complex-analytic---viewpoint, 
therefore, a link at infinity in the complex plane is nothing 
more (nor less) than a transverse \Clink{link} in $S^3$.  
But, like $\mathbb{H}^2(\R)$, $\mathbb{H}^2(\C)$ has a geometry
which makes it very different from the affine plane over the same
field.  \citet{Rudolph2002} exploited (in a very minor way)
the linear structure in $\mathbb{H}^2(\C)$.  What about
more subtle geometric features, like any of several different
variants on curvature?

\begin{question}\label{curvature question}
One way to interpret the main theorem of 
\citet{Osserman1981} is that, if $f\in\holo{\C^2}$ is entire,
then $\variety{f}$ has a well-defined link at infinity (i.e.,
is an algebraic curve) if and only if $\variety{f}$ has finite
total absolute curvature.  Given $f\in\holo{\mathbb{H}^2(\C)}$, 
possibly subject to further conditions (e.g., finite topological
type), are there reasonable geometric conditions---for instance, 
integral curvature conditons---on $\variety{f}$ that are equivalent 
to the existence of a well-defined transverse \Clink{link} ``at
infinity'' for $\variety{f}$?
\end{question}

\subsection{Spaces of \texorpdfstring{$\C$}{C}-links}

The sets of function defined in \dispeqref{function algebras}
can be used to describe sets of singular transverse and totally tangential
\Clink{links}, as follows.  In fact, it is obvious 
that the quotient set $(\diff{\holo{D^4}}{\{0\}})/\nzholo{D^4}$ 
is naturally identified with the set of ``singular transverse 
\Clink{link}s'' (with positive multiplicities), and likewise that 
the quotient set $\nizholo{D^4}/\nzholo{D^4}$ 
is naturally identified with the set of ``singular totally tangential 
\Clink{link}s''.  Moreover, both quotient sets are (more or less
naturally) partitioned by the severity of the ``singularity''
of the \Clink{link}.
It is perhaps less obvious what topologies are most appropriate
to turn these partitions (or suitable refinements thereof)
into good stratifications of some infinite-dimensional manifolds, 
in which (for instance) the \ncodimension{0} strata represent 
restricted isotopy types of ordinary non-singular transverse, or totally
tangential, \Clink{links}.  Nonetheless, I proclaim the following.

\begin{motto}
Naturally stratified spaces of singular {\Clink}links exists
and should be studied \foreign{\`a la} Vassiliev.
\end{motto}

A significant distinction between ``singular links'' in the present
style, and in the style of Vassiliev, is that because these
singular links are the $0$-loci of (restrictions of holomorphic)
functions $S^3\to\C$ rather than the images of functions 
$S^1\times\intsto{k}\to S^3$, at a generic bifurcation
(i.e., a \ncodimension{1} stratum in one of the proposed spaces of
singular links) the number of components changes (by one).

\subsection{Other questions.}\label{open questions}
\hypertarget{open questions}{}

Here are two last open questions ranging somewhat further afield
than the earlier ones.

\begin{question}
\label{Whitehead conjecture question}
It is a long-standing conjecture \citep[closely related 
to the Whitehead Conjecture; see][]{Howie1985} that,
if $S\sub D^4$ is a ribbon $2$-disk, then 
$\diff{D^4}{S}$ is aspherical.  A simple homological 
argument shows that this conjecture is equivalent to
the conjecture that $H^2(\widetilde{\diff{D^4}{S}};\Z)=0$, 
where for any $M$, $\widetilde{M}$ here denotes the universal 
covering space of $M$.  Now, as \citet{Rudolph1983b} points out,
if $S\sub D^4$ is any ribbon \ndisk{2}, then 
there is a piece of complex plane curve $\G\sub D^4$ (diffeomorphic
to $D^2$) such that $\diff{D^4}{S}$ is homotopy-equivalent to 
$\diff{D^4}{\G}$.
\ask{If a piece of complex plane curve $\G\sub D^4$ is diffeomorphic
to $D^2$, is $\diff{D^4}{\G}$ aspherical?}
Of course $\diff{D^4}{\G}$ is homotopy equivalent to $\Ext{\G}{D^4}$.  
For an appropriate choice of $\Nb{\G}{D^4}$, $\Ext{\G}{D^4}$ is a
Stein domain.  By \eqref{covering of Stein is Stein} and
\eqref{Cousin problem}, to answer 
this question affirmatively
it would suffice to show that every complex curve in 
the universal cover $X$ of $\Int{D^4}{\G}$ is
$\variety{f}$ for some $f\in\holo{X}$.
\end{question}

\begin{question}
\label{Smith conjecture question}
Call $K\sub S^3$ a \bydef{counter-Smith knot}%
    \index{knot!counter-Smith} in case $K$ is a non-trivial knot
and, for some $p>1$, the $p$-sheeted cyclic
branched cover of $S^3$ branched along $K$ is again $S^3$.
Assuming the Poincar\'e Conjecture, the Smith Conjecture 
is equivalent to the assertion that there are no counter-Smith
knots.  A \ndimensional{3} proof of the Smith Conjecture has 
of course been given by Thurston \foreign{et al.}; 
it is beautiful but lengthy, and depends on separate analyses
of a number of different cases.
In 1983 I noted that a simple appeal to Donaldson's results from
\ndimensional{4} gauge theory shows that, if the knot $K$ is a
non-slice transverse {\Clink}link, then $K$ is not counter-Smith.
Can techniques that have been developed in the meantime
be used to give a \ndimensional{4} proof of the Smith Conjecture
(modulo the Poincar\'e Conjecture) that avoids a case-by-case analysis?
\end{question}

\ack
It is a pleasure to be able to thank Claude Weber for his 
support and encouragement since my first trip to Switzerland
in 1982.  My continuing attempts to understand the knot 
theory of complex plane curves have been greatly furthered 
over the years by much work and talk with Michel Boileau, 
Walter Neumann, and Stepan Orevkov, and more recently
Mikami Hirasawa.  I have also had many useful conversations 
with Norbert A'Campo, Thomas Fiedler, Francoise Michel, and 
Andrei Pajitnov.  A number of graduate students (several of them 
advised by mathematicians already mentioned), having written theses 
contributing to the knot theory of complex plane curves, have been 
kind enough to explain their work to me; these include Laurence 
Fourrier, Laure Eichenberger, Arnaud Bodin, and others.
Many other workers in this and related fields have helped
me keep reasonably in touch with developments; I thank them all.

Completion of this survey was partially supported
by the Fonds National Suisse at the University of
Geneva, and by a National Science Foundation 
Interdisciplinary Grant in the Mathematical Sciences,
DMS-0308894.

\appendix
%
%

\section{And now a few words from our 
inspirations}\label{motivations}
\hypertarget{appendix}{}

\subsection{Hilbert's 16th Problem \textup{(1891, 1900)}}
\label{Hilbert's sixteenth problem}
\hypertarget{Hilbert's 16th Problem}{}

\begin{quotation}
{\selectlanguage{german}
\centerline{$16.$ Problem der Topologie algebraischer
Curven und Fl\"achen.}}

Die Maximalzahl der geschlossenen und getrennt liegenden
Z\"uge, welche eine ebene algebraische Curve 
\mbox{$n${\thinspace}ter} Ordnung
haben kann, ist von Harnack \dots~bestimmt worden; es
ensteht die weitere Frage nach der gegenseitigen Lage 
der Curvenz\"uge in der Ebene. \dots~ \emph{Eine 
gr\"undliche Untersuchung der gegenseitigen Lage
bei der Maximalzahl von getrennten Z\"ugen scheint
mir ebenso sehr von Interesse zu sein \dots.}
\end{quotation}
{\hfill \citep[p. 283]{Hilbert1900a}}

\begin{quotation}
\centerline{$16$. Problem of the topology of algebraic 
curves and surfaces}

The maximum number of closed and separate branches which 
a plane algebraic curve of the $n$th order can have has 
been determined by Harnack\dots .  There arises the further 
question as to the relative position of the branches in the 
plane. \dots~ \emph{A thorough investigation of the relative 
position of the separate branches when their number is the 
maximum seems to me to be of very great interest \dots .}
{\hfill \citep[p. 23]{Hilbert1900b}}
\end{quotation}

\subsection{Zariski's Conjecture \textup{(1929)}}
\label{Zariski's conjecture}
\hypertarget{Zariski's Conjecture}{}

\begin{quotation}
\textbf{Theorem 7.}
The fundamental group of an irreducible curve $f$ of order $n$,
possessing ordinary double points only, is cyclic of order $n$.
\noindent\dots~ 
To prove the theorem 7, we prove first the following lemma:

\textbf{Lemma.} \textit{The fundamental group of $n$ lines in 
arbitrary position is abelian.}

\noindent\dots~ 
Now let $f$ be an irreducible curve of order $n$ with
ordinary double points only.  The continuous system $\{f\}$
\dots~ contains in particular curves which degenerate into 
$n$ arbitrary lines.\textsuperscript{*}\dots

\noindent\rule[.5ex]{6em}{.33pt}

{\footnotesize{\textsuperscript{*} This follows 
from a noted principle, announced
by F.~Enriques in 1904, ``Sulla propriet\'a caratteristica delle
superficie algebriche irregolari,'' \textit{Rendiconto delle 
sessioni della R.\ Accademia delle Science dell'Instituto di 
Bologna}, nuova serie, Vol.\ 9 (1904--1905), pp.~5--13, and
completed later by F.~Severi, \textit{Vorlesungen \"uber 
algebraische Geometrie:} Anhang~F, No.~7.  Leipzig, Teubner,
1921. \dots .}}
\end{quotation}
{\hfill \citep[pp. 316, 317, \& 318]{Zariski1929}}

\subsection{Keller's Jacobian Problem \textup{(1939)}}
\label{Keller's Jacobian problem}
\hypertarget{Keller's Jacobian problem}{}
\begin{quotation}
{\selectlanguage{german}
Es sind nur die Faelle $6.$, $7.$ noch 
unentschieden, also die Frage, ob Polynome 
mit der Funktionaldeterminante $1$ sich stets durch 
Polynome umkehren lassen.  \dots~ 
Mir scheint die Frage eine Untersuchung 
sehr zu lohnen, sie scheint jedoch bereits 
im ebenen Fall sehr schwierig zu sein.} 

\end{quotation}
{\hfill \citep[p. 301]{Keller1939}}

[It is only cases $6.$, $7.$ that remain undecided, 
that is, the question of whether a polynomial mapping
with Jacobian determinant $1$ always has a 
polynomial inverse.  \dots~To me it seems that
this is a very worthwhile question, but it
seems to be very difficult already in the
case of the plane.]

\subsection{Chisini's Statement \textup{(1954)}}
\label{Chisini's statement}
\hypertarget{Chisini's statement}{}

\begin{quotation}
{\selectlanguage{italian}
\'E naturale che la esperienza abbia indotto i ricercatori 
ad ammettere (si pure nella sola fase di ricerca) che le 
trecce dedotte dalla treccia canonica con modifica dell'ordine 
dei tratti e con fusione di questi (curve che la treccia 
stessa mostra se irriducibili o non) corrispondano a curve 
algebriche esistenti, cio\`e siano trecce davvero algebriche.
La dimostrazione del fatto indicato \`e lo scopo di questo 
lavoro\dots~.}
\end{quotation}
{\hfill \citep[p. 144]{Chisini1954}}

[Naturally, experience has led researchers
to hypothesize (if only at this stage of research) that
the braids derived from the canonical braid by modification
of the order of the strings and by the fusion thereof (curves
which the braid itself shows to be irreducible or not)
correspond to existing algebraic curves, that is, they
are truly algebraic braids.  The proof of the indicated fact
is the purpose of this work\dots~.]

\subsection{Milnor's Question \textup{(1968)}} 
\label{Milnor's question}
\hypertarget{Milnor's question}{}

\begin{quotation}
Now let us look at the algebraic geometry of our singular point.
To any singular point ${\fontfamily{cmss}\fontseries{b}z}$ of a 
curve $V\sub\C^2$ there is associated an integer 
$\delta_{\fontfamily{cmss}\fontseries{b}z}>0$ which
intuitively measures the number of double points of $V$
concentrated at ${\fontfamily{cmss}\fontseries{b}z}$. \dots

\textbf{Remark 10.9.}  It would be nice to have a better
topological interpretation of the integer $\delta$.
It can be shown that the link $K=V\cap S_\epsilon$ bounds
a collection of $r$ smooth $2$-cells in the disk 
$D_\epsilon$ having no singularities other than $\d$
ordinary double points.  \textit{Question:} Is $\delta$
perhaps equal to the ``\"Uber\-schneidungs\-zahl'' of $K$:
the smallest number of times which $K$ must be allowed
to cross itself during a smooth deformation so as to
transform $K$ into a collection of $r$ unlinked and
unknotted circles? \dots~
\end{quotation}
{\hfill\citet[pp. 85 \& 92]{Milnor1968}}

\subsection{The Thom Conjecture \textup{(c.~1977)}}
\label{Thom conjecture}
\hypertarget{Thom Conjecture}{}

\begin{quotation}
\textbf{Problem 4.36.} (A) \textit{Conjecture}: The minimal genus 
of a smooth imbedded surface in ${CP}^2$ representing 
$n\in H_2({CP}^2;Z)=Z$ is $(n-1)(n-2)/2$.
\end{quotation}
{\hfill\citet[p. 309]{Kirby1978}}

According to \citet{Kirby1997}, ``\dots~when this problem was 
written in 1977, Thom's name was not generally associated 
with the conjecture; however the association developed later 
although it is not clear whether Thom ever made the conjecture.''
Indeed, \citet{Sullivan1973} states the conjecture with no name
attached.  On the other hand, \citet{Massey1973} 
calls it both ``the Thom--Atiyah conjecture'' and 
``the Atiyah--Thom conjecture'', a joint attribution which 
seems to be otherwise unattested in the literature.
Both \citet{Thom1982} and \citet{Massey1992} have disclaimed
being the first to state this conjecture.


\newpage\phantomsection{}
\listoffigures

\def\listfigurename{\protect\hypertarget{list of figures}%
{List of Figures}}
\addcontentsline{toc}{section}%
{\protect\hyperlink{list of figures}{List of Figures}}

\newpage\phantomsection{}
\makeatletter
\providecommand\MR[1]{\relax\ifhmode\unskip\spacefactor3000 \space\fi
  \href{http://www.ams.org/mathscinet-getitem?mr=#1}{MR #1}
}
\makeatother
\def\refname{\protect\hypertarget{references}%
{References\protect\hfill}
\protect\linebreak[4]
\protect\linebreak[4]
{\small\textmd{Items tagged with 
\textcircled{\scriptsize\textsf A}
\resp{\textcircled{\scriptsize\textsf G}}
are freely and permanently available in electronic form at 
\url{http://arxiv.org/math} \resp{\url{http://gdz.sub.uni-goettingen.de/}}; 
in the on-line version of 
this survey, each such tag is hyperlinked to the relevant
electronic resource.}}}
\addcontentsline{toc}{section}%
{\protect\hyperlink{references}{References}}
{\reversemarginpar

}

\def\indexname{\protect\hypertarget{index of definitions}%
{{Index of Definitions and Notations\protect\hfill}
\protect\linebreak[4]
\protect\linebreak[4]
{\small\textmd{Page references are to first, or defining, 
instances of terms and symbols.}}}}

\newpage\phantomsection{}
\addcontentsline{toc}{section}%
{\protect\hyperlink{index of definitions}{Index of Definitions and Notations}}
\cleardoublepage
\input{elsehand.ind}

\end{document}